\documentclass[11pt,
]
{article}
\usepackage{amssymb,amstext,amsthm,latexsym}
\usepackage{amsmath}
\usepackage{amsfonts}

\usepackage{mathrsfs}
\usepackage{mathbbol}

\usepackage{mparhack}

\input{q.sty}

\begin{document}

\title{On effective \bou ness and \sik ness}

\author{Vladimir Kanovei\thanks
{%Partially supported by RFBR. 
Contact author, 
{\tt kanovei@mccme.ru}}
\vyk{
\and
Vassily Lyubetsky\thanks
{%Partially supported by RFBR.  
{\tt lyubetsk@ippi.ru}}
}
}

%\institute{Institute for the information transmission problems, 
%Moscow 101447, Bolshoi Karetnyi, 19, Russia.}

\date{\today}
\maketitle

%\renek{\abstractname}{Abstract}

\begin{abstract}
Different generalizations of a known theorem by Kechris, 
saying that any $\is11$ set $A$ of the Baire space 
either is effectively sigma-bounded 
(that is, covered by a countable union of compact $\id11$ sets), 
or it contains a superperfect subset, are obtained, 
in particular, 
1) with  
covering by compact sets and equivalence classes of 
a given finite collection of $\id11$ \er s, 
2) generalizations to $\is12$ sets, 
3) generalizations true in the Solovay model. 

A generalization to $\is11$ sets $A$, 
of a theorem by Louveau, 
saying that any $\id11$ set $A$ of the Baire space 
either is effectively \sik\  
(that is, is equal to a countable union of compact 
$\id11$ sets), 
or it contains a relatively closed superperfect subset, 
is obtained as well.
\vspace{-3\baselineskip}
\def\contentsname{}
\tableofcontents
\end{abstract}

\newpage
%{\small\tableofcontents}

\punk*{Introduction}
\addcontentsline{toc}{subsection}{\hspace{\parindent}
\protect\phantom{1}Introduction}

%$\bc$

Effective descriptive set theory appeared in the 1950s
as a useful technique of simplification and clarification of 
constructions of classical descriptive set theory 
(see \eg\ \cite{add}, \cite{Shon}, or \cite{umnKL}). 
Yet it has become clear that development of effective 
descriptive set theory also leads to results having no direct 
analogies in classical descriptive set theory. 
As an example we recall the following well-known 
\rit{basis theorem}:
any countable $\id11$ set $A$ of the Baire space 
$\bn=\dN^\dN$ consists of $\id11$ points.
Its remote predecessor in classical descriptive set theory 
is the Luzin -- Novikov theorem on splitting of 
Borel sets with countable cross-sections into countable 
unions of uniform Borel sets.

In this paper, we focus on effectivity aspects 
of the properties of \sik ness and \bou ness of pointsets.
Our starting point will be a pair of classical dichotomy
theorems on
%\sik ness and \bou ness of
pointsets, together with their effective versions
obtained in the end of 1970s.

The first of them deals with the property of \bou ness. 
% instead of \sik ness.
Recall that a pointset is \rit{\bou} iff it is a subset 
of a \sik\ set.\snos
{For subsets of the Baire space $\bn=\om^\om$, the
property of \bou ness is equivalent to being bounded in
$\bn$ with the \rit{eventual domination} order, while 
the compactness is equivalent to being bounded in
$\bn$ with the \rit{termwise domination} order.}
Saint~Raymond~\cite{sr} proved that if $X$ is a $\fs11$ 
set then one and only one of 
the following two (obviously incompatible) 
conditions holds$:$
\ben
\Renu
\itla{san1}
the set\/ $X$ is\/ \bou$;$

\itla{san2}
there is a \sps\ $Y\sq X$.
\een
Recall that a \rit{superperfect} set is a closed set
homeomorphic to $\bn$. 

An effective version of this result
(Theorem~\ref{mt'} below),
by Kechris~\cite{K}, 
says that if $X$ is a $\is11$ set then
condition \ref{san1} can be strengthened to a
\ddd{\id11}effective \bou ness
(so that a given set $X$ is covered by a $\id11$
sequence of compact sets).
The proof in \cite{K} uses the determinacy-style technique.
A different proof of this result, based rather on methods 
of effective descriptive set theory, 
will be presented in Section~\ref{d2}, 
in particular, as a foundation for
%proofs  of
a more general dichotomy theorem  
in Section~\ref{gen+2}. 

The other background result, an immediate concequence of 
a theorem by Hurewicz~\cite{hur}, 
deals with the property of \sik ness instead of 
\bou ness.
It says that if $X$ is a $\fs11$ set
then one and only one of 
the following two (clearly incompatible as well)
conditions \ref{hu1}, \ref{hu2} holds$:$
\ben
\Renu
\itla{hu1}
the set\/ $X$ is\/ \sik$;$

\itla{hu2}
there is a set\/ $Y\sq X$ homeomorphic to\/ $\bn$ 
and relatively closed in\/ $X$.
\een
An effective version of this theorem
(Theorem~\ref{mt} below),
essentially by Louveau~\cite{L}
(see also 4F.18 in \cite{mos} 
which the author of \cite{mos} credits to Louveau),
shows that if $X$ is a lightface $\id11$ set then
condition \ref{hu1} can be strengthened to a
\ddd{\id11}effective \sik ness
(so that a given set is equal to the union of a $\id11$
sequence of compact sets).
We present here a somewhat different proof  
of this result in Section~\ref{gahaSK}, 
in particular, as a base for the proof of 
a similar but more complicated dichotomy theorem
on\/ $\is11$ sets in Section~\ref{mtS}. 
%(Theorem~\ref{pt}).

\vyk{The other classical theorem, by Saint~Raymond~\cite{sr},
deals with the property of \bou ness instead of \sik ness.
Recall that a pointset is \rit{\bou} iff it is a subset 
of a \sik\ set.
For subsets of the Baire space $\bn=\om^\om$, the
property of \bou ness is equivalent to being bounded in
$\bn$ with the \rit{eventual domination} order.
The theorem says if $X$ is a $\fs11$ set then one
and only one of 
the following two (incompatible)
conditions holds$:$
\ben
\Renu
\itla{san1}
the set\/ $X$ is\/ \bou$;$

\itla{san2}
there is a \sps\ $Y\sq X$.
\een
Recall that a \rit{superperfect} set is a closed set
homeomorphic to $\bn$. 

An effective version of this result
(Theorem~\ref{mt'} below),
by Kechris~\cite{K}, 
says that if $X$ is a $\is11$ set then
condition \ref{san1} can be strengthened to a
\ddd{\id11}effective \bou ness
(so that a given set $X$ is covered by a $\id11$
sequence of compact sets).
The proof in \cite{K} uses the determinacy-style technique.
A somewhat different proof of this result
will be presented in Section~\ref{d2}, 
in particular, as a base for
%proofs  of
a more general dichotomy theorem  
in Section~\ref{gen+2}. 
}

Some well-known classical results related to Theorems 
\ref{mt'} and \ref{mt} are discussed in Section~\ref{zz}.
We outline several counterexamples with sets   
more complicated than $\is11$ in Section~\ref{abo}.

Sections \ref{gen+1}, \ref{dig}, \ref{gen+2} contain  
a generalization of Theorem~\ref{mt'} 
(Theorem~\ref{fm}) 
which replaces \bou\ sets by \nsm\ sets, where
$\rF_1,\dots,\rF_{n}$ are given $\id11$ \er s and being
\nsm\ means being covered by the union of a \bou\ set
and countably many equivalence classes of 
$\rF_1,\dots,\rF_{n}$.
Accordingly the condition of existence of a \sps\ 
strengthens by the
requirement that the \sps\ is \pift i\ for $i=1,\dots,n$.
Section~\ref{gen+1} develops a necessary technique while 
the proof of the generalized dichotomy is presented in 
Section~\ref{gen+2}.
In the classical form, the case of a single \er\ $\rF$ 
in this dichotomy was earlier obtained by Zapletal, 
see \cite{ksz}. 

In parallel to this, we prove in Section~\ref{dig} that 
a \bou\  set and a countable union of equivalence classes 
as above can be defined so that they depend only on a 
given set $X$ 
(and the collection of \er s $\rF_j$), 
but are independent of the choice of a parameter $p$ such 
that $X$ is $\is11(p)$ and the relations are $\id11(p)$.  

\vyk{
As usual, the theorems remain true in the
relativized form, \ie\ when classes $\id11$ and $\is11$ are  
replaced by $\id11(p)$ and $\is11(p)$, where $p\in\bn$.  
%is a fixed parameter, with essentially the same proofs.

Conditions \ref{mt'1} and \ref{mt'2} of Theorem~\ref{mt'} 
are incompatible. 
Indeed if $Y$ is a set as is \ref{mt'2} then $Y$ is not 
\bou. 
Conditions \ref{mt1} and \ref{mt2} of Theorem~\ref{mt} are 
incompatible either, 
since $A$ is \ddd\fsg compact provided \ref{mt1} 
holds, so that any relatively closed subset of $A$ is 
\ddd\fsg compact itself, while 
the space $\bn$ is not \ddd\fsg compact, of course.
}

\vyk{
\bcor
\lam{C}
Suppose that\/ $A\sq\bn$. 
Then
\ben
\renu
\itla{C1}
if\/ $A$ is a \bou\/ $\is11$ set then it is\/ 
\ddd{\id11}effectively \bou\ in the sense 
of\/ \ref{mt'1} of Theorem~\ref{mt'}$;$

\itla{C2}
if\/ $A$ is a\/ \sik\/ $\id11$ set 
then it is\/ \ddd{\id11}effectively \sik\ in the sense 
of\/ \ref{mt1} of Theorem~\ref{mt}.\qed
\een
\ecor  

Theorem~\ref{mt} is a direct corollary of some well-known 
results in this field, especially of 4F.18 in \cite{mos} 
which the author of \cite{mos} credits to Louveau.\snos 
{Indeed suppose that a $\id11$ set $A\sq\bn$ does not 
satisfy \ref{mt2} of Theorem~\ref{mt}.
Then $A$ is \ddd\fsg compact by a theorem of Hurewicz 
(see Theorem~\ref{hur}). 
In this case, it is clear from 4F.18 that $A$ 
is equal to the union of all compact $\id11$ sets $A'\sq A$. 
Yet it follows from 4F.14 in \cite{mos} 
that if $X$ is a compact $\id11$ subset of $\bn$ then there 
is a compact $\id11$ tree $T\sq\nse$ such that $X=\bod T$. 
With a little further care this leads to 
condition \ref{mt1} of Theorem~\ref{mt}.} 
Nevertheless we present here somewhat different proofs 
of both theorems, 
in particular, as a base for the proofs of 
\rit{similar 
but more complicated dichotomy theorems on\/ $\is11$ sets\/} 
(\ref{pt} and \ref{nt}).
}

In the remaining parts of the paper, 
we prove, in Sections \ref{uu1}, \ref{uu2}, \ref{BC}, 
a generalization, along the same lines, 
of another Kechris' result of \cite{K}, 
related to $\is12$ sets, which by necessity involves 
uncountable unions of equivalence classes and \bou\ sets 
as well as coding by uncountable constructible Borel codes. 
In the course of the proof of this generalized theorem 
(Theorem~\ref{tk}), it will be shown 
(Theorem~\ref{resT}) that 
if a countable union of equivalence classes 
of a $\id11$ \er\ is $\id11(\xi)$, where $\xi<\omi$, 
then all classes in this union admit Borel coding by 
constructible (not necessarily countable) 
codes.

In the final Sections \ref{sm1}, \ref{sm2} 
we present generalizations of some of the abovementioned 
theorems to ordinal definable pointsets in the Solovay 
model.
Some questions here remain open.

%The authors thank anonymous referees for valuable remarks 
%and suggestions, including an essential improvement in the 
%proof of Theorem~\ref{tks}. 

We thank Alekos Kechris, Ben Miller, Marcin Sabok, and 
Jindra Zapletal for valuable remarks 
and suggestions.

\punk{Preliminaries}
\las{oo}

We use standard notation $\is11\yd\ip11\yd\id11$ for 
effective classes of points and pointsets in $\bn$, as well 
as $\fs11\yd\fp11\yd\fd11$ for corresponding projective classes. 

Let $\nse$ be the set of all finite strings of natural 
numbers, 
including the empty string $\La$. 
If $u,v\in\nse$ then $\lh u$ is the \rit{length} of $u$, 
and $u\su v$ means that $v$ is a 
\rit{proper extension} of $u$. 
If $u\in \nse$ and $n\in\dN$ then $u\we n$ is the string 
obtained by adding $n$ to $u$ as the rightmost term. 
Let, for $u\in\nse,$ 
$$
\ibn u= \ens{x\in\bn}{u\su x}\quad
\text{(\rit{a Baire interval} in $\bn$)}\,. 
$$
If a set $X\sq\bn$ contains at least two elements then there 
is a longest string $u=\stem X$ such that $X\sq\ibn u$.
We put $\diam X=\frac1{1+\stem X}$ in this case, 
and additionally $\diam X=0$ whenever $X$ has at most one 
element.

A set $T\sq\nse$ is a \rit{tree} if 
$u\in T$ holds whenever $u\we n\in T$ for at least one $n$, 
and a \rit{pruned} tree iff $u\in T$ implies $u\we n\in T$ 
for at least one $n$. 
%(so that $T$ has no \rit{endpoints}, or \ddd\su maximal elements).
% 
Any non-empty tree contains $\La$. 
%A string $s\in T$ is an \rit{endpoint} of a tree $T$ if 
%$s\we n\nin T$ for each $n\in\dN$. 
A string $u\in T$ is a \rit{branching point} of $T$ if 
there are 
%exist numbers 
$k\ne n$ such that $u\we k\in T$ 
and $u\we n\in T$; let $\bran T$ be the set of all 
branching points of $T$. 
The \rit{branching height} $\bh Tu$ of a string $u\in T$ in a 
tree $T$ is equal to the number of strings 
$v\in\bran T\yt v\su u$.
For instance, if $T=\nse$ then $\bh{\nse}u=\lh u$ 
for any string $u$.
%
%A tree $T$ is \rit{perfect} iff for any $u\in T$ there 
%is a string $v\in\bran T$ such that $u\su v$.

A tree $T\sq\nse$ is \rit{compact}, if
it is pruned and has \rit{finite branchings}, that is, 
if $u\in \bran T$ then  
$u\we n\in T$ holds for finitely many $n$.
Then
$$
\bod T=\ens{x\in\bn}{\kaz m\:(x\res m\in T)}, 
$$
the \rit{body} of $T$,
is a compact set.
Conversely, if $X\sq\bn$ is compact then 
$$
\der X=\ens{x\res n}{x\in X\land n\in\dN}
$$
is a compact tree.
Let $\ct$ be the $\id11$ set of all non-empty 
compact trees.

A pruned tree $T\sq\nse$ is \rit{perfect}, if
for each $u\in T$ there is a string 
$v\in\bran T$ with $u\su v$.   
Then $\bod T$ is a perfect set.
%Conversely, if $X\sq\bn$ is a perfect set then 
%$\der X$ is a perfect tree. 
%
A perfect tree $T$ is \rit{superperfect}, if
for each $u\in \bran T$ there are infinitely many 
numbers $n$ such that $u\we n\in T$.   
Then $\bod T$ is a superperfect set.
Conversely, if $X\sq\bn$ is a perfect set then 
$\der X$ is a perfect tree, while for any  
\supt\ set $X\sq\bn$ there is a \supt\ tree $T\sq\der X$. 
% is a (super-)perfect tree. 

If $\pX\yi\pY$ are any sets and $P\sq\pX\ti\pY$ then 
$$
\pr P=\ens{x\in\pX}{\sus y\:(\ang{x,y}\in P)}
\quad\text{and}\quad
\seq Px=\ens{y\in\pY}{\ang{x,y}\in P}
$$ 
are, resp., the \rit{projection} of $P$ to $\pX$, and the 
\rit{cross-section} of $P$ defined by $x\in\pX$.
A set $P\sq\pX\ti\pY$ is \rit{uniform} if every 
cross-section $\seq Px$ ($x\in\pX$) contains at most one 
element.
%If $P\sq Q\sq \pX\ti\pY$, $P$ is uniform, and 
%$\pr P=\pr Q$, then they say that $P$ \rit{uniformizes} $Q$. 

\punk{Some basic facts}
\las{some}

We'll make use of several known results of 
effective descriptive set theory. 
They are listed below, with a few proofs 
(of claims which are not in common use in this area) 
attached to make the text self-contained.
%including the following. 

\bdf
\lam{pros}
A \rit{product space} is any finite product 
of factors $\om\yi\nse\yi\bn\yi\pws\nse$.
A \rit{discrete} product space is a finite product 
of $\om\yi\nse$.
\edf   

\vyk{
\bfa
[\rm $\is11$ Separation]
\lam{21}
If\/ $X,Y\sq\bn$ are disjoint\/ $\is11$ sets then there is a\/ 
$\id11$ set\/ $Z\sq\bn$ such that\/ $X\sq Z$ and\/ 
$Y\cap Z=\pu$.\qed
\efa
}

\bfa
[\rm Kreisel selection, 4B.5 in \cite{mos}]
\lam{22}
If\/ $\pX$ is a discrete product space, 
$P\sq\bn\ti\pX$ is a\/ $\ip11$ set,
%the projection\/ $\pr P$ is a\/ $\ip11$ set either,
and\/ $A\sq\pr P$ is a\/ $\is11$ set,
then there is a\/ $\id11$ map\/ 
$f:\bn\to \pX$ such that\/ $\ang{x,f(x)}\in P$ for all\/ 
$x\in A$.\qed
\efa

%The next theorem provides a convenient effective enumeration  
%of all $\id11$ sets $X\sq\nse$. 

\bfa
[\rm 4D.3 in \cite{mos}]
\lam{BQ}
If\/ $P(x,y,z,\dots)$ is a\/ $\ip11$ relation\/  
on a product space then the following derived relations
are\/ $\ip11$, too$:$
$$
\sus x\in\id11\,P(x,y,z,\dots)
\qand
\sus x\in\id11(y)\,P(x,y,z,\dots)\,.
\eqno\qed
$$ 
\efa

\bfa
[\rm 4D.14 in \cite{mos}]
\lam{dp}
If\/ $\pX$ is a product space then the following 
two sets are\/ $\ip11\::$ 
$D=\ens{x\in\pX}{x\,\text{\rm\ is }\,\id11}$ and 
$$
\ens{\ang{p,x}}
{p\in\bn\land x\in\pX\,\land \,x \text{\rm\ is } \id11(p)}.
\eqno\qed
$$ 
\efa

For instance, 
${x\in D}\eqv \sus y\in\id11\,(x=y)$; 
then apply Fact~\ref{BQ}.

%\vyk{
\bfa
[enumeration of $\id11$, 4D.2 in \cite{mos}]
\lam{23}
Let\/ $\pX$ be a product space.
There exist\/ $\ip11$ sets\/ $E\sq\dN$ and\/ 
$W\sq\dN\ti\pX$, and a\/ $\is11$ set\/ $W'\sq\dN\ti\pX$ 
such that\/ 
\ben
\renu
\itsep
\itla{-enu1}
if\/ $e\in E$ then\/ $\seq We=\seq {W'}e$  
{\rm(where $\seq {W}e=\ens{x\in\pX}{\ang{e,x}\in W}$)}$;$ 

\itla{-enu2}
a set\/ $X\sq\pX$ is\/ $\id11$ iff there is\/ 
$e\in E$ such that\/ $X=\seq We$.\qed
\een
\efa
%}

There is a useful uniform version of Fact~\ref{23}.

\bfa
%[\rm uniform enumeration]
\lam{23+}
Let\/ $\pX$ be a product space.
There exist\/ $\ip11$ sets\/ $\bE\sq\bn\ti\dN$ and\/ 
$\bW\sq\bn\ti\dN\ti\pX$, 
and a\/ $\is11$ set\/ $\bW'\sq\bn\ti\dN\ti\pX$ 
such that\/ 
\ben
\renu
\itsep
\itla{+enu1}
if\/ $\ang{p,e}\in \bE$ then\/ $\sek \bW pe=\sek{\bW'}pe$  
{\rm(where, as above,  
$\sek \bW pe=\ens{x\in\pX}{\ang{p,e,x}\in \bW}$)}$;$ 

\itla{+enu2}
if\/ $p\in\bn$ then
a set\/ $X\sq\pX$ is\/ $\id11(p)$ iff there is a number\/ 
$e\in E$ such that\/ $T=\sek \bW pe=\sek {\bW'}pe$.\qed
\een
\efa

This result implies the following stronger version 
of Fact~\ref{22}.
%, also rather known in effective descriptive set theory.

\bfa
[\rm 4D.6 in \cite{mos}]
\lam{22+}
Suppose that\/ $\pX$ is a product space, 
$Q\sq\bn\ti\pX$ is\/ $\ip11$, 
%the projection\/ $A=\pr Q$ onto\/ $\bn$ is\/ $\ip11$, 
$A\sq\pr Q$ is\/ $\is11$, 
and for each\/ $a\in A$ there is a point\/ $x\in\id11(a)$ 
such that\/ $\ang{a,x}\in Q$. 
Then there is a\/ $\id11$ map\/ 
$f:\bn\to \pX$ such that\/ $\ang{a,f(a)}\in Q$ 
for all\/ $a\in A$. 
\efa
\bpf
Assume that $\pX=\bn$, for the sake of brevity. 
Then any $x\in\pX$ satisfies $x\sq\pY=\om\ti\om$.
Making use of the sets $\bE\sq\bn\ti\om$ and 
$\bW,\bW'\sq\bn\ti\om\ti\pY$ as in Fact~\ref{23+}, 
we let 
$$
P=\ens{\ang{a,e}\in\bE}
{\sek \bW ae\in\bn\land \ang{a,\sek \bW ae}\in Q}.
$$
Easily the set $P$ and its projection $\pr P$ both are
$\ip11$, and $A\sq\pr P$. 
By Fact~\ref{22}, there is a $\id11$ map $f:\bn\to \dN$ 
such that $\ang{a,f(a)}\in P$ for all\/ $a\in A$. 
It remains to define $f(a)=\sek \bW {a,}{f(a)}$ for   
$a\in A$; to prove that $f$ is $\id11$  
use both sets $\bW$ and $\bW'$.
\epf

\bfa
[4F.17 in \cite{mos}]
\lam{25}
Let\/ $\pX\yi\pY$ be product spaces, 
$P\sq\pX\ti\pY$ be a\/ $\id11$ set, and every 
cross-section\/ $\seq Px$ ($x\in\pX$) be at most countable. 
Then\/ 
\ben
\renu
\itla{25.1}
$X=\pr P$ is a\/ $\id11$ set, 

\itla{25.2} 
there is a\/ $\id11$ 
set\/ $Q\sq \om \ti\pX\ti\pY$ such that if\/ $n<\om$ then 
the set\/ $Q_n=\ens{\ang{x,y}}{\ang{n,x,y}\in Q}$ is a 
uniform subset of\/ $P$, $\pr Q_n=X$, and\/ 
$P=\bigcup_nQ_n$, and hence

\itla{25.3} 
$P$ is a countable union of\/ $\id11$ sets each of which 
uniformizes\/ $P$.\qed
\een
\efa

\bfa
[a corollary of \ref{25}]
\lam{25c}
If\/ $X\ne\pu$ is a countable\/ $\id11$ set then there
is a\/ $\id11$ map defined on\/ $\om$ such that\/
$X=\ens{f(n)}{n<\om}$.\qed
\efa

\bfa
[4F.14 in \cite{mos}]
\lam{ks=t}
If\/ $F\sq\bn$ is a closed\/ $\id11$ set and\/  
$X\sq F$ is a compact\/ $\is11$ set then there is a
compact\/ $\id11$ tree\/ $T\sq\nse$ such that\/ 
$X\sq\bod T\sq F$. 
In particular, in the case\/ $X=F$, any compact\/ 
$\id11$ set\/ $X\sq\bn$ has the form\/ $X=\bod T$ for 
some compact\/ $\id11$ tree\/ $T\sq\nse$.
\qed 
\efa
\vyk{
\bpf
We claim that there is a $\id11$ tree $S$ 
(not necessarily pruned) such that 
$F=\bod S$. 
Indeed the $\id11$ set $G=\bn\bez F$ is open. 
Therefore the set
$$
P=\ens{\ang{x,t}\in \bn\ti\nse}{x\in\ibn t\land\ibn t\sq G}
$$
satisfies $\dom P=G$.
Moreover $P$ is $\ip11$.  
Fact~\ref{22} 
%9.3.1(iii)*** 
yields a $\id11$ map $f:G\to\nse$ such that 
$x\in\ibn{f(x)} \sq G$ for each $x\in G$. 
Then the set 
$$
U=\ens{f(x)}{x\in G}=\ens{t\in\nse}{\sus x\in G\:(f(x)=t)}
$$
belongs to $\is11$ and satisfies 
$G=\bigcup_{t\in U}\ibn t$. 
However
$$
V=\ens{t\in\nse}{\ibn t\sq G}=
\ens{t\in\nse}{\kaz x\:(x\in\ibn t\imp x\in G)}
$$
is a $\ip11$ set and $U\sq V$. 
By $\is11$ separation there is a 
$\id11$ set $W$ such that $U\sq W\sq V$ --- and then 
$G=\bigcup_{t\in W}\ibn t$. 
Now put 
$$
S=\ens{s\in\nse}{\kaz t\:(t\in W\imp t\not\sq s)}.
$$ 
Easily $S$ is a $\id11$ tree 
(possibly with endpoints), and $\bod S=F$, as required.

\vyk{
The $\ip11$ set
$$
P=\ens{\ang{s,m}}
{s\in\nse\land m<\om\land 
\kaz x\in X\:(s\su x\imp x(\lh s)<m)}
$$
satisfies $\dom P=\nse$ by the compactness of $X$.
(Note that if $s\nin\der X=\ens{x\res m}{x\in X\land m<\om}$ 
then $\ang{u,0}\in P$.)
}

The set $P$ of all pairs $\ang{s,u}$ such that $s\in\nse$, 
$\pu\ne u\sq\dN$ is finite, and 
$$
\kaz k\in u\:(s\we k\in S)\;\land\;  
\kaz x\in X\:
\big(s\su x\imp \sus k\in u\:
(s\we k\su x)
\big),
$$
is $\ip11$ and $\dom P=\nse$. 
(Let $s\in\nse$, $n=\lh s$.  
As $X$ is compact, the set 
$u=\ens{x(n)}{s\su x\in X}$ is finite, 
thus $\ang{s,u}\in P$.) 
Therefore there is a $\id11$ map $f:\nse\to\pwf\dN$ 
such that $\ang{s,f(s)}\in P$ for each $s\in\nse$. 

The tree $T=\ens{s\in\nse}{\kaz n<\lh s\:(s(n)\in f(s\res n))}$ 
has finite branchings since all values of $f$ are finite sets. 
To show that $T$ pruned,  
let $s\in T$. 
Then $s(n)\in f(s\res n)$ for all $n<\lh s$,  
and $f(s)$ is non-empty. 
Pick any $k\in f(s)$; 
the extended string $s'=s\we k$ belongs to $T$. 
Thus $T$ is a compact $\id11$ tree. 

To prove that $X\sq\bod T$ let $x\in X$. 
As $\La=x\res0$ (the empty string) obviously belongs to $T$, 
it suffices to show that if $m<\om$ and $s=x\res m\in T$ 
then the extended string $t=x\res{(m+1)}=s\we x(m)$ belongs to 
$T$ as well, or, equivalently, that $x(m)\in f(s)=f(x\res m)$ 
--- but this holds because $\ang{s,f(s)}\in P$ by the choice 
of $f$.
Finally to prove that $\bod T\sq F$ note that by construction 
$T\sq S$ and $\bod S=F$.
\epf
}

\vyk{
\bpf
(Added for the sake of the reader's convenience.) 
Then the set
$$
H=\ens{\ang{t,n}}{t\in\nse\land n\in\dN\land
\kaz k\:(t\we k\in T(Y)\imp k \le n)}
$$
belongs to $\ip11$ 
(as the $\is11$ set $T(Y)$ stands to the left of an implication), 
and $\dom H=\nse$. 
Fact~\ref{22} implies that there is a $\id11$ map 
$f:\nse\to\dN$ such that $\ang{t,f(t)}\in H$ for all $t\in\nse$. 
We have $y(n)\le f(y\res n)$ for all $y\in Y$ and $n$
by the definition of $H$, and hence $Y\sq[T']$, 
where $T'$ consists of all strings $t\in\nse$ satisfying 
$t(n)\le f(t\res n)$ for every $n<\lh t$.
Yet $T'$ is a compact tree, and clearly $T'$ is $\id11$ 
by the choice of $f$. 
Therefore $[T']\sq U$, and this is a contradiction.
\epF{Lemma}
}

\vyk{
\bfa
[4F.11 in \cite{mos}]
\lam{DinD}
Any compact\/ $\id11$ set\/ $\pu\ne A\sq\bn$ contains a\/ 
$\id11$ element\/ $x\in A$.\qed
\efa 
}

Facts \ref{22}, \ref{BQ},  
\ref{dp} (the first set), 
\ref{23}, \ref{22+}, \ref{25}, \ref{ks=t} 
%, \ref{DinD}  
%\ref{23+} , 
remain true for relativized lightface classes  
$\is11(p)\yd\ip11(p)\yd\id11(p)$, where $p\in\bn$ is 
an arbitrary fixed parameter.
Therefore Facts \ref{22}, \ref{25} also 
hold with lightface classes replaced by boldface 
projective classes $\fs11\yd\fp11\yd\fd11$.

\punk{The Gandy -- Harrington topology}
\las{gaha1}

The \rit{\gh\ topology} on the Baire space $\bn$ consists of  
all unions of $\is11$ sets $S\sq\bn$.
This topology includes the Polish topology on $\bn$ but 
is not Polish. 
Nevertheless the \gh\ topology satisfies a condition  
typical for Polish spaces.

\bdf
\lam{genb'}
Let $\cF$ be any family of sets, \eg\ sets in a given 
background space $\dX$.
A set $\cD\sq\cF$ is \rit{open dense}   
iff \,
$\kaz F\in\cF\:\sus D\in\cD\:(D\sq F)$, \, and %
$$
\kaz F\in\cF\:\kaz D\in\cD\:(F\sq D\imp F\in\cD)\,.
$$
Sets $\cD$ satisfying only the first requirement  
are called \rit{dense}.
If $\cD\sq\cF$ is dense then the set 
$\cD'=\ens{F\in\cF}{\sus D\in\cD\:(F\sq D)}$ is open dense. 
The notions of \rit{open} and \rit{dense} are related to 
a certain topology which we'll not discuss, but not necessarily 
with the topology of the background space $\dX$. 

A \rit{Polish net} for $\cF$ is any   
collection $\ens{\cD_n}{n\in\dN}$ of open dense sets 
$\cD_n\sq\cF$ such that we have $\bigcap_nF_n\ne\pu$ 
for every sequence of sets 
$F_n\in\cD_n$ satisfying the finite intersection property  
(\ie\ $\bigcap_{k\le n}F_k\ne\pu$ for all $n$). 
\edf

For instance the family of all non-empty closed sets of a 
complete metric space $\pX$ admits a Polish net: 
let $\cD_n$ contain all closed sets of diameter 
$\le n\obr$ in $\pX$.
The next theorem is less elementary. 
This theorem and the following corollary are well-known, see \eg\  
\cite{hkl,h-ban,umnS,k-ams}. 

\bte
\lam{s11p}
The collection\/ $\dP$ of all non-empty\/ $\is11$ sets  
in\/ $\bn$ admits a Polish net.\qed
\ete

\punk{Effective \bou ness dichotomy for $\is11$ sets}
\las{d2}

Here we present a proof of the following theorem by methods
of effective descriptive set theory (including the \gh\
topology).
The original proof in \cite{K} 
was based rather on determinacy ideas.

\bte
[Kechris~\cite{K}, p.~198]
\lam{mt'}
If\/ $A\sq\bn$ is a\/ $\is11$ set then one and only one of 
the following two claims \ref{mt'1}, \ref{mt'2} holds$:$
\ben
\tenu{{\rm(\Roman{enumi})}}
\itla{mt'1}\msur
$A$ is\/ \ddd{\id11}effectively \bou, so that
there is a\/ $\id11$ sequence\/ $\sis{T_n}{n\in\dN}$ of 
compact trees\/ 
$T_n\sq\nse$ such that\/ $A\sq\bigcup_n\bod{T_n}\,;$ 

\itla{mt'2}
there is a superperfect set\/ $Y\sq A$.  
\een
\ete

\bcor
\lam{Cbou}
If\/ $A\sq\bn$ is a\/ \bou\/ $\is11$ set then it is\/ 
\ddd{\id11}effectively \bou\ in the sense 
of\/ \ref{mt'1} of Theorem~\ref{mt'}.\qed
\ecor  

\bpf[theorem]
Recall that $\ct$ is the set of all compact trees 
$\pu\ne T\sq\nse$; $\ct$ is $\id11$, of course.
Let $U$ be the union of all sets of the form $\bod T$, 
where $T\sq\nse$ is a compact tree. 
Formally,
$$
x\in U
\leqv
\sus T\in\id11\:
(T\in\ct\land x\in \bod T)\,,
$$
and hence $U$ is $\ip11$ by Fact~\ref{BQ}, and  
the difference $A'=A\dif U$ is a $\is11$ set.

\ble
\lam{tkm*}
Under the conditions of Theorem~\ref{mt'}, 
if\/ $Y\sq A'$ is a non-empty\/ $\is11$ set then 
its topological closure\/ $\clo Y$ in\/ $\bn$ 
is not compact, \ie, the tree\/  
$\der Y =\ens{y\res n}{y\in Y\land n\in\dN}$ 
has at least one infinite branching.
\ele
\bpf
Suppose otherwise: $\clo Y$ is compact. 
Then by Fact~\ref{ks=t} (with $F=\bn$) 
there is a compact $\id11$ tree $T$ such that 
$\clo Y\sq\bod{T}$.
Therefore $Y\sq\clo Y\sq \bod{T}\sq U$, and this contradicts 
to the assumption $\pu\ne Y\sq A'$.
\epF{Lemma}

\rit{Case 1}: 
$A'=\pu$, that is, $A\sq U$. 
To prove \ref{mt'1} of Theorem~\ref{mt'}, 
note that 
%the set
%
\vyk{
$$
\bay{cclcr}
Q&=&\ens{\ang{x,T}}
{x\in\bn\land T\in \ct\cap\id11\land x\in\bod{T}}\,,
&&\text{and\,}\\[1ex]
Z&=&\ens{x\in\bn}
{\sus T\in\id11\,(T\in\ct\land x\in\bod{T})}&=&\pr Q
\eay
$$
}%
$$
Q=\ens{\ang{x,T}}
{x\in\bn\land T\in \ct\cap\id11\land x\in\bod{T}}
$$
is a $\ip11$ set by Facts~\ref{BQ} and \ref{dp}, and obviously
$U=\pr Q$.
By $\is11$ separation there is a $\id11$ set 
$X$ such that $A\sq X\sq U$. 
Then 
$$
P=\ens{\ang{x,T}\in Q}{x\in X}
$$ 
is still a $\ip11$ set, 
and $\pr P=X$ is a $\id11$ set. 
Therefore by Fact~\ref{22+} there is a $\id11$ function 
$\tau:X\to\ct$ such that $\ang{x,\tau(x)}\in Q$ for all $x\in X$. 

Note that $\tau(x)\in\ct\cap\id11$ and $x\in \bod{\tau(x)}$
for all 
$x\in A$ by the construction. 
Thus the full image $R=\ens{\tau(x)}{x\in A}$ is a $\is11$ 
subset of the $\ip11$ set $\ct\cap\id11$, 
and hence there is a $\id11$ set $D$ such that 
$R\sq D\sq \ct\cap\id11$. 
By Fact~\ref{25}, 
there is a $\id11$ map $\delta:\dN\onto D$.
Now put $T_n=\delta(n)$ for all $n$, 
getting \ref{mt'1} of Theorem~\ref{mt'}.\vom

\rit{Case 2}: $A'=A\dif U\ne\pu$.
To prove that \ref{mt'2} of Theorem~\ref{mt'} holds,
we'll define a system of  
$\is11$ sets $\pu\ne Y_u\sq A'$ satisfying the following
conditions: 
\ben
\itsep
\tenu{(\arabic{enumi})}
\itla{gan1}
if $u\in\nse$  and $i\in\dN$ then 
$Y_{u\we i}\sq Y_u$;

\itla{gan2} 
$\diam{Y_u}\le2^{-\lh u}$; 

\itla{gan3}
if $u\in\nse$ and $k\ne n$ then
$Y_{u\we k}\cap Y_{u\we n}=\pu$, and moreover, sets 
$Y_{u\we k}$ are covered by pairwise disjoint (clopen) 
Baire intervals $J_{u\we k}$;

\itla{gan4} 
$Y_s\in\cD_{\lh u}$, where by Theorem~\ref{s11p} 
$\ens{\cD_n}{n\in\dN}$ is a fixed Polish net for the family  
$\dP$ of all non-empty $\is11$ sets $Y\sq\bn$;

\itla{han5} 
if $u\in\nse$ and $x_k\in Y_{u\we k}$ for all $k\in\dN$ 
then the sequence of points $x_k$ does not have 
convergent subsequences in $\bn$.
\een

If such a construction is accomplished then \ref{gan4} implies 
that $\bigcap_mY_{a\res m}\ne\pu$ for each $a\in\bn$. 
On the other hand by \ref{gan2} every such an intersection 
contains a single point, which we denote by $f(a)$, and the  
map $f:\bn\na Y=\ran f=\ens{f(a)}{a\in\bn}$ is   
a homeomorphism by clear reasons. 

Prove that $Y$ is closed in $\bn$. 
Consider an arbitrary sequence of points $a_n\in\bn$ such 
that the corresponding sequence of points $y_n=f(a_n)\in Y$ 
converges to a point $y\in\bn$; we have to prove that 
$y\in Y$.
If the sequence $\sis{a_n}{n\in\dN}$ contains a subsequence 
of points $b_k=a_{n(k)}$ convergent to some 
$b\in\bn$ then quite obviously the sequence of points 
$z_k=f(b_k)$ (a subsequence of $\sis{y_n}{n\in\dN}$)
converges to $z=f(b)\in Y$, as required. 
Thus suppose that the sequence $\sis{a_n}{n\in\dN}$ 
has no convergent subsequences. 
Then it cannot be covered by a compact set, and it easily 
follows that there is a string $u\in\nse$, an infinite 
set $K\sq\dN$, and for each $k\in K$ --- a number 
$n(k)$ such that $u\we k\su a_{n(k)}$. 
But then $y_{n(k)}\in Y_{u\we k}$ by construction. 
Therefore the subsequence $\sis{y_{n(k)}}{k\in\dN}$ 
diverges by \ref{han5}, which is a contradiction.

Thus $Y$ is closed,
%but condition \ref{han5} is stronger than \ref{gan5}. 
and hence we have \ref{mt'2} of Theorem~\ref{mt'}.  

As for the construction of sets $Y_u$, 
if a $\is11$ set $Y_u\sq A'$ is defined then 
by Lemma~\ref{tkm*} there is a string $t\in T(Y_u)$ such 
that $t\we k\in T(Y_u)$ for all $k$ in an infinite set 
$K_u\sq\dN$. 
This allows us to define a sequence of pairwise different 
points $y_k\in Y_u$ ($k\in\dN$) having no convergent 
subsequences. 
We cover these points by Baire intervals $U_k$ small enough 
for \ref{han5} to be true for the $\is11$ sets  
$Y_{u\we i}=Y_u\cap U_i$, and then shrink these sets if 
necessary to fulfill \ref{gan2} and \ref{gan4}.

\epF{Theorem~\ref{mt'}}

\punk{Effective \sik ness dichotomy for $\id11$ sets}
\las{gahaSK}

Here we present a proof of the following result.

\bte
[essentially Louveau~\cite{L}] 
\lam{mt}
If\/ $A\sq\bn$ is a\/ $\id11$ set then one and only one of 
the next two claims holds$:$
\ben
\tenu{{\rm(\Roman{enumi})}}
\itla{mt1}\msur
$A$ is\/ \ddd{\id11}effectively \sik, so that
there is a\/ $\id11$ sequence\/ $\sis{T_n}{n\in\dN}$ 
of compact trees\/ $T_n\sq\nse$ such that\/ 
$A=\bigcup_n\bod{T_n}\,;$

\itla{mt2}
there is a set\/ $Y\sq A$ homeomorphic to\/ $\bn$ 
and relatively closed in\/ $A$.
\een
\ete

\bcor
\lam{Ckom}
If\/ $A\sq\bn$ is a\/ \sik\/ $\id11$ set 
then it is\/ \ddd{\id11}effectively \sik\ in the sense 
of\/ \ref{mt1} of Theorem~\ref{mt}.\qed
\ecor  

\bpf[theorem]
By Theorem~\ref{mt'}, we can \noo\ assume 
that $A$ is \bou, and hence 
if $F\sq A$ is a closed set then $F$ is \ddd\fsg compact.
Further, the union $U$ of all sets $\bod{T}\sq A$, 
where $T$ is a compact $\id11$ tree, is $\ip11$:
$$
x\in U
\leqv
\sus T\in\id11\:
(T\,\text{ is a compact tree and }\,x\in \bod{T}\sq A)\,,
$$
and the result follows from Fact~\ref{BQ}.
We conclude that $A'=A\bez U$ is $\is11$. 

\ble
\lam{tkm-l}
If\/ $F\sq A'$ is a non-empty\/ $\is11$ set then\/ 
$\clo F\not\sq A$. 
% is not\/ \ddd\fsg compact.
\ele
\bpf
We first prove that if $X\sq A$ is a compact $\is11$ set 
then $A'\cap X=\pu$.
Suppose towards the contrary that $A'\cap X$ is non-empty. 
We are going to find a closed $\id11$ set $F$ satisfying 
$X\sq F\sq A$ --- this would imply $X\sq U$ by 
Fact~\ref{ks=t}, which is a contradiction. 

Since the complementary $\ip11$ set $C=\bn\bez X$ is open, 
the set 
$$
H=\ens{\ang{x,u}}
{u\in\nse\land x\in C\cap\ibn u\land \ibn u\cap X=\pu}
$$ 
is $\ip11$ and $\pr H=C$. 
Thus the $\id11$ set $D=\bn\bez A$ is included 
in $\pr H$. 
By Fact~\ref{22},
there is a $\id11$ map $\nu:D\to\nse$ such that 
$x\in D\imp \ang{x,\nu(x)}\in H$, or equivalently, 
$x\in \ibn{\nu(x)}\sq C$ for all $x\in D$. 
Then the set $\Sg=\ran\nu=\ens{\nu(x)}{x\in D}\sq\nse$ 
is $\is11$ and $D\sq \bigcup_{u\in\Sg}\ibn u\sq C$.

But $\Pi=\ens{u\in\nse}{\ibn u\sq C}$ is a $\ip11$ set
and $\Sg\sq \Pi$.
It follows that there exists a 
$\id11$ set $\Da$ such that $\Sg\sq\Da\sq\Pi$. 
Then still $D\sq \bigcup_{s\in\Da}\ibn s\sq C$,
and hence the closed set 
$F=\bn\bez \bigcup_{u\in\Da}\ibn u$ satisfies $X\sq F\sq A$. 
But $x\in F$ is equivalent to $\kaz u\:(u\in\Da\imp x\nin \ibn u)$, 
thus $F$ is $\id11$, as required.

Now suppose towards the contrary that 
$\pu\ne F\sq A'$ is a $\is11$ set but $\clo F\sq A$. 
By the \noo\ assumption above,  
$\clo F=\bigcup_nF_n$ is \ddd\fsg compact, 
where all $F_n$ are compact. 
There is a Baire interval $\ibn u$ such that the set 
$X=\ibn u\cap \clo F$ is non-empty and $X\sq F_n$ for some $n$. 
Thus $X\sq A$ is a non-empty compact $\is11$ set, hence 
$X\cap A'=\pu$ by the first part of the proof. 
In other words, $\ibn u\cap \clo F\cap A'=\pu$.
It follows that $\ibn u\cap F=\pu$ (because $F\sq A'$), 
contrary to $X=\ibn u\cap \clo F\ne \pu$. 
\epF{Lemma}

We return to the proof of the theorem.\vom

\rit{Case 1}: $A'=\pu$, that is, $A=U$. 
This implies \ref{mt1} 
of Theorem~\ref{mt}, exactly as 
in the proof of Theorem~\ref{mt'} above.\vom
      
\rit{Case 2}: 
$A'=A\bez U\ne\pu$. 
To get a set $Y\sq A'$, \rit{relatively} 
closed in $A$ and homeomorphic to $\bn$, as 
in \ref{mt2} of Theorem~\ref{mt}, we'll define  
a system of non-empty $\is11$ sets $Y_u\sq A'$ satisfying 
conditions \ref{gan1}, \ref{gan2}, \ref{gan3}, \ref{gan4} of 
Section~\ref{d2}, along with the next requirement 
instead of \ref{han5}:
\ben
\itsep
\tenu{$(\arabic{enumi}')$}
\atc\atc\atc\atc
\itla{gan5} 
if $u\in\nse$ then there is a point $y_u\in \clo{Y_u}\bez A$ 
such that any sequence of points 
$x_k\in Y_{u\we k}$ ($k\in\dN$) converges to $y_u$.
\een

If we have defined such a system of sets, then the associated 
map $f:\bn\to A'$ is $1-1$ and is a homeomorphism 
from $\bn$ onto its full image 
$Y=\ran f =\ens{f(a)}{a\in\bn}\sq A'$,
%by the same reasons 
as in the proof of Theorem~\ref{mt'}.  

Let's prove that $Y$ is relatively closed in $A$. 
Consider a sequence of points $a_n\in\bn$ such 
that the corresponding sequence of $y_n=f(a_n)\in Y$ 
converges to a point $y\in\bn$; 
we have to prove that $y\in Y$ or $y\nin A$.
If the sequence $\sis{a_n}{}$ contains a subsequence 
convergent to $b\in\bn$ then, as in the proof of 
Theorem~\ref{mt'}, 
$\sis{y_n}{}$ converges to $f(b)\in Y$. 
If the sequence $\sis{a_n}{}$ 
has no convergent subsequences, 
then there exist a string $u\in\nse$, 
an infinite set $K\sq\dN$, 
and for each $k\in K$ --- a number $n(k)$, 
such that $u\we k\su a_{n(k)}$. 
But then $y_{n(k)}\in Y_{u\we k}$ by construction. 
Therefore the subsequence $\sis{y_{n(k)}}{k\in\dN}$ 
converges to a point $y_u\nin A$ by \ref{gan5}, 
as required.  

Finally on the construction of sets $Y_s$. 

Suppose that a $\is11$ set $\pu\ne Y_u\sq A'$ is defined. 
Then its closure $\clo{Y_u}$ is a $\is11$ set, too, 
therefore $\clo{Y_u}\not\sq A$ by Lemma~\ref{tkm-l}.
There is a sequence of 
pairwise different points $x_n\in Y_u$ which converges to 
a point $y_u\in \clo{Y_u}\bez A$. 
Let $U_n$ be a neighbourhood of $x_n$ (a Baire interval) 
of diameter less than $\frac13$ of the least distance 
from $x_n$ to the points $x_k\yt k\ne n$.
Put $Y_{u\we n}=Y_u\cap U_n$, and shrink the sets 
$Y_{u\we n}$ so that they satisfy \ref{gan2} 
and \ref{gan4}.\vom

\epF{Theorem~\ref{mt}}

\punk{Effective \sik ness dichotomy: 
generalization to $\is11$}
\las{mtS}

There is a difference between Theorem~\ref{mt'} and 
Theorem~\ref{mt}:
the first theorem deals with $\is11$ sets $A$ 
while the other one --- with $\id11$ sets only.
The proof of Theorem~\ref{mt} in Section~\ref{gahaSK} 
does not work in the case when $A$ is a $\is11$ set.
Indeed then $A'$ is a set in $\fs11$ and $\is12$, but,
generally speaking, it cannot be expected to be a $\is11$ set, 
so the rest of the proof does not go through.

As a matter of fact, Theorem~\ref{mt} {\sl per se\/} 
fails for $\is11$ sets $A$, 
as the following counterexample shows. 

\vyk{
The following counterexample is just a slight modification 
of Counterexample~3 in \cite{sti} which we render here 
closely to the text of \cite{sti}. 

\bex
\lam{ex1}
For $x\in\bn$ let $w^x\in2^\dN$ be the characteristic function of 
$$
W^x=\ens{e\in\dN}
{e\,\text{ is a G\"odel number of a wellorder of }\,\dN\,
\text{ recursive in }\,x}.
$$  
Let $\bo\in\bn$ be the constant $0$, and let $y=w^{w^\bo}$. 
The singletons $\ans{w^\bo}$ and $\ans y$ are $\ip11$.
Therefore the set $A= 2^\dN\bez\ans y$ is $\is11$ and is 
an open subset of $2^\dN$, hence, \ddd\fsg compact. 
Suppose towards the contrary that Theorem~\ref{mt} holds for 
$A$. 
Then \ref{mt1} of Theorem~\ref{mt} must be true. 
Let $\sis{T_n}{n\in\dN}$ be a $\id11$ 
sequence of compact trees such that $A=\bigcup_n\bod{T_n}$.
The sequence of trees $T_n$ is then recursive in $w^\bo$.
Therefore $y=w^{w^\bo}$, 
as the only point in $2^\dN$ which does not 
belong to $\bigcup_n\bod{T_n}$, is $\id11$ in $w^\bo$, 
a contradiction.\qed
\eex
}

\bex
\lam{ex1}
Let $\ans{y}$ be a $\ip11$ singleton such that $y\in\dn$ 
is not $\id11$. 
The set $A= 2^\dN\bez\ans y$ is then $\is11$ and  
an open subset of $2^\dN$, hence, \ddd\fsg compact. 
Suppose towards the contrary that Theorem~\ref{mt} holds for 
$A$. 
Then \ref{mt1} of Theorem~\ref{mt} must be true. 
Let $\sis{T_n}{n\in\dN}$ be a $\id11$ 
sequence of compact trees such that $A=\bigcup_n\bod{T_n}$.
Therefore $y$ is $\id11$, 
as the only point in $2^\dN$ which does not 
belong to $\bigcup_n\bod{T_n}$, 
a contradiction.\qed
\eex

\vyk{
\bqe
\lam{t2?}
Does Theorem~\ref{mt} hold for all $\is11$ sets $A$ in a 
modified form, where condition \ref{mt1} is weakened to 
its ``therefore'' part\,?

Note that the set $A= 2^\dN\bez\ans y$ of Example~\ref{ex1} 
satisfies the ``therefore'' part of \ref{mt1} of 
Theorem~\ref{mt}.
\qed 
\eqe 
}
 
Our best result in the direction of Theorem~\ref{mt} for 
$\is11$ sets with still 
some effectivity in \ref{mt1} is the following theorem:

\bte
\lam{pt}
%Пусть\/ $p\in\bn$. 
If\/ $A\sq\bn$ is a\/ $\is11$ set then one and only one of 
the following two claims holds$:$
\ben
\tenu{{\rm(\Roman{enumi})}}
\itla{pt1} 
$A$ is\/ \ddd{\id13}effectively \sik, so that 
there exists a\/
\vyk{
\/$:$ 
a countable\/ $\id13$ ordinal\/ $\la$ and an\/
}% 
$\id13$ 
sequence\/ $\sis{T^n}{n<\om}$ of compact\/ 
$\id13$ trees\/ $T^n\sq\nse$ 
such that\/ $A=\bigcup_{n<\om}\bod{T^n}\,;$ 

\itla{pt2}
there is a set\/ $Y\sq A$ homeomorphic to\/ $\bn$ and relatively 
closed in\/ $A$.
\een
\ete

\vyk{
We'll not try to estimate the level and character of the 
effectivity condition in \ref{pt1}, since we don't think that 
our construction gives a result even close to optimal. 
But it will be quite clear from the construction that it 
is absolute for all transitive models 
containing the true $\omi$, and lies 
within the projective hierarchy and probably within $\id13$.
}

\bpf
Given a tree $S\sq\nnse$, define a \rit{derived tree} 
$S'\sq S$ so that
\ben
\fenu
\itla{+1A}\msur
$S'$ consists of all nodes $\ang{u,v}\in S$ 
such that $\clo{\pr{\bod{\tes Suv}}}\not\sq A$, where
$\tes Suv=\ens{\ang{u',v'}\in S}
{(u\su u'\land v\su v')\lor(u'\sq u\land v'\sq v)}$.
\een
Note that $S'$ can contain maximal nodes even if 
$S$ contains no maximal nodes. 
Yet if $\ang{u,v}$ is a maximal node in $S$, 
or generally a note in the well-founded part of $S$ 
(so $\bod{\tes Suv}=\pu$), 
then definitely $\ang{u,v}\nin S'$.

\ble
\lam{ptl}
The set\/ 
$\ens{\ang{S,u,v}}{\ang{u,v}\in S'}$ 
is\/ $\is12$. 

In addition,\/ $S'\sq S$, and if\/ $S\sq T$ then\/ 
$S'\sq T'$.

Moreover, if\/ $\gM$ is a countable transitive model 
of a large enough fragment of\/ $\ZFC$ and\/ 
$S\in\gM$ then\/ $(S')^\gM\sq S'$.
\ele
\bpf
As $A$ is $\is11$, the key condition 
$\clo{\pr{\bod{\tes Suv} }}\not\sq A$ 
is $\is12$.
\epf

Beginning the proof of Theorem~\ref{pt}, 
we \noo\ assume, by Theorem~\ref{mt'},  
that $A$, the given set, 
is \bou, and hence 
if $F\sq A$ is a closed set then $F$ is \ddd\fsg compact.
Let $P\sq\bn\ti\bn$ be a $\ip01$ set such that 
$A=\pr P$. 
Let
$$
S=\ens{\ang{x\res n,y\res n}}
{n\in\dN\land \ang{x,y}\in P} \sq\nse\ti\nse,
$$
so that $P=\bod{S}$.
A decreasing sequence of derived trees 
$\sa S\al,\;\al\in\Ord$, is 
defined by transfinite induction so that 
$\sa S0=S$, if $\la$ is a limit ordinal then naturally 
$\sa S\la=\bigcap_{\al<\la}\sa S\al$, and 
$\sa S{\al+1}=(\sa S{\al})'$ for any $\al$.

Obviously there is a countable ordinal $\la$ such that 
$\sa S {\la+1}=\sa S\la$. 
\vom

\rit{Case 1}: $\sa S\la=\pu$. 
Then, if $x\in A=\pr P$ then by construction there exist an 
ordinal $\al<\la$ and a node $\ang{u,v}\in \sa S\al$ such 
that 
$$
x\in \aog \al uv\sq \clo{\aog \al uv}\sq A\,,
\quad
\text{where}
\quad
\aog \al uv=\pr{\bod{\sog\al uv} }\,,
$$ 
and hence $A$ is a countable union of sets $F\sq A$ of the 
form $\clo{\aog \al uv}$, where $\al<\la$ and 
$\ang{u,v}\in \sa S\al$, closed, therefore 
\ddd\fsg compact by the above. 

Let us show how this leads to \ref{pt1} of the theorem.

It easily follows from Lemma~\ref{ptl} 
that both the ordinal $\la$, and each ordinal $\al<\la$, 
and the sequence 
$\sis{\sa S\al}{\al<\la}$ itself, are $\id13$. 
Therefore there is a $\id13$ sequence 
$\sis{\sa Un}{n<\om}$ of the same trees, that is, 
$$
\ens{\sa S\al}{\al<\la}=\ens{\sa Un}{n<\om}.
$$
Each tree $\sa Un\yt n<\om$, is $\id13$ either, 
as well as all restricted subtrees of the form $\uog nuv$ 
(where $\ang{u,v}\in \sa Un$) 
and their \lap{projections} 
$$
\tog nuv=\ens{u}{\sus v\,(\ang{u,v}\in \uog nuv)}
\sq\nse.
$$
On the other hand, if $\al<\la$ and $\ang{u,v}\in \sa S\al$ 
then we have $\clo{\aog \al uv}=\bod{\tog nuv}$
for some $n=n(\al)$
by construction.

To conclude, if $x\in A$ then there is a 
%pruned 
$\id13$ tree $\tog nuv\sq\nse$ 
such that $x\in \bod{\tog nuv}\sq A$ --- 
and $\bod{\tog nuv}$ is \sik\ in this case.
Then by Theorem~\ref{mt} (relativized version) 
there is a $\id11(\tog nuv)$ sequence 
of compact trees $\tog nuv(k)$ 
such that $\bod{\tog nuv}=\bigcup_k\bod{\tog nuv(k)}$. 
This easily leads to \ref{pt1} of the theorem.\snos
{Class $\id13$ in \ref{pt1} of the theorem looks too weird. 
One may want to improve it to $\id12$ at least. 
This would be the case if the ordinal $\la$ in the 
argument of Case 1 could be shown to be $\id12$. 
Yet by Martin~\cite{mar} closure ordinals of 
inductive constructions of this sort may exceed the domain of 
$\id12$ ordinals.} 
\vom

\rit{Case 2}: $S^{(\la)}\ne\pu$, and then $S^{(\la)}\sq S$ 
is a pruned tree. 

\ble
\lam{tt''}
If\/ $\ang{u,v}\in \sa S\la$, $u'\in\nse$, $u\su u'$, and\/ 
$\aog \la uv\cap\ibn{u'}\ne\pu$ then there is a string\/ 
$v'\in\nse$ such that\/ $v\su v'$ and\/ 
$\ang{u',v'}\in \sa S\la$.\qed
\ele

We'll define a pair $\ang{u(t),v(t)}\in \sa S\la$ 
for each $t\in\nse$, such that  
%$\ang{t,v(t)}\in S'$ and 
%
\ben
\itsep
\tenu{(\arabic{enumi})}
\itla{px1}
if $t\in\nse$ then $t\sq u(t)$;

\itla{px2}
if $s,t\in\nse$ and $s\sq t$ then 
$u(s)\sq u(t)$ and $v(s)\sq v(t)$;

\itla{px3}
if $t\in\nse$ and $k\ne n$ then 
$u(t\we k)$ and $u(t\we n)$ are \ddd\sq incomparable;

\itla{px4}
if $s\in\nse$ then there exists a point 
$y_s\in \clo{\aog \la {u(s)}{v(s)}}\bez A $ 
%_{u(s),v(s)}$ 
such that
any sequence of points 
$x_k\in \aog \la {u(s\we k)}{v(s\we k)}$ 
converges to $y_s$.
\een

Suppose that such a system of sets is defined. 
Then the associated map 
$f(a)=\bigcup_n u(a\res n):\bn\to A$ is $1-1$ and 
is a homeomorphism 
from $\bn$ onto its full image 
$Y=\ran f =\ens{f(a)}{a\in\bn}\sq A$.  

Let's prove that $Y$ is relatively closed in $A$. 
Consider a sequence of points $a_n\in\bn$ such that 
the corresponding sequence of points $y_n=f(a_n)\in Y$ 
converges to a point $y\in\bn$; 
we have to prove that $y\in Y$ or $y\nin A$.
If the sequence $\sis{a_n}{}$ contains a subsequence 
convergent to $b\in\bn$ then 
$\sis{y_n}{}$ converges to $f(b)\in Y$. 
So suppose that the sequence $\sis{a_n}{}$ 
has no convergent subsequences. 
Then there exist a string $s\in\nse$, 
an infinite set $K\sq\dN$, 
and for each $k\in K$ --- a number $n(k)$, 
such that $s\we k\su a_{n(k)}$. 
Then $y_{n(k)}\in \aog\la{u(s\we k)}{v(s\we k)}$ 
by construction. 
Therefore the subsequence $\sis{y_{n(k)}}{k\in\dN}$ 
converges to a point $y_s\nin A$ by \ref{px4}, 
as required.
 
Finally on the construction of sets $Y_s$. 

Suppose that a pair $\ang{u(t),v(t)}\in \sa S\la$ is defined. 
Then $\clo{\aog\la{u(t)}{v(t)}}\not\sq A$ by the choice of $\la$.
There is a sequence of pairwise different points 
$x_n\in\aog\la{u(t)}{v(t)}$ which converges to 
a point $y_s\in \clo{\aog\la{u(t)}{v(t)}}\bez A$. 
We can associate a string $u_n\in\nse$ with each $x_n$ 
such that $u(t)\su u_n\su x_n$, the strings $u_n$ are 
pairwise \ddd\sq incompatible, and $\lh{u_n}\to\infty$.
Then, by Lemma~\ref{tt''}, for each $n$ there is a 
matching string $v_n$ such that $v(t)\su v_n$ and  
$\ang{u_n,v_n}\in \sa S\la$. 
Put $u(t\we n)=u_n$ and $v(t\we n)=v_n$ for all $n$.
\vom

\epF{Theorem~\ref{pt}}

\punk{Related classical results}
\las{zz}

The ``effective'' results presented above can be compared 
with some known theorems of classical descriptive set theory, 
including the following two.

\bte
[\rm Saint~Raymond~\cite{sr} or 21.23 in \cite{dst}]
\lam{hur2}
If\/ $A$ is a\/ $\fs11$ set in a Polish space then either\/
$A$ is \bou\ or 
there is a superperfect set\/ $P\sq A$.\qed
\ete

\bte
[\rm Hurewicz~\cite{hur}]
\lam{hur}
If\/ $A$ is a\/ $\fs11$ set in a Polish space then either\/
$A$ is \sik\ or 
there is a subset\/ $Y\sq A$ homeomorphic to the 
Baire space\/ $\bn$ 
and relatively closed in\/ $A$.\qed
\ete

Arguments in \cite{dst} show that it's sufficient to prove 
each of these theorems in the case $\pX=\bn$; then the 
results can be generalized to an arbitrary Polish space $\pX$ 
by purely topological methods. 
In the case $\pX=\bn$, Theorem~\ref{hur2} immediately follows  
from our Theorem~\ref{mt'} 
(in relativized form, \ie, 
for classes $\is11(p)$, where $p\in\bn$ is arbitrary),  
while Theorem~\ref{hur} follows from  
Theorem~\ref{pt} (relativized). 
On the other hand, Theorem~\ref{hur} also follows from  
Theorem~\ref{mt} (relativized) 
for sets $A$ in $\fd11$ (that is, Borel sets).\snos
{See \cite{dst,ls} for another modern approach to those 
classical theorems, based mainly on infinite games rather 
than methods of effective descriptive theory.} 

Theorem~\ref{mt} implies yet another theorem, 
which combines several classical results of descriptive 
set theory by Arsenin, Kunugui, Saint~Raymond, 
She\-gol\-kov, see references in \cite{dst} or 
in \cite[\S\,4]{umnL}.

\bte
[\rm compare with Fact~\ref{25}]
\lam{tks}
Suppose that\/ $\pX,\pY$ are Polish spaces, $P\sq\pX\ti\pY$  
is a\/ $\fd11$ set, and all cross-sections\/ 
$\seq Px=\ens{y}{\ang{x,y}\in P}$ {\rm($x\in\pX$)} 
are\/ \ddd\fsg compact. 
Then\/ 
\ben
\renu
\itla{tks2}
the projection\/ $\pr P$ is a\/ $\fd11$ set$;$

\itla{tks3} 
$P$ is a countable union of\/ $\fd11$ sets with compact 
cross-sections$;$

\itla{tks1}
$P$ can be uniformized by a\/ $\fd11$ set.
\een
\ete
\bpf[a sketch for the case $\pX=\pY=\bn$]
\ref{tks2}
Assume, for the sake of simplicity, that $P\sq\bn\ti\bn$ is   
a $\id11$ set. 
The set 
$$
H=\ens{\ang{x,T}} 
{x\in\bn\land T\in\ct\land T\in\id11(x)\land
\bod{T}\sq \seq Px}
$$
is $\ip11$ by Fact~\ref{dp}. 
It follows from Theorem~\ref{mt} that if $\ang{x,y}\in P$ 
then there is a tree $T$ such that $\ang{x,T}\in H$ 
and $y\in\bod{T}$. 
Therefore the $\ip11$ set
$$
E=\ens{\ang{x,y,T}}
{\ang{x,y}\in P\land \ang{x,T}\in H\land y\in \bod{T}}
\sq\bn\ti\bn\ti 2^{(\nse)}
$$
satisfies $\pr_{xy}E=P$, 
that is, if $\ang{x,y}\in P$ then there is a tree $T$ such that 
$\ang{x,y,T}\in E$. 
There is a uniform $\ip11$ set $U\sq E$ which 
\rit{uniformizes} $E$, 
\ie, if $\ang{x,y}\in P$ then there is a unique $T$ such that 
$\ang{x,y,T}\in U$. 
Yet $U$ is $\is11$ as well by Fact~\ref{BQ}, 
since $\ang{x,y,T}\in U$ is 
equivalent to:
$$
\ang{x,y}\in P\land y\in \bod{T}\land \kaz T'\in\id11(x)\:
(\ang{x,y,T'}\in U\imp T=T')\,.
$$ 
%and quantifiers of the form $\kaz x\in\id11(y)$ are known to 
%preserve class $\is11$. 
Thus the $\is11$ set 
$F=\ens{\ang{x,T}}{\sus y\:(\ang{x,y,T}\in U)}$ 
is a subset of the $\ip11$ set $H$. 
By $\is11$ separation, there is a $\id11$ set $V$ such 
that $F\sq V\sq H$.
Then
$$
\ang{x,y}\in P\leqv \sus T\:
(\ang{x,T}\in V\land y\in\bod{T}) 
$$
by definition.
Finally all cross-sections of $V$ are at most countable: 
indeed if $\ang{x,T}\in V$ then $T\in\id11(x)$ 
(since $V\sq H$). 
Note that $\pr P=\pr V$, and hence the projection $D=\pr P$ 
is $\id11$ (hence Borel) by Fact~\ref{25}.\vom 

\ref{tks3} 
It follows from Fact~\ref{25} that $V$ is equal to a  
union $V=\bigcup_nV_n$ of uniform $\id11$ sets $V_n$, 
and then each projection $D_n=\pr V_n\sq D$ is $\id11$. 
Each $V_n$ is basically the graph of a $\id11$ map  
$\tau_n:D_n\to \ct$, and 
$\seq Px=\bigcup_{x\in D_n}\bod{\tau_n(x)}$. 
If $n\in\dN$ then we put  
$$
P_n=\ens{\ang{x,y}}
{x\in D_n\land  y\in\bod{\tau_n(x)})}\,. 
$$
Then $P=\bigcup_nP_n$ by the above, each set $P_n$ has only  
compact cross-sections, and each $P_n$ is a $\id11$ set,  
since the sets $D_n$ and maps $\tau_n$ belong to $\id11$.\vom

\ref{tks1} 
Still by Fact~\ref{25}, the set $V$ can be uniformized by  
a uniform $\id11$ set, that is, there exists a  
$\id11$ map $\tau:D\to \ct$ 
such that $\ang{x,\tau(x)}\in V$ for all $x\in D$.
To uniformize the original set $P$, let $Q$ consist of all pairs 
$\ang{x,y}\in P$ such that $y$ is the lexicographically leftmost 
point in the compact set $\bod{\tau(x)}$. 
Clearly $Q$ uniformizes $P$.
To check that $Q$ is $\id11$, note that 
``$y$ is a the lexicographically leftmost point in $\bod{T}$'' 
is an arithmetic relation in the assumption that $T\in\ct$.
\epf

Similar arguments, this time based on Theorem~\ref{mt'}, also 
lead to an alternative proof of the following known result. 

\bte
[Burgess, Hillard, 35.43 in \cite{dst}] 
\lam{buhi}
If\/ $P$ is a\/ $\fs11$ set in the product\/ $\pX\ti\pY$ of 
two Polish spaces\/ $\pX$, and every section\/ $\seq Px$  
is a \bou\ set, then 
there is a sequence of Borel sets\/ $P_n\sq \pX\ti\pY$ 
with compact sections\/ $\seq{P_n}{x}$ such that\/    
$P\sq\bigcup_nP_n$.\qed
\ete

But at the moment it seems that no conclusive theory 
of $\fs11$ sets with \ddd\fsg compact sections 
(as opposed to those with \bou\ sections) 
is known. 
For instance what about effective decompositions of such 
sets into countable unions of definable sets with compact 
sections? 
Our Theorem~\ref{pt} can be used to show that such a 
decomposition is possible, but the decomposing sets 
with compact sections appear to be excessively complicated 
(3rd projective level by rough estimation). 
It is an interesting {\ubf problem} to improve this result 
to something more reasonable like Borel combinations of 
$\fs11$ sets.

On the other hand, it is known from \cite{sri,sti} that 
$\fs11$ sets with \ddd\fsg compact sections are not 
necessarily decomposable into countably many $\fs11$ 
sets with compact sections.

\vyk{
\bte
\lam{ert}
If\/ $\rE$ is a\/ $\id11$ \er\ on\/
%a\/ $\id11$ set\/ $D\sq
$\bn$ then one and only one of 
the following two claims \ref{ert1}, \ref{ert2} holds$:$
\ben
\tenu{{\rm(\Roman{enumi})}}
\itla{ert1}\msur
$\bn$ is a union of a \bou\/  
$\id11$ set and at most 
countably many\/ $\id11$ \ddd\rE classes$;$ 

\itla{ert2}
there is a 
%pairwise\/ \ddd\rE inequivalent 
\spp\/ \pis\rE\/ $Y\sq \bn$.  
\een
\ete 
}

\punk{Counterexamples above $\is11$}
\las{abo}

Here we outline several counterexamples to Theorems \ref{mt'} 
and \ref{mt} with sets $A$ more complicated than $\is11$.

\bex
\lam{abo1}
Suppose that the universe is a Cohen real extension 
$\rL[a]$ of the constructible universe $\rL$.
The set $A=\bn\cap\rL$ is $\is12$ and it is 
not \bou\  in $\rL[a]$. 
On the other hand, it is known from \cite{gs} that $A$ has 
no perfect subsets, let alone \spp\ ones. 
Thus $A$ is a $\is12$ counterexample to both 
Theorem~\ref{mt'} and Theorem~\ref{mt} in $\rL[a]$.
We then immediately obtain a similar $\ip11$ counterexample, 
using the $\ip11$ uniformization theorem.\qed
\eex

%We reiterate that it is likely open whether Theorem~\ref{mt} 
%holds for $\is11$.

\bex
\lam{abo2}
Suppose that the universe is a dominating real extension 
$\rL[d]$ of $\rL$.
%This leads to another counterexample. 
The set $A=\bn\cap\rL$ is then \bou\ in $\rL[d]$. 
The dominating forcing is homogeneous 
enough for any $\od$ (ordinal-definable) real in $\rL[d]$ 
to be constructible, and hence it is true in $\rL[d]$ that 
$A$ cannot be covered by a countable union of $\od$ compact 
sets in $\rL[d]$. 
Thus $A$ is a $\is12$ counterexample to 
Corollary~\ref{Cbou}.\qed
\eex 

Yet it is not clear how a similar 
$\ip11$ counterexample, or even a 
$\is12$ counterexample to Corollary~\ref{Ckom}, 
can be produced.

\bex
\lam{abo3}
Let $A=\ans y$ be a $\ip11$ singleton such that $y$ is not 
a $\id11$ real. 
Then conditions \ref{mt'1}, \ref{mt'2} 
of Theorem~\ref{mt'} obviously fail for $A$. 

The same for Theorem~\ref{mt}.

Moreover $A$ is a $\ip11$ counterexample to 
Corollary~\ref{Cbou} as well, although not as strong as 
those given in Example~\ref{abo2}.
\qed
\eex

It is known that there is a countable $\ip11$ set $A\sq\bn$ 
containing at least one non-$\id12$ element.
Can it serve as a more profound $\ip11$ counterexample than 
the singleton $A$ of Example~\ref{abo3}\,?

\vyk{

\punk{Generalizing the \bou{}ness dichotomy theorem}
\las{gen}

Let $\cI\sq\pws\dN$ be an ideal on $\dN$. 
A tree $T\sq\nse$ is:
\bde
\item[\it\ddi small,]
if for any $u\in T$ the set 
$\suc Tu=\ens{n}{u\we n\in T}$ belongs to $\cI$;

\item[\it\ddi positive,]
if 1) it is perfect, and 
2) if $u\in\bran T$ then the set 
$\suc Tu$ does \rit{not} belong to $\cI$.
\ede
Accordingly, a set $X\sq\bn$ is:
\bde
\item[\it\ddi small,]
if $\der X=\ens{x\res n}{n\in\dN\land x\in X}$ is an 
\ddi small tree;

\item[\it\dsi small,]
if it is a countable union of \ddi small sets;

\item[\it\ddi positive,]
if it contains a subset of the form $\bod{T}$, 
where $T\sq\nse$ is an \ddi positive tree.
\ede
For instance, if $\cI=\fin$ is the Frechet ideal of all 
finite sets $x\sq\dN$ then \ddi small trees and sets are 
exactly compact trees, resp., subsets of compact sets, 
\dsi small sets are \bou\ sets,
while \ddi positive trees are \spt s, and 
if $T$ is a \ddd\fin positive tree then the 
set $\bod{T}$ is superperfect, hence, 
non-\ddd\fsg compact.
Thus condition \ref{mt'2} of Theorem~\ref{mt'} can be 
reformulated as follows:
\rit{$A$ is a\/ \ddd\fin positive set.}

Here we prove the following theorem 
(compare with Theorem~\ref{mt'}).

\bte
\lam{nt}
Let\/ $\cI$ be a\/ $\ip11$ ideal on\/ $\dN$.
If\/ $A\sq\bn$ is a\/ $\is11$ set then one and only one of 
the following two claims holds$:$
\ben
\tenu{{\rm(\Roman{enumi})}}
\itla{nt1}\msur
$A$ is\/ \dsi small$;$ 

\itla{nt2}
$A$ is an\/ \ddi positive set.
\een
\ete

Condition \ref{nt1} of this theorem is notably weaker than 
a true generalization of Theorem~\ref{mt'} would require:
no \ddd{\id11}effectivity!
Unfortunately such a stronger version is not accessible 
so far. 
The key element in the proof of Theorem~\ref{mt'}, which 
allows to strengthen \ref{mt'1} this way, 
%from $\is11$ to $\id11$, 
is Lemma~\ref{tkm*} based on Fact~\ref{ks=t}.
We don't know whether the latter is true in the context 
of Theorem~\ref{nt}, \eg, at least in the form: 
\rit{any \ddi small\/ $\is11$ set is covered by an\/ 
\ddi small\/ $\id11$ set}. 
It would be sufficient to assume that $\cI$ satisfies the 
following property, obviously true for $\cI=\fin$: 
\rit{
if\/ $p\in\bn$ and\/ $x\in\cI$ is a\/ $\is11(p)$ set then 
there is a\/ $\id11(p)$ set\/ $y\in\cI$ such that\/ $x\sq y$.}

\bpf
As covering of small $\is11$ sets by small $\id11$ sets is 
not available, we'll follow a line of arguments which 
somewhat differs 
from the proof of Theorem~\ref{mt'} above.
First of all, let $P\sq\bn\ti\bn$ be a $\ip01$ set such that 
$A=\pr P=\ens{x\in\bn}{\sus y\:P(x,y)}$. 
Consider the tree
$$
S=\ens{\ang{x\res n,y\res n}}
{n\in\dN\land \ang{x,y}\in P} \sq\nse\ti\nse,
$$
so that 
$P=\bod{S}=\ens{\ang{x,y}\in\bn^2}
{\kaz n\:(\ang{x\res n,y\res n}\in S)}$.
If $u,v\in\nse$ then let 
$P_{uv}=\ens{\ang{x,y}\in P}{u\su x\land v\su y}$ and 
$A_{uv}=\pr P_{uv}$. 
Thus, in particular, $P_{\La\La}=P$ and $A_{\La\La}=A$. 
If the subtree 
$$
S'=\ens{\ang{u,v}\in S}
{A_{uv}\text{ is not \dsi small}}
$$ 
of $S$ is empty then $A=A_{\La\La}$ is \dsi small, 
getting \ref{nt1} of the theorem. 
Therefore we assume that $S'\ne\pu$, and the goal is 
to get \ref{nt2} of the theorem.  

Note that $P_{uv}=\bigcup_{k,n}P_{u\we k,v\we n}$, and hence 
the tree $S'$ has no maximal nodes: if $\ang{u,v}\in S'$ then 
$\ang{u\we k,v\we n}\in S'$ for some $k\yi n$.
We consider the corresponding closed set  
$$
P'=\bod{S'}=\ens{\ang{x,y}\in\bn^2}
{\kaz n\:(\ang{x\res n,y\res n}\in S')}
$$
and the $\fs11$ set $A'=\pr{P'}$. 
If $\ang{u,v}\in S'$ then let 
$$
P'_{uv} =\ens{\ang{x,y}\in P'}
{u\su x\land v\su y}
\quad\text{and}\quad
A'_{uv}=\pr P'_{uv}\,,
$$
%so that 
%$S'_{uv}$ is a subtree of $S'$, $P'_{uv}$ is 
so that $A'_{uv}$ is a non-empty $\fs11$ subset of $A'$, 
not \dsi small by the definition of $S'$. 
The next lemma is quite obvious.

\ble
\lam{tt'}
If\/ $\ang{u,v}\in S'$, $u'\in\nse$, $u\su u'$, and\/ 
$A'_{uv}\cap\ibn{u'}\ne\pu$ then there is a string\/ 
$v'\in\nse$ such that\/ $v\su v'$ and\/ 
$\ang{u',v'}\in S'$.\qed
\ele

We are going to define a pruned tree $T\sq\nse$ and a 
string $v(t)\in\nse$ for all $t\in T$, such that  
%$\ang{t,v(t)}\in S'$ and 
%
\ben
\itsep
\tenu{(\arabic{enumi})}
\itla{lx1}
if $t\in T$ then $\ang{t,v(t)}\in S'$; 
% $t\sq u(t)$, therefore 
%$A'_{u(t),v(t)}\sq \ibn t=\ens{a\in\bn}{t\su a}$;

\itla{lx2}
if $s,t\in T$ and $s\sq t$ then 
%$u(s)\sq u(t)$ and 
$v(s)\sq v(t)$;

\itla{lx3}
if $s\in T$ then there exists a string $t\in T$ such that
$s\su t$ and the set $\ens{k}{t\we k\in T}$ does not 
belong to $\cI$.
\een
If such construction is accomplished then $T$ is an 
\ddi positive tree by \ref{lx3}, 
and on the other hand $\bod{T}\sq A'\sq A$, 
so that \ref{nt2} of the theorem holds. 

Thus it remains to carry out the construction. 

To begin with we define $\La\in T$, of course, and let 
$v(\La)=\La$.

Suppose that $t\in T$, so that $\ang{t,v(t)}\in S'$ and  
%, and we have 
%$t\sq u(t)$ and $A'_{u(t),v(t)}\sq \ibn t$ by \ref{lx1}. 
the set $A'_{t,v(t)}$ is not \dsi small,
% by \ref{lx3}, 
% since $\ang{u(t),v(t)}\in S'$, 
in particular, not \ddi small, hence the tree 
$\der{A'_{t,v(t)}}$ is not \ddi small. 
We conclude that there is a string $s\in\nse$ such that 
$t\sq s$ and the set 
$K=\ens{k}{\sus a\in A'_{t,v(t)}\,(s\we k\su a)}$
does not belong to $\cI$. 

We let every string $t'$ with $t\su t'\sq s$ belong to $T$, 
and choose $v(t')$ for any such $t'$ so that \ref{lx1} and 
\ref{lx2} hold, using Lemma~\ref{tt'}. 
Then let every string $s\we k\yt k\in K$, belong to $T$, 
and let $v(s\we k)=v$, where $v$ is any string such that 
$v(s)\sq v$ and $\ang{s\we k,v}\in S'$.
(The existence of at least one such string $v$ follows from 
Lemma~\ref{tt'}.)\vom
 
\epF{Theorem~\ref{nt}}

}

\punk{Generalizing the \bou\ dichotomy: preliminaries}
\las{gen+1}

Below in Section~\ref{gen+2},
we establish a generalization of Theorem~\ref{mt'} 
for a certain system of pointset ideals which include 
the ideal of \bou\ sets along with equivalence classes 
of a given finite or countable family of \er s. 
The next definition introduces a necessary framework.

\bdf
\lam{efs}
%Let $\mil$ be the ideal of all \bou\ sets $X\sq\bn$. 
%
Let $\cF$ be a family of \er s on a set $X_0\sq\bn$.
A set $X\sq X_0$ is \rit{\ebou\cF}, 
iff it is covered by a union of the form 
$B\cup\bigcup_{n\in\om}Y_n$, where $B$ is a \bou\/ set 
and each $Y_n$ is an \ddd{\rF}equivalence class for 
an \er\ ${\rF}={{\rF}(n)}\in\cF$ which depends on $n$.

A set $X\sq X_0$ is \rit{\bol\cF},
if it is a \spp\ \pis\rF\ 
(\ie, a \tra{\rF})
for every ${\rF}\in\cF$.
\vyk{
Let $\mil(\cF)$ 
(the \rit{\ddd\cF extension of\/ $\mil$}) 
be the \ds ideal of all \ebou\cF\ sets $X\sq\bn$.
}%
\edf

Clearly \ebou\cF\ sets form a \ds ideal containing
all \bou\ sets, and
no \ebou\cF\ set can be \bol\cF.
What are properties of these ideals? 
Do they have some semblance of the superperfect ideal itself? 
We begin with a lemma and a corollary afterwards, 
which show that this is indeed the 
case \poo\ the property of being \psii.
The lemma is a generalization of Corollary~\ref{Cbou}, 
of course.

\vyk{
As usual, $\mil(\rE)$ denotes $\mil(\ans\rE)$, and this is 
the \ds ideal of all \ebou\rE\ sets $X\sq\bn$, 
that is, those covered by countably 
many \de classes and a set in $\mil$. 
Notation like $\mil(\rF_1,{\rF_2},{\rF_3})$ is understood 
similarly.
}

\ble
\lam{ico}
Suppose that\/  $\sis{\rF_n}{n<\om}$ 
is a\/ $\id11$ sequence of \er s on\/ $\bn$, and a\/ 
$\is11$ set\/ $X\sq\bn$ is\/ \sm{\sis{\rF_n}{n<\om}}.
Then\/ $X$ is\/ 
\ddd{\id11}effectively\/ \sm{\sis{\rF_n}{n<\om}}, in 
the sense that there exist$:$
\ben
\aenu
\item 
a\/ $\id11$ sequence of compact trees\/ $T_k$, 

\item
a\/ $\id11$ sequence of numbers\/ $n_k$, and 

\item 
a\/ $\id11$ set\/ $H\sq\dN\ti\bn$ 
\een
such that, 
for every\/ $k<\om$ the cross-section\/
$\seq Hk=\ens{a}{\ang{k,a}\in H}$ is an\/ 
\ddd{\rF_{n_k}}equivalence class and\/ 
$X\sq\bigcup_k\bod{T_k}\cup\bigcup_k \seq Hk$.
\ele

In particular, if a $\is11$ set $X\sq\bn$ is
\sm{\sis{\rF_n}{n<\om}}
%\ddd{\id11}effective 
%\sm{\sis{\rF_n}{n<\om}}ness implies that the given set 
then $X$ is covered by the union of 
all\/ $\id11$ \ddd{\rF_0}classes, 
all\/ $\id11$ \ddd{\rF_1}classes, 
all\/ $\id11$ \ddd{\rF_2}classes, et cetera, 
and all\/ $\id11$ compact sets.

\bpf
The set $C=\ct\cap\id11$ of all $\id11$ compact trees
is $\ip11$, and hence so is $K=\bigcup_{T\in C}\bod{T}$. 
If $n<\om$ then let $U_n$ be the union of all $\id11$ 
\ddf nclasses. 
Let's show that $U={\bigcup_nU_n}$ is $\ip11$ either.
We make use of sets $E\sq\dN$ and  
$W,W'\sq\dN\ti\bn$ as in Fact~\ref{23}. 
The $\ip11$ formula $\vpi(e,n)\::=$
$$
e\in E 
\:\land\:
\kaz y,z\in\seq{W'}e\:(y\rF_n z)
\:\land\:
\kaz y\in\seq{W'}e\:\kaz z\:({y\rF_n z}\imp z\in\seq{W}e)
$$
says that $e\in E$ and  
$\seq{W'}e=\seq{W}e$ is a \ddd{\rF_{n}}equivalence class.
Moreover
$$
x\in U\leqv
\sus n\:\sus e\:(\vpi(e,n)\land x \in\seq{W}e).
$$
\vyk{\big(
e\in E\land 
%x\in\seq We\land 
\underbrace
{\kaz y\:(y\in\seq{W'}e\imp {x\rF_n y}\imp y \in\seq{W}e)}
_{\vpi(n,e,x)}  
\big), 
$$
which is a $\ip11$ formula, of course. 
(If $e\in E$ then the subformula $\vpi(n,e,x)$ says that the set 
$\seq{W'}e=\seq{W}e$ is a \ddd{\rF_{n}}equivalence class 
containing $x$.)
\vom
}%

\rit{Case 1}: 
$X\sq K\cup U$.
Then the set $S$ of all pairs $\ang{x,h}$ such that  
\bit
\item[$-$]
either $h=T\in C$ and $x\in\bod{T}$,

\item[$-$]
or $h=\ang{e,n}\in\Phi=
\ens{\ang{e,n}\in E\ti\dN}{\vpi(e,n)}$ and $x\in\seq We$,
\eit
is a $\ip11$ set satisfying $X\sq\pr S$.
\vyk{
Arguing as in the proof of Theorem~\ref{mt'} (Case 1),
we obtain a (countable) $\id11$ set $\tau$ which consists
of compact $\id11$ trees and a $\id11$ set $G\sq E\ti\dN$
such that if $x\in X$ then either $x\in\bod{T}$ for some
$T\in\tau$ 
}%
By Fact~\ref{22+} there is a $\id11$ map $f$ defined on
$\bn$ and such that $\ang{a,f(a)}\in S$ for each $a\in X$.
The sets
$$
X'=\ens{x\in X}{f(x)\in\ct}
\qand
X''=\ens{x\in X}{f(x)\in \Phi}
$$
are $\is11$ as well as their images
$$
R'=\ens{f(x)}{x\in X'}\sq C
\qand
R''=\ens{f(x)}{x\in X''}\sq\Phi\,,
$$
and $X'\cup X''=X$, $R'\cup R''=\ens{f(x)}{x\in X}$.
By the $\is11$ Separation theorem there is a $\id11$ set $\tau$
such that $R'\sq\tau\sq C$, and by Fact~\ref{25c} we have
$\tau=\ens{T_k}{k<\om}$, where $k\mto T_k$ is a $\id11$ map.
By similar reasons, there is a $\id11$ map
$k\mto\ang{e_k,n_k}$ such that 
$R''\sq \rho=\ens{\ang{e_k,n_k}}{k<\om}\sq \Phi$. 
To finish the proof in Case 1, it remains to define 
$$
H
=\ens{\ang{k,x}\in\dN\ti\bn}{x\in\seq{W}{e_k}}  
=\ens{\ang{k,x}\in\dN\ti\bn}{x\in\seq{W'}{e_k}}\,.
$$ 

\rit{Case 2}:  
$A=X\bez{(K\cup U)}\ne \pu$. 
Then $A$ is a non-empty $\is11$ set. 
We are going to derive a contradiction.
By definition, we have 
$X\sq \bigcup_kC_k\cup \bigcup_n\bigcup_k E_{nk}$, 
where each $C_k$ is compact and each $E_{nk}$ is an 
\ddf nclass.
Let $M$ be a countable elementary substructure of a 
sufficiently large structure, containing, in particular, 
the whole sequence of covering sets $C_k$ and $E_{nk}$. 
Below ``generic'' will mean \gh\ generic over $M$.

As $A\ne\pu$ is $\is11$, there is a perfect set $P\sq A$ 
of points both generic and pairwise generic. 
It is known that then $P$ is a \pis{\rF_n} for every $n$, 
hence, definitely a set not covered by a countable union 
of \ddf nclasses for all $n<\om$. 
Thus to get a contradiction it suffices to prove that 
$P\cap C_k=\pu$ for all $k$. 
In other words, we have to prove that if $k<\om$ and 
$x\in A$ is any generic real then $x\nin C_k$. 

Suppose towards the contrary that a non-empty $\is11$ 
condition $Y\sq A$ forces that $\ja\in C_k$, where 
$\ja$ is a canonical name for the \gh\ generic real. 
We claim that $Y$ is not \bou. 
Indeed otherwise we have $Y\sq \bigcup_n \bod{T_n}$ by 
Theorem~\ref{mt'}, where all trees $T_n\sq\nse$ are 
$\id11$ and compact, which contradicts the fact that $A$ 
%(and hence $X$ as well) 
does not intersect any compact $\id11$ set. 

Therefore $Y\not\sq C_k$. 
Then there is a point $x\in Y$ and a number $m$ such 
that the set $I=\ens{y\in\bn}{y\res m=x\res m}$ 
does not intersect $C_k$. 
But then the $\is11$ condition $Y'=Y\cap I$ 
forces that $\ja\nin C_k$, a contradiction.
\epf

\bcor
\lam{mpsii}
If\/  $\sis{\rF_n}{n<\om}$ 
is a\/ $\id11$ sequence of \er s on\/ $\bn$ 
%If\/  $\cF$ is a countable family of Borel \er s on\/ $\bn$ 
then\/
%$\mil(\cF)$
the ideal of\/ \ebou{\sis{\rF_n}{n<\om}}\ sets
is\/ \lsii\ and\/ \psii.
%\ and\/ \psii. 
\ecor

See \cite[section 3.8]{id} on \lsii\ and \psii\ ideals.

\bpf
Consider a $\is11$ set\/ $P\sq\bn\ti\bn$. 
We have to prove that 
$$
X=\ens{x\in\bn}
{\seq Px=\ens{y}{\ang{x,y}\in P}\,
\text{ is \ebou{\sis{\rF_n}{n<\om}}}}
$$
is a $\ip11$ set. 
%Let $T_k\yi n_k\yi H$ be as in Lemma~\ref{ico}. 
%
%, where, we recall, $\seq Px=\ens{y}{\ang{x,y}\in P}$.
By the relativized version of Lemma~\ref{ico}, 
$x\in X$ iff 
\ben
\fenu
\itla{ico*}
there exist $\id11(x)$ sequences 
$\sis{T_k}{k<\om}$ (of compact trees) 
and $\sis{n_k}{k<\om}$ and a $\id11(x)$ 
set $H\sq\om\ti\bn$ such that, 
for every\/ $k<\om$ the cross-section\/ $\seq Hk$ is an\/ 
\ddd{\rF_{n_k}}equivalence class and\/ 
$\seq Px\sq\bigcup_k\bod{T_k}\cup\bigcup_k \seq Hk$.
\een
A routine analysis (as in the proof of Lemma~\ref{ico}) 
shows that this is a $\ip11$ 
description of the set $X$.
\epf

\punk{Digression: another look on the effectivity}
\las{dig}

 As usual,
Lemma~\ref{ico} and Corollary~\ref{mpsii} remain true for
relativized classes.
In particular, if $p\in\bn$, $\rF_n$ are $\id11(p)$ \er s, 
and a $\is11(p)$ set $X\sq\bn$
is \sm{\sis{\rF_n}{n<\om}}
then $X$ is covered by the union of 
all\/ $\id11(p)$ \ddd{\rF_n}classes, $n=0,1,2,\dots$, 
and all $\id11(p)$ compact sets.
If now $p\ne q\in\bn$ is a different parameter, but still 
$\rF_n$ are $\id11(q)$ and $X$ is $\is11(q)$ and  
\sm{\sis{\rF_n}{n<\om}}
then accordingly $X$ is covered by the union of 
all\/ $\id11(q)$ \ddd{\rF_n}classes, $n=0,1,2,\dots$, 
and all $\id11(q)$ compact sets.
Those two countable coverings of the same set $X$ can be 
different, of course. 
This leads to the question: 
is there a covering of $X$ of the type indicated, which 
depends on $X$ and $\rF_n$ themselves, 
but not on the choice of a 
parameter $p$ such that $X$ is $\is11(p)$ and $\rF_n$ are 
$\id11(p)$. 
We are able to answer this question in the positive at least 
in the case of finitely many \er s. 
The next theorem will be instrumental in the proof of a 
theorem in Section~\ref{uu2}.

\bte
\lam{ico+}
Suppose that\/ $n\ge1$, $\rF_1,\dots,\rF_n$ are  
Borel \er s on\/ $\bn$, and a\/ 
$\fs11$ set\/ $X\sq\bn$ is\/ \sm{\ans{\rF_1,\dots,\rF_n}}.
Then there exist Borel sets\/ $Y_1,\dots,Y_n,X_{n+1}\sq\bn$ 
such that\/ 
\ben
\renu
\itla{ico+1}\msur
$X\sq Y_1\cup\dots\cup Y_n\cup X_{n+1}$, 

\itla{ico+2}
each set\/ $Y_j$ is a countable union of\/ 
\ddd{\rF_{j}}equivalence classes while the set\/ $X_{n+1}$ 
is\/ \bou,

\itla{ico+3}
%there is a\/ $\id12$ map\/ $p\mto\bap$ such that 
if\/ $p\in\bn,$ $X$ is\/ $\is11(p)$, 
and all relations\/ $\rF_m$ are\/ $\id11(p)$, then 
there is a parameter\/ $\bap\in\bn$ in\/ $\id12(p)$
such that both\/ $X_{n+1}$ and all sets\/ $Y_j$ are\/ 
$\id11(\bap)$ --- hence, $\id12(p)$.
\een
\ete

This, under the assumptions of the theorem, there is a Borel 
covering of $X$ satisfying \ref{ico+1} and \ref{ico+2}, and 
effective as soon as $X$ and $\rF_j$ are granted some 
effectivity. 
It is a challenging problem to get rid of $\bap$ in 
\ref{ico+3}
(so that $X_{n+1}$ and all $Y_j$ are just $\id11(p)$
with the same $p$), but this remains open.

\bpf
We define sets 
$X=X_1\qs X_2\qs X_3\qs\dots\qs X_n\qs X_{n+1}$ 
so that $X_{j+1}=X_j\bez Y_j$, where by induction
$$
Y_j=\ens{x\in\bn}{\text{the set }\,X_j\cap \fn xj\,
\text{ is not \sm{\ans{\rF_j,\dots,\rF_n}}}}.
\eqno(1)
$$
In particular,
$$
\bay{rcl}
Y_1 &=& 
\ens{x\in\bn}{\text{the set }\,X_1\cap \fn x1\,
\text{ is not \sm{\ans{\rF_2,\dots,\rF_n}}}}\,,\\[1ex]
Y_2 &=&
\ens{x\in\bn}{\text{the set }\,X_2\cap \fn x2\,
\text{ is not \sm{\ans{\rF_3,\dots,\rF_n}}}}\,,\\[1ex]
&\dots\\[1ex]
Y_{n-1} &=&
\ens{x\in\bn}{\text{the set }\,X_{n-1}\cap \fn x{n-1}\,
\text{ is not \sm{\ans{\rF_n}}}}\,,\\[1ex]
Y_{n} &=&
\ens{x\in\bn}{\text{the set }\,X_{n}\cap \fn xn\,
\text{ is not \sm{\pu}}}\,, 
\eay
$$
where \sm\pu\ is the same as just \bou.

\ble
\lam{icoS1}
If\/ $1\le j\le n$ then\/ $Y_j$ is a countable union of\/ 
\ddd{\rF_j}equivalence classes and the set\/ 
$X_{j+1}=X_j\bez Y_j$ is\/ \sm{\ans{\rF_{j+1},\dots,\rF_n}}.
\ele
\bpf
Let $\cY_j$ be the family of all sets $Y$ such that $Y$ 
is a union of at most countably many \ddd{\rF_j}classes 
and $X_j\bez Y$ is \sm{\ans{\rF_{j+1},\dots,\rF_n}}.
Note that $\cY_j$ is a non-empty 
(since $X_j$ is \sm{\ans{\rF_j,\dots,\rF_n}} by induction) 
\ds filter 
(since the collection of all \sm{\ans{\rF_{j+1},\dots,\rF_n}} 
sets is a \ds ideal). 
Therefore $Y'_j=\bigcap\cY_j$ is a set in $\cY_j$, in fact, 
the \ddd\sq least set in $\cY_j$.

It remains to show that $Y_j=Y'_j$.
We claim that if $C$ is an \ddd{\rF_j}class then $C\sq Y'_j$ 
iff $C\sq Y'_j$. 
%the set $X_1\cap C$ is {\ubf not} \sm{\ans{\rF_{j+1},\dots,\rF_n}}. 
Indeed if $C\cap Y_j=\pu$ then  
$X_j\cap C$ is \sm{\ans{\rF_{j+1},\dots,\rF_n}}, thus  
$Y'_j\bez C$ is still a set in $\cY_j$, 
therefore $C\cap Y'_j=\pu$.
Conversely if $C\cap Y'_j=\pu$ then 
$(X_j\cap C)\sq (X_j\bez Y'_j)$, 
and hence $X_j\cap C$ is \sm{\ans{\rF_{j+1},\dots,\rF_n}}, 
so $C\cap Y_j=\pu$, as required. 
\epF{Lemma}

Thus by the lemma the sets $Y_j$ and $X_{n+1}$ satisfy 
\ref{ico+1} and \ref{ico+2} of the theorem. 
To verify \ref{ico+3}, assume that $p\in\bn$, 
$X$ is $\is11(p)$, and all $\rF_m$ are $\id11(p)$.
The main issue is that the sets $Y_j$, albeit Borel 
(as countable unions of Borel equivalence classes) 
do not seem to be $\id11(p)$, at least straightforwardly. 
For instance, $Y_1$ is $\is11(p)$ 
by Corollary~\ref{mpsii} (relativized), 
and accordingly $X_2$ is $\ip11(p)$ 
(instead of $\id11(p)$), 
which makes it very 
difficult to directly estimate the class of $Y_2$ at 
the nest step. 
This is where a new parameter appears.

We precede the last part of the proof of the theorem with 
the following auxiliary fact on equivalence relations, 
perhaps, already known.

\ble
\lam{cer}
Let\/ $\rE$ be a\/ $\id11$ \er\ on\/ $\bn$, and\/ 
$X\sq\bn$ be a\/ $\is11$ set which
intersects only countably many\/ \de classes. 

Then all\/ \de classes\/ $\ek x\yt x\in X,$ are\/ 
$\id11$ sets, and there is an\/ \de invariant\/ 
$\id11$ set\/ $Y\sq\bn$ such that\/ $X\sq Y$ 
and all\/ \de classes\/ $\ek y$, $y\in Y$, are\/ 
$\id11$ sets
(therefore\/ $Y$ still contains only countably
many\/ \de classes).
\ele
\bpf 
The union $C$ of all $\id11$ \de classes is
an \de invariant $\ip11$ set. 
(See, \eg, 10.1.2 in \cite{k-ams}.) 
Thus, if $X\not\sq C$ then $H=X\bez C$ is a non-empty 
%\de invariant
$\is11$ set which does not intersect $\id11$ \de classes.
Then  
(see, \eg, Case 2 in the proof of Theorem 10.1.1 in \cite{k-ams}) 
$H$ contains a perfect \pis\rE, 
% pairwise \de inequivalent subset, 
which contradicts our assumptions. 
Therefore $X\sq C$, so indeed 
all \de classes\/ $\ek x$, $x\in X$, are $\id11$. 
To prove the second claim apply the invariant 
$\is11$ separation theorem
(see, \eg, 10.4.2 in \cite{k-ams}), 
which yields an \de invariant $\id11$ set $Y$ satisfying 
$X\sq \ek X\sq Y\sq C$.
\epF{Lemma}

We continue the proof of Theorem~\ref{ico+}.
The next goal is to find a parameter $q_1\in\bn$ in $\id12(p)$
such that the $\is11(p)$ set $Y_1$ is $\id11(q_1)$.
Let $\ip11$ sets\/ $\bE\sq\bn\ti\dN$ and  
$\bW\sq\bn\ti\dN\ti\bn$, 
and a $\is11$ set\/ $\bW'\sq\bn\ti\dN\ti\bn$ 
be as in Lemma~\ref{23+}.
Let $E(p)=\ens{e}{\ang{p,e}\in\bE}$ and, for all $e<\om$,
$$
W_e(p)=\ens{x}{\ang{p,e,x}\in\bW}\,,
\quad
W'_e(p)=\ens{x}{\ang{p,e,x}\in\bW'}\,,
$$
so that $E(p)$ and all sets $W_e(p)$ are $\ip11(p)$ while all 
sets $W'_e(p)$ are $\is11(p)$. 
By Lemma~\ref{cer} (relativized), a 
point $x\in\bn$ belongs to $Y_1$ iff 
$$
\bay{l}
\sus e\:
\big(
e\in E(p)\,\land\, x\in W_e(p)\,\land\, 
W_e(p)\,\text{ is an \ddd{\rF_1}class}\,\land\,
{\hspace*{15ex}}\\[1ex]
{\hspace*{15ex}}
\,\land\,W'_e(p) \cap X_1\,
\text{ is not \sm{\ans{\rF_2,\dots,\rF_n}}}
\big).
\eay
$$
The first line is $\ip11(p)$. 
(Note that $W_e(p)=W'_e(p)$ for all $e\in E(p)$.) 
The second line is only $\is11(p)$ by Corollary~\ref{mpsii}.
However the set 
$$
Q_1(p)=\ens{e\in E(p)}
{W'_e(p) \cap X_1\,\text{ is not \sm{\ans{\rF_2,\dots,\rF_n}}}}
\sq\om
$$
is $\id12(p)$ 
(more precisely, an intersection of $\ip11(p)$ and $\is11(p)$), 
and 
$$
x\in Y_1
\leqv 
\sus e\in E(p)\cap Q_1(p)\:
\big(
x\in W_e(p)\,\land\, 
W_e(p)\,\text{ is an \ddd{\rF_1}class}\big).
$$
We conclude that $Y_1$ is $\id11(p,Q_1)$, hence, $\id11(q_1)$, 
where $q_1\in\bn$ is a ``concatenation'' of $p$ and $Q_1$ 
(so that $q_1$ is $\id12(p)$).

Arguing the same way, we find parameters $q_2\yi q_3\yi\dots$ 
such that each $Y_j$ is $\id11(q_j)$ and each $q_{j+1}$ is 
$\id12(q_j)$, and hence $\id12(p)$ by induction. 
Wrapping this construction up in a parameter $\bap$ 
as in \ref{ico+3} is a routine.
%This ends the proof of the theorem.
\epf

We don't know whether the theorem still holds for   
countably infinite sequences of \er s. 
Yet the proof miserably fails in this case. 
Indeed, let, for any $n$, $\rF_n$ be an \er\ on $\bn$ 
whose classes are $I_k=\ens{x\in\bn}{x(0)=k}$, 
$k=0,1,\dots,n$, and all singletons outside of these 
large classes. 
The whole space $\bn=\bigcup_nI_n$ is   
\sm{\ans{\rF_0,\rF_1,\rF_2,\dots}}, of course.
But running the construction as above, we'll obviously 
have $Y_0=Y_1=Y_2=\dots=\pu$ 
(as each \ddd{\rF_n}class is covered by an appropriate 
\ddd{\rF_{n+1}}class), 
which results in nonsense. 

There is another interesting problem. 
Under the assumptions of the theorem, the covering of $X$ 
by sets $Y_1,\dots,Y_n,X_{n+1}\sq\bn$ depends on $X$ 
but is independent of the choice of a parameter $p$ as 
in \ref{ico+3}. 
On the other hand, if such a parameter $p$, and accordingly 
$\bar p$ as  in \ref{ico+3}, is given then not only each 
$Y_j$ but also a representation of $Y_j=\bigcup_mY_{jm}$ 
as a countable union of \dd{\rF_j}classes $Y_{jm}$, 
can be obtained in $\id11(\bap)$ by Lemma~\ref{ico}.
One may ask whether such a decomposition of each $Y_j$ is 
available in a way independent of the choice of $p$ 
(as the sets $Y_j$ themselves). 
The answer in the negative is expected, but it may likely 
take a lot of work. 
On the other hand, Theorem~\ref{resT} below will show that, 
under some restrictions, if a countable union of 
equivalence classes 
of a $\id11$ \er\ is $\id11(\xi)$, where $\xi<\omi$, 
then all classes in  
this union admit constructible (not necessarily countable) 
Borel codes.

\vyk{
\rit{Step 1}. 
Let $X_1=X$.
Let $\cY_1$ be the family of all sets $Y$ such that $Y$ 
is a union of at most countably many \ddd{\rF_1}classes 
and $X_1\bez Y$ is \sm{\ans{\rF_2,\dots,\rF_n}}.
Note that $\cY_1$ is a non-empty 
(since $X$ is \sm{\ans{\rF_1,\dots,\rF_n}}) \ds filter 
(since the collection of all \sm{\ans{\rF_2,\dots,\rF_n}} 
sets is a \ds ideal). 
Therefore $Y_1=\bigcap\cY_1$ is a set in $\cY_1$, in fact, 
the \ddd\sq least set in $\cY_1$.

We claim that if $C$ is an \ddd{\rF_1}class then $C\sq Y_j$ 
iff $C\sq Y'_j$.
%the set $X_1\cap C$ is {\ubf not} \sm{\ans{\rF_2,\dots,\rF_n}}. 
Indeed if $C\cap Y_1=\pu$ then 
$X_j\cap C$ is \sm{\ans{\rF_{j+1},\dots,\rF_n}}, thus   
$Y'_j\bez C$ is still a set in $\cY_j$, 
therefore $C\cap Y'_j=\pu$.
Conversely if $C\cap Y'_j=\pu$ then 
$(X_j\cap C)\sq (X_j\bez Y_j)$, 
and hence $X_j\cap C$ is \sm{\ans{\rF_{j+1},\dots,\rF_n}}, as 
required. 
To conclude, 
$$
Y_1=\ens{x\in\bn}{\text{the set }\,X_1\cap \fn x1\,
\text{ is not \sm{\ans{\rF_2,\dots,\rF_n}}}}.
\eqno(1)
$$

\rit{Step 2}. 
The difference $X_2=X_1\bez Y_1$ is 
\sm{\ans{\rF_2,\dots,\rF_n}} 
by construction. 
Let $\cY_2$ = all sets $Y$ such that $Y$ 
is a union of at most countably many \ddd{\rF_2}classes 
and $X_2\bez Y$ is \sm{\ans{\rF_3,\dots,\rF_n}}.
Once again, the set $Y_2=\bigcap\cY_2$ belongs to 
$\cY_2$ and is the least set in $\cY_2$.
Moreover,\pagebreak[0] 
$$
Y_2=\ens{x\in\bn}{\text{the set }\,X_2\cap \fn x2\,
\text{ is not \sm{\ans{\rF_3,\dots,\rF_n}}}}.
\eqno(2')
$$
We claim that in addition 
\ben
\aenu
\atc
\itla{aa1}
if $x\in\bn$ then $x\in Y_2$ if and only if 
for each collection $\cZ$ of \ddd{\rF_1}classes $Z\sq Y_1$, 
the set $(X\bez \bigcup\cZ)\cap \fn x2$ 
is not \sm{\ans{\rF_3,\dots,\rF_n}}. 
%where as usual $\fn A1=\bigcup_{a\in A}\fn a1$.
\een
Indeed let $x$ belong to the right-hand side of \ref{aa1}. 
Let $\cZ=\cY_1$, so that $\bigcup\cZ=Y_1$. 
Then $X_2=X_1\bez\bigcup\cZ$ and 
easily $x\in Y_2$. 
Conversely, let $x\in Y_2$, and let $\cZ$ be a 
collection of \ddd{\rF_1}classes $Z\sq Y_1$. 
Then $\bigcup\cZ\sq Y_1$, therefore 
$X\bez \bigcup\cZ \qs X_2=X_1\bez Y_1$, and hence
%the set 
$$
\textstyle
(X_1\bez \bigcup\cZ)\cap \fn x2 \qs
X_2\cap \fn x2\,.
$$ 
But $X_2\cap \fn x2$ is not \sm{\ans{\rF_3,\dots,\rF_n}} 
by $(2')$ since $x\in Y_2$.\vom

\rit{Step 3}. 
Once again, the set $X_3=X_2\bez Y_2$ is 
\sm{\ans{\rF_3,\dots,\rF_n}}. 
Let $\cY_3$ = all sets $Y$ such that $Y$ 
is a union of at most countably many \ddd{\rF_3}classes 
and $X_3\bez Y$ is \sm{\ans{\rF_4,\dots,\rF_n}}.
As above, the set $Y_3=\bigcap\cY_3$ belongs to 
$\cY_3$ and is the least set in $\cY_3$.
Moreover,\pagebreak[0] 
$$
Y_3=\ens{x\in\bn}{\text{the set }\,X_3\cap \fn x3\,
\text{ is not \sm{\ans{\rF_3,\dots,\rF_n}}}},
\eqno(3')
$$
and in addition 
\ben
\aenu
\atc
\atc
\itla{aa2}
if $x\in\bn$ then $x\in Y_3$ iff for 
any collections $\cZ_1$ of \ddd{\rF_1}classes $Z\sq Y_1$ 
and $\cZ_2$ of \ddd{\rF_2}classes $Z\sq Y_2$, the set 
$(X\bez (\bigcup\cZ_1\cup\bigcup\cZ_2))\cap \fn x3$ 
is not \sm{\ans{\rF_4,\dots,\rF_n}}.
\een
}

\punk{Generalizing the \bou\ dichotomy: the theorem}
\las{gen+2}

Coming back to the content of Section~\ref{gen+1},
we'll prove the following theorem 
in this section.

\bte
[common with Marcin Sabok and Jindra Zapletal]
\lam{fm}
Suppose that\/ $n<\om$, 
$\rF_1,\dots,\rF_n$ are\/ $\id11$ \er s on\/ $\bn,$
and\/ $A\sq\bn$ is a\/ $\is11$ set.
Then one and only one of 
the following two claims holds$:$
\ben
\tenu{{\rm(\Roman{enumi})}}
\itla{fm1}
the set\/ $A$ is\/ \nsm\ ---
and therefore 
\ddd{\id11}effectively\/ \sm{\sis{\rF_n}{n<\om}}
as in Lemma \ref{ico}\/$;$  
\vyk{
$A$ is covered by a countable 
union of all\/ $\id11$ \ddd{\rF_1}classes, 
all\/ $\id11$ \ddd{\rF_2}classes, 
\dots , all\/ $\id11$ \ddd{\rF_n}classes, 
and all\/ $\id11$ compact sets$;$
}

\itla{fm2}
there exists an\/ \nbol\ set\/ $P\sq A$. 
%which is an\/ \tra{\rF_i} 
%for each\/ $i=1,\dots,n$.  
\een
\ete

If $n=0$ then this theorem is equivalent to Theorem~\ref{mt'}: 
indeed, if $\cF=\pu$ then \ebou\pu\ sets are just \bou,
while \bol\pu\ sets are just \spp.

\vyk{
We precede the proof of Theorem~\ref{fm} with 
the following auxiliary lemma on equivalence relations, 
perhaps, already known.

\ble
\lam{cer}
Let\/ $\rE$ be a\/ $\id11$ \er\ on\/ $\bn$, and\/ 
$X\sq\bn$ be a\/ $\is11$ set which
intersects only countably many\/ \de classes. 

Then all\/ \de classes\/ $\ek x$, $x\in X$, are\/ 
$\id11$ sets, and there is an\/ \de invariant\/ 
$\id11$ set\/ $Y\sq\bn$ such that\/ $X\sq Y$ 
and all\/ \de classes\/ $\ek y$, $y\in Y$, are\/ 
$\id11$ sets
(therefore\/ $Y$ still contains only countably
many\/ \de classes).
\ele
\bpf[Lemma]
The union $C$ of all $\id11$ \de classes is
an \de invariant $\ip11$ set. 
(See, \eg, 10.1.2 in \cite{k-ams}.) 
Thus, if $X\not\sq C$ then $H=X\bez C$ is a non-empty 
%\de invariant
$\is11$ set which does not intersect $\id11$ \de classes.
Then  
(see, \eg, Case 2 in the proof of Theorem 10.1.1 in \cite{k-ams}) 
$H$ contains a perfect \pis\rE, 
% pairwise \de inequivalent subset, 
which contradicts our assumptions. 
Therefore $X\sq C$, so indeed 
all \de classes\/ $\ek x$, $x\in X$, are $\id11$. 
To prove the second claim apply the invariant 
$\is11$ separation theorem
(see, \eg, 10.4.2 in \cite{k-ams}), 
which yields an \de invariant $\id11$ set $Y$ satisfying 
$X\sq \ek X\sq Y\sq C$.
\epF{Lemma}
}

The following key result of
Solecki -- Spinas \cite[Theorem 2.1 and Corollary 2.2]{ss} 
will be an essential pre-requisite in the proof of 
Theorem~\ref{fm}.

\bte
\lam{pss}
Suppose that\/ $E\sq\bn\ti\bn$ and\/ {\rm(*)}
there is a decomposition\/ $E=\bigcup_nE_n$ such that
\ben
\renu
\itla{pss1} 
if\/ $n<\om$ and\/ $U\sq\bn\ti\bn$ is open then the 
projection\/ 
$\pr{(E_n\cap U)}$ has the Baire property in\/ $\bn\;;$

\itla{pss2} 
if\/ $n<\om$ and\/ $a\in\bn$ then the cross-section\/ 
$\seq{E_n}a=\ens{x}{\ang{a,x}\in E_n}$ is bounded\/ 
{\rm(= covered by a compact set)}.
\een
%$\fs11$ set such that every
%cross-section\/ $\seq Ex\yt x\in\bn,$ is\/ \bou.
Then there is a \sps\/ $P\sq\bn$ {\bfit free for\/ $E$}
in the sense that if\/ $x\ne y$ belong to\/ $P$ then\/ 
$\ang{x,y}\nin E$.
\qed
\ete

\bcor
\lam{pssC}
If\/ $E\sq\bn\ti\bn$ is a\/ $\fs11$ set and each
cross-section\/ $\seq Ea\yt a\in\bn,$ is\/ \bou,
then there is a \sps\/ $P\sq\bn$ free for\/ $E$.
In particular, if\/ $\rE$ is a\/ $\fs11$ \er\ on\/ $\bn$ 
with all\/ \de\ec es\/ \bou\ then there is a \spp\ \pis\rE. 
\ecor
\bpf[see \cite{ss}]
By Theorem~\ref{buhi} $E$ admits a decomposition satisfying 
(*) of Theorem~\ref{pss}.
We also note that if $\rE$ is an \er\ then 
a set free for $\rE$ is the same as a \pis\rE.
\epf

\bpf[Theorem~\ref{fm}]
We argue by induction on $n$. 
The case $n=0$ (then $\ans{\rF_1,\dots,\rF_{n}}=\pu$)
is covered by Theorem~\ref{mt'}. 
Now the step $n\to n+1$. 

Let $\rF_1,\dots,\rF_n,\rF_{n+1}$ be $\id11$ \er s on $\bn,$
and $A\sq\bn$ be a $\is11$ set.
The set
$$
\bU=\ens{x\in A}{\fn x1\,\text{ is \nila}}  
$$
is $\is11$ by Corollary~\ref{mpsii}.
\vom

\rit{Case 1}: 
the $\is11$ set 
$\bU$ has only countably many \ddf1classes.
Then by Lemma~\ref{cer}, there is 
an \ddf1invariant $\id11$ set $D$ such that $\bU\sq D$, 
$D$ contains only countably 
many \ddf1classes, and all of them are $\id11$.\vom

\rit{Subcase 1.1}: 
the complementary  
$\is11$ set $B=A\bez D$ is \nism.
Then 
%, by the inductive hypothesis,
the whole domain $A= D\cup B$ is \insm,
%covered by the union of 
%all $\id11$ \ddf kclasses, $k=1,\dots,n,n+1$, 
%and all $\id11$ compact sets, 
hence we have \ref{fm1} for 
$\rF_1,\dots,\rF_n,\rF_{n+1}$.\vom

\rit{Subcase 1.2}: 
$B$ is \nila. 
By the inductive hypothesis 
there is an \nibol\ set $P\sq B$.
%which is an \tra{\rF_i} 
%for each\/ $i=2,\dots,n+1$. 
Let $x\in P$. 
Then the class $\fn x1$ is \nism. 
We claim that the set $P_x=\fn x1\cap P$ is just \bou. 
Indeed by definition $P_x\sq Y\cup\bigcup_kX_k$, where 
$Y$ is \bou\ while each $X_k$ is an 
\ddd{\rF_{n(k)}}equivalence class for some 
$n(k)=2,3,\dots,n+1$.
By construction $P$ has at most one common point with each 
$X_k$. 
Therefore the set $P_x\bez Y$ is at most countable, hence, 
\bou, and we are done. 

Thus all \ddf1classes inside $P$ are \bou.
By Corollary~\ref{pssC},  
there is a \spp\ \pis{\rF_1}\ $Q\sq P$ --- 
then the set $Q$ is \inbol\ by construction. 
Thus \ref{fm2} holds.\vom 

\rit{Case 2}: 
$\bU$ has uncountably many \ddf1classes.
Then by the Silver dichotomy \cite{sil} 
there exists a perfect \pis{\rF_1}\ 
$X\sq\bU$. 
If $x\in X$ then by definition the class $\fn x1$ 
is not \nism.
Therefore by the inductive hypothesis 
%(Silver for the \bou\ ideal) 
there exists an \nibol\ set $Y\sq \fn x1$,  above
and hence a \supt\ tree $T\sq\nse$ such that $\bod{T}=Y$. 
The next step is to get such a tree $T$ by means of a Borel 
function defined on a smaller domain.

\ble
\lam{sau}
In our assumptions, there is a perfect set\/ $X'\sq X$ and
a Borel map\/ $x\mto T_x$ defined on $X'$, such that if\/
$x\in X'$ then\/ $T_x$ is a \supt\ tree, $\bod{T_x}\sq\fn x1$,
and\/ $\bod{T_x}$ is \nibol.
\ele
\bpf[Lemma]
Let $p\in\bn$ be a parameter such that $X$ is $\ip01(p)$. 

Let $\rV$ be the set universe considered, and let $\rV^+$ 
be a generic extension of $\rV$ such that $\omi^{\rL[p]}$ 
is countable in $\rV^+$. 
Let $X^+$ be the \ddd{\rV^+}extension of $X$, so that $X^+$ 
is $\ip01(p)$ in $\rV^+$ and $X=X^+\cap\rV$. 
Let $\rF_i^+$ be a similar extension of $\rF_i$. 
It is true then in $\rV^+$ by the Shoenfield absoluteness
that each 
%$D^+$ is a $\is11$ set, 
$\rF_i^+$ is a $\id11$ \er\ on $\bn$, 
and $X^+$ is a perfect set in $\ip01(p)$.
Moreover, it is true in $\rV^+$ by the Shoenfield
absoluteness that 
\ben
\fenu
\itla X
if $x\in X^+$ then 
the \ddd{\rF^+_1}class 
$\fnp x1$ is
%non-\bou\
not \nisp\  
\een
--- simply because the formula 
$$
\kaz x\in X\:
(\text{$\fn x1$ is  not \nism})
%\kaz A\:\sus y\in\ek x\:
%(A\text{ is \sik}\imp y\nin A)
$$
is essentially $\ip12$ by Corollary~\ref{mpsii}, 
and is true in $\rV$. 
It follows by the inductive hypothesis
(applied in $\rV^+$)
that, in $\rV^+$, the $\ip11(p)$ set
$W^+$ of all pairs $\ang{x,T}$ such that $x\in X^+$,
$T\sq\nse$ is a \spt, and
$$
{\bod{T}\sq\fnp x1 \;\land\;
\text{the set $\bod{T}$ is \nibop}},
$$
--- satisfies $\pr {W^+}=X^+$. 
Therefore by the Shoenfield absoluteness theorem 
the set $W=W^+\cap\rV$ is $\ip11(p)$ 
and satisfies $\pr {W}=X$ in $\rV$.

Applying the Kondo --- Addison uniformization  
in $\rV^+$, we get a $\ip11(p)$ set 
$U^+\sq W^+$ which uniformizes $W^+$, in particular, 
$\pr {U^+}=\pr {W^+}=X^+$.
The corresponding set $U=U^+\cap\rV$ of type $\ip11(p)$ 
in $\rV$ then uniformizes $W$ 
and satisfies $\pr {U}=\pr {W}=X$ still by Shoenfield. 

Now, by the choice of the universe $\rV^+$, the 
uncountable $\ip11(p)$ set $U^+$ must contain a perfect 
subset $P^+\sq U^+$ of class $\ip01(q)$ for a parameter 
$q\in\rL[p]$, hence, $q\in\rV$. 
The according set $P=P^+\cap\rV$ is then a perfect 
subset of $U$ in $\rV$, and hence $X'=\pr P\sq X$ 
is a perfect set. 

Finally, if $x\in X'$ then let $T_x$ be the only element 
such that $\ang{x,T_x}\in P$. 
The map $x\mto T_x$ is Borel. 
On the other hand, still by the 
Shoenfield absoluteness, if $x\in X'$ then 
%$T_x$ is a \supt\ tree and 
$\bod{T_x}\sq\fn x1$ and the set $\bod{T_x}$ is \nibol.
\epF{Lemma}

We continue the proof of Theorem~\ref{fm}.

Let $X'\sq X$ and a Borel map\/ $x\mto T_x$ be as in the 
lemma.
If $x\in X'$ and $i=2,\dots,n+1$, 
then every \ddf iclass $\fn yi$ has at most 
one point common with the set $Y_x=\bod{T_x}$. 
Thus if $C$ is a \nism\ set then 
the intersection $C\cap Y_x$ is \bou\ 
and hence $C\cap Y_x$ 
is meager in $Y_x$. 

There is a Borel set $W\sq X'\ti\bn$ such that the 
collection of all cross-sections $\seq Wx$, $x\in X'$, 
is equal to 
the family of all countable unions of \ddf iclasses,
$i=2,\dots,n+1$, plus a \bou\ $\Fs$ set. 
(Note that \bou\ $\Fs$ sets is the same as \sik\ sets, and 
that every \bou\ set is a subset of a \bou\ $\Fs$ set.) 
Thus if $x\in X'$ then $\seq Wx\cap Y_x$ is meager in $Y_x$ 
by the above. 
Therefore, by a version of ``comeager uniformization'',  
there is a Borel map $f$ defined on $X'$ such 
that if $x\in X'$ then $f(x)\in Y_x\bez \seq Wx$.   
Clearly $f$ is $1-1$, hence 
the set $R=\ens{f(x)}{x\in X'}$ is Borel.

Moreover $R$ is \pit{\rF_1}\ by construction.
We assert that $R$ is \nila, in particular, 
not \bou! 

Indeed suppose otherwise. 
Then there is $x\in X'$ such that $R\sq\seq Wx$. 
But then $f(x)\in \seq Wx$, which contradicts the choice of 
$f$. 

Thus indeed $R$ is \nila. 
It follows by the inductive hypothesis 
that there exists a \nibol\ set $P\sq R$.
And $P$ is \pit{\rF_1}\ 
%pairwise \de inequivalent 
since so is $R$.
We conclude that $P$ is even \inbol, which leads 
to \ref{fm2} of the theorem.\vom

\epF{Theorem~\ref{fm}}

It is an interesting problem to figure out whether 
Theorem~\ref{fm} is true for a countable infinite family 
of \er s (as in Lemma~\ref{ico}). 
The inductive proof presented above is of little help, 
of course.

Another problem is to figure out whether the theorem still 
holds for $\fp11$ \er s, as the classical Silver dichotomy 
does. 
This is open even for the case of one $\fp11$ \er, since the 
background result, Corollary~\ref{pssC}, does not cover 
this case.

And finally we don't know whether Theorem~\ref{fm} can be 
strengthened to yield the existence of sets free 
(as in Corollary~\ref{pssC}) 
for a given (finite or countable) collection of Borel sets. 

It remains to note that Theorem~\ref{fm}
(in its relativized form)
implies the following theorem, perhaps, not
known previously in such a generality.

\bte
\lam{fmC}
Suppose that\/ $\rF_1,\dots,\rF_n$ are\/ Borel \er s
on a Polish space\/ $\pX$,
and\/ $A\sq\pX$ is a\/ $\fs11$ set.
Then either\/ $A$ is\/ \nsm, or
there exists an\/ \nbol\ set\/ $P\sq A$.\qed 
\ete

Yet the case $n=1$ is known in the form of the following
(not yet published) superperfect
dichotomy theorem of Zapletal:

\bte
\lam{zt}
If\/ $\rE$ be a Borel\/ \er\ on\/ $\bn$
and\/ $A\sq\bn$ is a\/ $\fs11$ set 
then either\/ $A$ is covered by 
countably many\/ \de classes and a \/ \bou\ set or 
there is a \spp\/ \pit\rE\ set\/ $P\sq A$. \qed
\ete

Theorem~\ref{zt} can be considered as 
a \lap{superperfect} version of Silver's dichotomy
(see \cite{sil} or \cite[10.1]{k-ams}), 
saying that if $\rE$ is a Borel \er\ then either 
the domain of $\rE$ is a union of  
countably many \ddd\rE classes or there is a 
perfect\/ \pis\rE\/ $Y\sq D$.

\punk{The case of $\is12$ sets: preliminaries}
\las{uu1}

In view of the counterexamples in Section~\ref{abo}, one can
expect 
that positive results for $\is12$ sets similar to Theorems 
\ref{mt'}, \ref{fm}, \ref{mt} should be expected in terms of 
\ddd\omi unions of compact sets. 
And indeed using a determinacy-style argument, 
Kechris proved in \cite
%[pp.\ 198--199]
{K}
that if $A\sq\bn$ is a $\is12$ set then 
(in a somewhat abridged form) 
one of the following two claims holds:
\ben
\tenu{{\rm(\Roman{enumi})}}
\itla{kes1} 
$A$ is \rit{\lbou}, 
in the sense that 
it is covered by the union of all sets\/ $[T]$,
where\/ $T\in\rL$ is a compact tree\snos
{$\rL$ is the constructible universe.} 
(hence not necessarily a countable union) 
--- or equivalently, 
for each $x\in A$ there is $y\in\bn\cap\rL$ with 
$x\doo y$, where $\doo$ is the
eventual domination order on $\bn$, 

\itla{kes2} 
there is a \sps\/ $P\sq A$.
%or there is a \spt\/ $\tau\in\rL$ such that\/
%$\bod\tau\sq A$.
\een

Our next goal is to generalise this result 
in the directions of Theorem~\ref{fm}. 
%(See Theorem~\ref{tk} below.) 
The logic of such a generalization forces us to change 
\sps s in \ref{kes2} by 
\nbol\ sets, where $\rF_1,\dots,\rF_n$ 
is a given collection of $\id11$ \er s. 
As for a corresponding change in \ref{kes1}, 
one would naturally look for a 
condition like: 
\bqu
for each $x\in A$, 
either there is $y\in\bn\cap\rL$ with $x\doo y$, 
or there is $j=1,\dots,n$ and an \lde\  
\dd{\rF_j}equivalence class containing $x$, 
\equ
whatever being \lde\  would mean.
The following example shows that the most elementary 
definition of \lde\  as 
``containing a constructible element'' fails.

\bex
\lam{uue1}
Let $\rF$ be the \er\ of equality of countable sets of 
reals. 
That is, its domain is the set $\bn^\om$ of all infinite 
sequences of reals, and for $x,y\in\bn^\om$, $x\rF y$ iff 
$\ran x=\ran y$. 
Let 
%us work in a \dd{\coll(\omi^\rL)}generic extension 
%$\rL[f]$ of $\rL$, where 
$f:\om\onto\omi^\rL$ is a generic 
collapse map. 
In $\rL[f]$, let $A$ be the $\is12$ set of all $x\in\bn^\om$ 
such that $\ran x$ (a set of reals) belongs to $\rL$ 
(but $x$ itself does not necessarily belong to $\rL$). 
Then, if $x\in A$ then the \df class $\fk x$ is not \bou, 
and the quotient $A/{\mathord{\rF}}$ 
(the set of all \df classes inside $A$) 
is uncountable in $\rL[f]$. 

We believe that there is no perfect (let alone \spp) 
\pit\rF set $P\sq A$ in $\rL[f]$, which is quite a safe 
conjecture in view of the results in \cite{gs}.
Yet to make the example self-contained let us add to 
$\rL[f]$ a set $C$ of $\aleph_3^\rL=\aleph_2^{\rL[f]}$ 
Cohen reals.
By a simple cardinality argument, there are no perfect 
\pit\rF\ sets $P\sq A$ in $\rL[f,C]$.
%Therefore by Theorem~\ref{tk}  
%the set $A$ is \lebo{\rF} in $\rL[f,C]$.

However, in $\rL[f,C]$, the quotient $A/{\mathord{\rF}}$ 
has uncountably many particular \df classes which are non-\bou\ 
and even non-\lbou\ in the sense of \ref{kes1} above, 
but contain no constructible elements. 
Thus $A$ neither contains an \dd\rF superperfect subset nor 
satisfies the condition that 
\rit{for each\/ $x\in A$, 
either there is\/ $y\in\bn\cap\rL$ with\/ $x\doo y$, 
or there is an\/ \df equivalence class containing\/ $x$ and 
containing a constructible element}.
%and not equal to Borel sets coded in $\rL$. 
%Therefore condition \ref{tk1} of Theorem~\ref{tk} cannot 
%be strengthened along the lines outlined just before 
%\ref{uue1}.
\qed
\eex 

Our model for \lde\ will be somewhat more 
complex than just ``containing a constructible element''.
In fact we'll consider two (connected) models, 
one being based on a certain uniform version of $\id11$, 
with ordinals as background parameters, the other one 
being based on Borel coding.
They are introduced in the following definitions.

\bdf
[coding ordinals]
\lam{kodo}
Let $\wo\sq\bn$ be the $\ip11$ set of all codes of countable 
(including finite) ordinals, and if $\xi<\omi$ then we define 
$\wo_\xi=\ens{w\in\wo}{w\text{ codes }\xi}$. 
If $w\in\wo_\xi$ then put $\abs w=\xi$.
\edf

\bdf
\lam{mod1}
A $\is12$ map $h:\bn\to\bn$ is \rit{absolutely total} 
if it remains total  
in any set-generic extension of the universe. 
In other words, it is required that there is a $\is12$ 
formula $\sg(\cdot,\cdot)$ such that 
$h=\ens{\ang{x,y}}{\sg(x,y)}$ and the sentence 
$\kaz x\:\sus y\:\sg(x,y)$ is forced by any set 
forcing. 
(Note that a total but not absolutely total map can be defined 
in $\rL$ by $h(x)$= the Goedel-least $w\in\wo$ such 
that $x$ appears at the \dd\xi th step of the Goedel 
construction, where $\xi=\abs w<\omi$ is the ordinal 
coded by $w$.)

%\lam{unxi}
Suppose that $\xi<\omi$. 
A set $X\sq\bn$ is \rit{\ess\/ $\is1n(\xi)$} if 
there is a $\is1n$ formula $\vpi(x,w)$ such that 
$X=\ens{x\in\bn}{\vpi(x,w)}$ for every $w\in\wo_\xi$. 
\Ess\ $\ip1n(\xi)$ sets are defined 
similarly, while an \ess\ $\id1n(\xi)$ set is any set 
both \ess\  $\is1n(\xi)$ and \ess\  $\ip1n(\xi)$.

A set $X$ is \rit{\ess\/ $\DD\xi$} if 
there is an absolutely total $\is12$ map $h$, 
%with $\wo_\xi\sq\dom h$, 
a $\is11$ formula $\chi(\cdot,\cdot)$, and 
a $\ip11$ formula $\chi'(\cdot,\cdot)$, such that 
if $w\in\wo_\xi$ then 
$X=\ens{x\in\bn}{\chi(x,h(w))}=\ens{x\in\bn}{\chi'(x,h(w))}$. 
\edf

Thus 
\ess\ $\DD\xi$ sets belong in between \ess\ $\id11(\xi)$ 
and \ess\ $\id12(\xi)$. 
%We'll show (see \ref{sepp} and \ref{sepq}) 
%that \ess\  $\DD\xi$ sets admit a  
%Borel coding with (not necessarily countable) codes in $\rL$. 
%
Each \ess\  $\DD\xi$ set $X$ is Borel, 
hence, it admits a Borel code. 
Moreover, if $X$ is \ess\  $\DD\xi$ via an absolutely total 
$\is12$ map $h$, and $w\in\wo_\xi$, 
then $X$ admits a Borel code in $\rL[w]$. 
We'll show (see \ref{sepp} and \ref{sepq}) that such a set 
$X$ admits a Borel code, even in $\rL$, 
in some generalized sense which allows 
uncountable Borel operations.  

\bdf
\lam{mod2}
Let $\olom$ be the class of all strings 
(finite sequences) 
of ordinals. 
If $s\in\olom$ and $\xi\in\Ord$ then $s\we\xi$ denotes 
the string $s$ extended by $\xi$.
If $s\in\olom$ then $\lh s$ is the length of $s$.
% 
%If $m<\lh s$ then $s\res m$ is the restricted string.
% 
$\La$ is the empty string. 
%; $\lh \La=0$ and $\La=s\res0$ for any $s\in\olom$.

A set $T\sq\olom$ is a \rit{tree} if $T\ne\pu$, and 
for any $s\in T$ and 
$m<\lh s$ we have $s\res m\in T$.
Then 
let $\sup T$ be the least ordinal $\la$ such that 
$T\sq\la^{<\om}$, and let $\mat T$ be the set of all 
\dd\sq maximal elements $s\in T$.

%Obviously $\La\in T$ for any tree $T$.

If a tree $T$ is \rit{well-founded} 
%iff it contains no infinite branches. In this case, 
then a \rit{rank function} $s\mto \rak Ts\in\Ord$ 
can be associated with $T$ so that 
$\rak Tt=\tsup_{t\we\xi\in T}(\rak T{t\we\xi}+1)$ 
(the least ordinal strictly bigger than 
all ordinals of the form $\rak T{t\we\xi}$, where 
$\xi\in\Ord$ and $t\we\xi\in T$) 
for each $t\in T$.
In particular $\rak Ts=0$ for any $s\in\mat T$. 

Let $\ra T=\rak T\La$ (the \rit{rank} of $T$).

%\bdf
%\lam{bc}
Let $\bc$ be the class of all 
\rit{generalized Borel codes} in $\rL$, that is, all pairs 
$c=\ang{T,d}=\ang{T_c,d_c}\in\rL$, 
where $T\sq\olom$ is a well-founded tree and 
$d\sq T\ti \nse$.
In this case, a set $[T,d,s]\sq\bn$ can be defined for each 
$s\in T$ by induction on $\rak Ts$ so that
$$
[T,d,s]= 
\left\{
\bay{rcl}
\bn\bez\bigcup_{\ang{s,u}\in d}\ibn u 
&\text{, whenever}& s\in\mat T\,;\\[1ex]
\bn\bez\bigcup_{s\we\xi\in T}[T,d,s\we\xi]
&\text{, whenever}& \rak Ts>0\,. 
\eay
\right.
$$
Recall that $\ibn u=\ens{a\in\bn}{u\su a}$ is a Baire interval.

Finally we put $[T,d]=[T,d,\La]$.
%\edf

%\bdf
%\lam{kro}
If $\ro<\omi$ then let $\bc_\ro\in\rL$ be the set of all codes 
$\ang{T,d}\in\bc$ such that $\ra T\le\ro$ and 
$\sup T\le \om_\ro^\rL$. 
(Not necessarily $\sup T<\omi$.)

%Note that the set of codes $\bc_\ro$ belongs to $\rL$. 

Accordingly let $[\bc_\ro]=\ens{[T,d]}{\ang{T,d}\in\bc_\ro}$.
\edf

If $\ang{T,d}\in\bc$ and $\sup T<\omi$ then $[T,d]$ is a 
Borel set in $\fp0{1+\ra T}$.

We underline that only \rit{constructible} codes are 
considered.

\punk{The case of $\is12$ sets: the result}
\las{uu2}

We'll prove the next theorem which generalizes the result 
of Kechris in \cite{K} sited above.
If $\rF$ is an \er\ on $\bn$ then let a \rit{\ds\df class} 
be any 
finite or countable union of \df equivalence classes. 

\bte
\lam{tk}
Let\/ $n<\om$, 
%$\rho<\omil$,
$\rF_1,\dots,\rF_n$ be\/ $\id11$ \er s on\/ $\bn,$
%simultaneously and\/ $\fp0{1+\ro}$,
and\/ $A\sq\bn$ be a\/ $\is12$ set.
Then one of the following\/ \ref{tk1}, \ref{tk2} 
holds$:$
\ben
\tenu{{\rm(\Roman{enumi})}}
\itla{tk1}
$A$ is\/ \lebo{\ans{\rF_1,\dots,\rF_n}}, 
%below\/ $\rho+2\;;$  
in the sense that for each $x\in A$:
\bit
\item[\rm --] 
either there is\/ $y\in\bn\cap\rL$ such that\/ $x\doo y$, 

\item[\rm --] 
or (non-exclusively) 
we have both\/ \ref{tk1a} and\/ \ref{tk1b}, where
\ben
\def\theenumii{{\rm\alph{enumii}}}
\def\labelenumii{{\rm(\theenumii)}}
\itla{tk1a} 
there is\/ $j=1,\dots,n$ and a \ds\dd{\rF_j}class $C$ 
which contains\/ $x$ and is\/ 
%both $\fp0{1+\ro}$ and 
\ess\  $\DD\xi$ for some\/ $\xi<\omi\;,$\vom

\itla{tk1b}
if, in addition, $\ro<\omil$ 
and all\/ $\rF_j$ are\/ $\fp0{1+\ro}$  
then there is\/ $j=1,\dots,n$ 
and an\/ \dd{\rF_j}class\/ $X$ in\/ $[\bc_\ro]$
which contains\/ $x\;;$
\een
\eit

\itla{tk2}
there exists an\/ \nbol\ set\/ $P\sq A$. 
%which is an\/ \tra{\rF_i} 
%for each\/ $i=1,\dots,n$.  
\een
\ete

A point of certain dissatisfaction is 
$\om^\rL_\ro$ as the measure of borelness 
in Definition~\ref{mod2} 
%and \ref{ef-2} 
and subsequently in \ref{tk1b} of Theorem~\ref{tk}.
Can it be reduced, to present the borelness involved by 
considerably narrower trees (of the same height)? 
Examples given in \cite{sa} and more resently in \cite{dsr}
allow to conjecture that the value $\omro$ cannot be 
reduced in any essential way.
A similar question can be addressed to 
the inequality 
$\omr1<\omi$ in the next remark.%

Suppose that $\omr1<\omi$. 
Then both $\bn\cap\rL$ and $\bc_\ro$ are countable sets, 
and hence the number of points $y$ and classes $X$ 
involved in \ref{tk1} of Theorem~\ref{tk}. 
Thus, assuming $\omr1<\omi$, condition \ref{tk1} 
of Theorem~\ref{tk} can be replaced by just the 
\sm{\ans{\rF_1,\dots,\rF_n}}ness of  $A$. 

\vyk{
The ``or'' option in 
%the definition of \lebo\cF\ sets 
%and subsequently the formulation of 
\ref{tk1} of Theorem~\ref{tk} 
leaves a certain sense of dissatisfaction as one would 
rather look for coverage by \dd{\rF_j}classes themselves 
than \ds\dd{\rF_j}classes. 
See Theorem~\ref{resT} and Corollary~\ref{resTC} on 
resolution of \ds\dd{\rF_j}classes into appropriately definable 
\dd{\rF_j}classes, 
in the context of \ref{tk1} of Theorem~\ref{tk}.
See also Remark~\ref{tkC} on a special form of 
Theorem~\ref{tk}.
}

The proof of Theorem~\ref{tk} will consist of two major 
parts. 
First of all, we prove, in this section, the version 
sans condition \ref{tk1b} in \ref{tk1} of the theorem. 
Then we prove, in Section~\ref{BC}, that 
\ref{tk1a} implies \ref{tk1b} in \ref{tk1}. 

\bpf[Theorem~\ref{tk} sans \ref{tk1b}]
We'll make use of  
Theorems \ref{ico+} and \ref{fm} in key arguments. 
To begin with, we reveal a certain uniformity in 
Theorem~\ref{ico+}\ref{ico+3}, which was not convenient 
to deal with in Section~\ref{dig}.

\bpro
\lam{+ico}
Under the conditions of Theorem~\ref{ico+}, 
if\/ $\xi<\omi$, 
all relations\/ $\rF_j$ are \ess\/ $\id11(\xi)$, and\/ 
$X$ is \ess\/ $\is11(\xi)$, then the sets\/ $Y_j$ and\/ 
$X_{n+1}$ satisfying\/ \ref{ico+1}, \ref{ico+2} can 
be chosen to be \ess\/ $\DD\xi$. 
%hence, uniform\/ $\id12(\xi)$. 
\epro
\bpf[Sketch]
We come back to the proof of Theorem~\ref{ico+}. 
The map $h(p)=\bap$ is defined in finitely many steps, 
such that each step is governed by a combination of $\is11$ 
and $\ip11$ formulas, so it is absolutely total $\is12$. 
\epf

%\bpf[Theorem~\ref{tk}]
In continuation, note that, by Kondo's uniformization, 
$A$ is the projection of a uniform
$\ip11$ set $B\sq\bn\ti 2^\dN$. 
% and hence we can assume from
%the beginning that $A$ itself is a $\ip11$ set.
Let $B=\bigcup_{\xi<\omi}B_\xi$ be an ordinary decomposition
of $B$ into pairwise disjoint Borel sets $B_\xi$
(called \rit{constituents}).
There is a $\is11$ formula $\ba(w,x,y)$ and a $\ip11$ 
formula $\ba'(w,x,y)$ such that 
\ben
\Aenu
\itla{13.1}
if $\xi<\omi$ and $w\in\wo_\xi$ 
then\\[1ex]
\phantom{1}\hfill$
B_\xi=\ens{\ang{x,y}}{\ba(w,x,y)}=
\ens{\ang{x,y}}{\ba'(w,x,y)}. 
%\quad\text{whenever $\xi<\omi$ and $w\in\wo_\xi$}
%\eqno(1)
$\hfill\phantom{1}
\een
We put $A_\xi=\pr B_\xi$; 
then $A=\bigcup_{\xi<\omi}A_\xi$ 
and all sets $A_\xi$ are Borel, and moreover $A_\xi$ is 
$\id11(w)$ whenever $w\in\wo_\xi$, because 
\ben
\Aenu
\atc
\itla{13.2}
if $\xi<\omi$ and $w\in\wo_\xi$ 
then\\[1ex]
\phantom{1}\hfill$
x\in A_\xi
\leqv  
\underbrace{\sus y\:\ba(w,x,y)}_{\al(w,x)}
\leqv  
\underbrace{\sus y\in\id11(x,w)\:\ba'(w,x,y)
}_{\al'(w,x)}\,,
$\hfill\phantom{1}\\[1ex]
where $\al(\cdot,\cdot)$ is a $\is11$ formula while 
$\al'(\cdot,\cdot)$ is a $\ip11$ formula.
\een

\rit{Case 1\/}: 
There is an ordinal $\xi<\omi$ such that $A_\xi$ is 
{\ubf not} \sm{\ans{\rF_1,\dots,\rF_n}}. 
Then we have 
\ref{tk2} of the theorem by Theorem~\ref{fm}.\vom 

\rit{Case 2\/}: 
All sets $A_\xi$ are \sm{\ans{\rF_1,\dots,\rF_n}}.
We claim that, 
\rit{under this assumption, if\/ $\xi<\omi$ 
then the set\/ $A_\xi$
is\/ \lebo{\ans{\rF_1,\dots,\rF_n}} 
%below\/ $\ro+2$, 
in the sense of\/ \ref{tk1} of Theorem~\ref{tk},
hence\/ \ref{tk1} of Theorem~\ref{tk} holds for\/ $A$.}

To prove the claim, fix $\xi<\omi$.

The set $A_\xi=\ens{x}{\al(w,x)}=\ens{x}{\al'(w,x)}$
(for any $w\in\wo_\xi$) is \ess\  $\id11(\xi)$. 
The relations $\rF_j$ are just $\id11$.
By Theorem~\ref{ico+} and Proposition~\ref{+ico} 
there exist Borel sets 
$Y_1,\dots,Y_n,X_{n+1}\sq\bn$ satisfying 
\ref{ico+1}, \ref{ico+2} 
%, \ref{ico+3} 
of Theorem~\ref{ico+} for $X=A_\xi$, and \ess\  $\DD\xi$. 
Thus $X_{n+1}$ is \bou,  
each set\/ $Y_j$ is a \ds\dd{\rF_j}class, and 
$A_\xi\sq Y_1\cup \dots\cup Y_n\cup X_{n+1}$.

Now the claim and Theorem~\ref{tk} immediately follow 
from: 

\ble
\lam{tkl1}
If\/ $y\in X_{n+1}$ then there is a real\/ 
$a\in\rL$ such that\/ $y\doo a$. 
\ele

\vyk{
\ble
\lam{tkl2}
If\/ $1\le j\le n$ then the set\/ $Y_j$ belongs to\/ 
$[\bc_{\ro+2}]$.
%is\/ \lebo{\ans{\rF_1,\dots,\rF_n}}  below\/ $1+\rho+2$. 
%absolutely\/ $\od$ over\/ $\rL$. 
\ele
}

\bpf[Lemma]
Let $\rV$ be the whole set universe in which we prove the
lemma.
Thus $\xi<\omi^\rV$ but not necessarily $\xi<\omi^\rL$.\vom

{\it Case L1}: $\xi<\omi^\rL$.
%Then $(A_\xi)^\rL$ is a Borel set coded in $\rL$.
We assert that the set $(X_{n+1})^\rL=X_{n+1}\cap\rL$ 
is \bou\ in $\rL$.
Indeed otherwise by Theorem~\ref{mt'}
(relativized)
there is a \spt\ $T\in\rL$ such that
$[T]\sq(X_{n+1})^\rL$ in $\rL$. 
Then by Shoenfield %we have
$[T]\sq X_{n+1}$ in the universe,
contrary to the \bou ness of $X_{n+1}$.

Thus $(X_{n+1})^\rL$ is \bou\ in $\rL$,
so there is a real $a\in\rL$ such that
if $y\in (X_{n+1})^\rL$ then $y\doo a$.
Then, again by Shoenfield, it is true in the universe that
if $y\in X_{n+1}$ then $y\doo a$, as required.\vom

{\it Case L2}: $\xi\ge\omi^\rL$.
Recall that the set $X_{n+1}$ is \ess\  $\DD\xi$, via a 
certain absolutely total $\is12$ map $h(w)=\baw$ and formulas 
$\chi(\cdot,\cdot)$ of type $\is11$ and 
$\chi'(\cdot,\cdot)$ of type $\ip11$, 
as in Definition~\ref{mod2}, so that 
\ben
\Aenu
\atc\atc
\itla{13.3}
if $w\in\wo_\xi$ then
$X_{n+1}=\ens{x}{\chi(\baw,x)}=\ens{x}{\chi'(\baw,x)}$.
\een
Let $f,g\in\xi^\om$ be collapse functions generic over $\rV$,
such that the pair $\ang{f,g}$ is
generic over $\rV$ as well.
Then $\xi<\omi^{\rL[f]}$ and $\xi<\omi^{\rL[g]}$.
For $X_{n+1}$ being \bou\ is a $\is12$ formula 
(make use of \ref{13.3} and Corollary~\ref{mpsii}), 
hence by Shoenfield 
it is true in $\rL[f]$ that the set 
$$
(X_{n+1})^{\rL[f]}=\ens{x\in\rL[f]}{\chi(\baw,x)} 
=\ens{x\in\rL[f]}{\chi'(\baw,x)}
$$
($w\in\wo_\xi\cap\rL[f]$ is arbitrary) is \bou.
Thus there is a real $x\in\bn\cap\rL[f]$ such that
$y\doo x$ for all $y\in (X_{n+1})^{\rL[f]}$.
Then again by the Shoenfield absoluteness 
$y\doo x$ holds even for all
$y\in (X_{n+1})^{\rL[f,g]}$.

In particular if $y\in (X_{n+1})^{\rL[g]}$ then $y\doo x$.

On the other hand, as $f,g$ are mutually generic, one can
show that if $x\in\bn\cap\rL[f]$,
$y\in\bn\cap\rL[g]$, and $y\doo x$ then there is a real
$a\in\bn\cap\rL$ with $y\doo a\doo x$.
We conclude that if $y\in (X_{n+1})^{\rL[g]}$ then there is
a real $a\in\rL$ such that $y\doo a$.
Now, as $\omi^\rL\le\xi<\omi^{\rL[g]}$, 
there exists a sequence $\sis{a_n}{n<\om}\in \rL[g]$ 
such that $\bn\cap\rL=\ens{a_n}{n<\om}$.
Therefore we have
$$
\kaz y\in (X_{n+1})^{\rL[g]}\:\sus n\:(y\doo a_n)
$$
in $\rL[g]$, and hence by Shoenfield 
$\kaz y\in (X_{n+1})^{\rV[g]}\:\sus n\:(y\doo a_n)$
in $\rV[g]$.
But as $\xi<\omi^\rV$, there is a sequence
$\sis{b_n}{n<\om}\in \rV$ such that  
$\bn\cap\rL=\ens{b_n}{n<\om}$.
Then we have
$\kaz y\in (X_{n+1})^{\rV[g]}\:\sus n\:(y\doo b_n)$
in $\rV[g]$, and by Shoenfield 
$\kaz y\in X_{n+1}\:\sus n\:(y\doo b_n)$ in $\rV$, and so on.
\epF{Lemma}

\epF{Theorem~\ref{tk} sans \ref{tk1b}}

\vyk{
\bpf[Lemma~\ref{tkl2}]
\vyk{Similarly to the above, if $j=1,\dots,n$ then there exists 
a $\is11$ formula $\et_j(\cdot,\cdot)$ 
and a $\ip11$ formula $\et'_j(\cdot,\cdot)$
such that 
\ben
\fenu
\atc\atc\atc
\itla{13.4}
if $w\in\wo_\xi$ then 
$
Y_{j}=\ens{x}{\et_j(\baw,x)}=\ens{x}{\et_j'(\baw,x)}
$
\een 
%whenever $w\in\wo_\xi$. 
Combining these formulas with the formulas $\al,\al'$ 
as above, and the formulas defining the map 
$p\mto\bap$, we easily obtain $\is12$ formulas 
$\psi_j(\cdot,\cdot)$ such that 
$Y_{j}=\ens{x}{\psi_j(w,x)}$ whenever $w\in\wo_\xi$. 
It follows that $Y_{j}$ is $\ish1(\xi)$, 
where $\hc$ is the set of all hereditarily countable sets. 
By a similar argument, $Y_{j}$ is $\iph1(\xi)$ as well, 
so that $Y_{j}$ is $\idh1(\xi)$.
}%
Each set $Y_j$ is \ess\  $\DD\xi$, see above, and hence 
\ess\  $\id12(\xi)$ and $\idh1(\xi)$.
Moreover $Y_j$ is a \ds\dd{\rF_j}class, hence, 
a set in $\fs0{1+\ro+1}$, while its complement 
$Z_j=\bn\bez Y_j$ is accordingly a set in 
$\idh1\cap\fp0{1+\rho+1}$.
It follows that $Z_j\in [\bc_{\ro+1}]$ by 
Proposition~\ref{sepp}, and hence $Y_j\in[\bc_{\ro+2}]$, 
as required.
\vyk{
Pick any $w\in\wo_\xi$.
As above, by \ref{ico+3} of Theorem~\ref{ico+} 
the set $Y_{j}$ is $\is11(\baw)$, and moreover, 
$$
Y_{j}=\ens{x}{\eta(\baw,x)}=\ens{x}{\eta'(\baw,x)}.\eqno(1)
$$ 
where $\eta(\cdot,\cdot)$ is a $\is11$ formula and 
$\eta'(\cdot,\cdot)$ is a $\ip11$ formula. 
We claim that 
$$
\vpi(\xi,x)\;:=\;
x\in\bn\land \xi<\omi\land \sus w\in\wo_\xi\:\eta(\baw,x)
$$
is a the formula which proves the lemma.
Assume that $P\in\rL$ is a forcing notion, 
a set $G\sq P$ is \dd Pgeneric 
over the basic set universe $\rV$, and $\xi<\omi^{\rL[G]}$. 
Pick any $v\in\wo_\xi\cap\rL[G]$. 
We have to prove that the set 
$$
Y'_j =
\ens{y\in\rL[G]}{\vpi(\xi,y)\text{ in }\rL[G]}= 
\ens{y\in\rL[G]}{\eta(\bav,y)\text{ in }\rL[G]}
\eqno(2)
$$ 
satisfies $Y_j\rF^+_j Y'_j$ in $\rV[G]$, where $\rF^+_j$ 
is the natural \dd{\rV[G]}extension of $\rF$.
(Then being a \ds\dd{\rF_j^+}class in $\rL[G]$ follows by   
Shoenfield.) 

%We can \noo\ assume that $G$ is \dd Pgeneric over the whole 
%universe $\rV$. 
Consider the bigger universe $\rV[G]$ and the 
\dd{\rV[G]}extension 
$$
Y^+_j =
\ens{y\in\rV[G]}{\eta(\bav,y)\text{ in }\rV[G]} 
\eqno(3)
$$
of $Y'_j$. 
As parameters $w,v$ belong to the same set $\wo_\xi$ 
in $\rV[G]$, 
we have 
$Y^+_j = \ens{y}{\eta(\baw,y)}$ in $\rV[G]$
as well by (3). 
Now let, in $\rV$, $S=\ens{a_k}{k<\om}\sq Y_j$ be a full 
transversal, that is, $S$ meets every \dd{\rF_j}class in 
$Y_j$ in exactly one point. 
By the Shoenfield absoluteness, 
this is true in $\rV[G]$, too, and hence 
$Y_j\rF_j^+ Y^+_j$.
By a similar argument, we have $Y'_j\rF_j^+ Y^+_j$ as well.
It follows that $Y_j\rF_j^+ Y'_j$, as required.
}%
\epF{Lemma~\ref{tkl2}}
}

\punk{The case of $\is12$ sets: proof of the second part}
\las{BC}

To demonstrate that the abridged version of Theorem~\ref{tk} 
implies the full version, it suffices to prove the 
following theorem.

\bte
\lam{resT}
Assume that, in the ground set universe $\rV$, 
\ben
\fenu
\itla{resF}\msur
\imar{resF}%
$\ro<\omil$, $\xi<\omi$, $\rE$ is an\/ \er\ on\/ $\bn$ 
in\/ $\id11\cap \fp0{1+\ro}$,  
$\pu\ne C\sq\bn$ is a\/ \ds\de class and a set \ess\/ $\DD\xi$. 
\een
Then each\/ \de class\/ $X\sq C$ 
is a set in\/ $[\bc_\ro]$.
\ete 

Any \ess\  $\DD\xi$ set is \ess\  $\id12(\xi)$, and 
hence $\idh1(\xi)$.
(Recall that $\hc$ is the set of all hereditarily 
countable sets.)
This simple fact 
will allow us to make use of the following result, 
explicitly proved 
in \cite{matsb} (Lemma~4) 
on the base of ideas and technique developed in 
\cite{ster',ster}.

\bpro
\lam{sepp}
Let\/ $X,Y\sq\bn$ are two disjoint sets in\/ 
$\ish1(\omi)$ 
{\rm(that is, $\ish1$ with any finite number of parameters 
in $\omi$)}.
Suppose that\/ $\ro<\omil$ and\/ $X$ is\/  
\dd{\fp0{1+\ro}}separable from\/ $Y$. 
Then there is a separating set in\/ 
$[\bc_\ro]$. 

In particular if\/ $Z\sq\bn$ is a set in\/ 
$\idh1\cap\fp0{1+\rho}$ then\/ $Z\in [\bc_\ro]$.\qed 
\epro

For instance, if $\ro=0$, so that $\fp0{1+\ro}=$ closed 
sets, then the result takes the form: any closed $\idh1$ 
set $Z\sq\bn$ has a code in the set\pagebreak[0]
$$
\bc_0=\ens{\ang{T,d}\in\bc}
{\ra T=0\text{ (hence just $T=\ans\La$) }\land\;\sup T\le \om},
$$
but this can be easily established directly.

Thus 
%if $\xi<\omi$ then 
sets \ess\  $\DD\xi$, $\xi<\omi$, 
even those \ess\ $\id12(\xi)$, 
admit a straight Borel coding by 
(not necessarily \rit{countable}) 
codes in $\rL$. 

\bpf[Theorem~\ref{resT}]
Assume that $\ro\yi\xi\yi\rE\yi C$ 
are as in \ref{resF} above.
Then $C$ is $\fs0{1+\ro+1}$, 
therefore by Lemma~\ref{sepq} $C$ belongs to 
$[\bc_{\ro+2}]$ in $\rV$. 
We'll show now that an appropriate coding can be chosen in 
{\ubf absolute} manner.

\vyk{
\bdf
\lam{add}
Let $\xi,\ro<\omi$. 
An \ess\ $\DD\xi$ set $X\sq\bn$ 
\rit{absolutely belongs to\/ $[\bc_\ro]$} 
if there is a code $\ang{T,d}\in\bc_\ro$ such that we have 
$X^{\rV[G]}=[T,d]$ in any set generic extension
${\rV[G]}$ of the universe $\rV$.

Note that then by Shoenfield the equality 
$X^{\rL[G]}=[T,d]$ also holds in any generic
extension ${\rL[G]}$ of $\rL$ such that $\xi<\omi^{\rL[G]}$.
\edf
}

\bre
\lam{xt}
Suppose that our set $C$ is \ess\ $\DD\xi$, 
via an absolutely total $\is12$ map $h$ and   
formulas $\chi\yd\chi'$ as in 
Definition~\ref{mod1}.
Then the following   
is true in the ground universe $\rV$:
\ben
\fenu
\atc
\itla{12c1} 
if $v,w\in\wo_\xi$ and $x\in\bn$ then\\[0.5ex]
\phantom{,}
\hfill
$\chi(x,h(v))\eqv\chi(x,h(w))\eqv\chi'(x,h(v))\eqv\chi'(x,h(w))$.
\hfill\phantom{,}
\een
If we eliminate $h$ by a formula $\sg$ as in 
Definition~\ref{mod1} 
then \ref{12c1} becomes a $\ip12$ sentence. 
Therefore \ref{12c1} is true in any
extension $\rV[G]$ of $\rV$ by Shoenfield, and
moreover, in any generic extension $\rL[G]$ of $\rL$
such that $\xi<\omi^{\rL[G]}$.
This allows us to unambiguously define extensions
$h^{\rV[G]}$ of $h$ (a total map)
and $C^{\rV[G]}$ of $C$ to $\rV[G]$, using the 
same formulas, so that $C^{\rV[G]}$ is an
\ess\ $\DD\xi$ set in $\rV[G]$ still
via $h^{\rV[G]}\yi\chi\yi\chi'$.
Then, assuming $\xi<\omi^{\rL[G]}$,
we define associated restrictions 
$h^{\rL[G]}=h^{\rV[G]}\cap\rL[G]$  
and $C^{\rL[G]}=C^{\rV[G]}\cap\rL[G]$ to  
$\rL[G]$, so that $C^{\rL[G]}$ is  
\ess\ $\DD\xi$ in $\rL[G]$  
via $h^{\rL[G]}\yi\chi\yi\chi'$, too.

And as $\rE$ is a $\id11$ \er\ in $\rV$, then, even easier, 
we define an extension $\EVG$ of $\rE$ to $\rV[G]$, using the 
same formulas which define $\rE$, so that $\EVG$ is a
$\id11$ \er\ in $\rV[G]$ by Shoenfield, and then define 
$\ELG=\EVG\cap\rL[G]$ 
(a $\id11$ \er\ in $\rL[G]$).
\ere

\ble
\lam{sepq}
%Suppose that\/ $\xi<\omi$, $\ro<\omil$, and a set\/ 
%$X\sq\bn$ is \ess\/ $\DD\xi$. 
%Then\/ 
$C$ {\bfit absolutely\/} belongs to\/ $[\bc_\rod]$, in the 
sense that there is a code\/ $\ang{T,d}\in\bc_\rod$ 
such that we have\/ 
$C^{\rV[G]}=[T,d]$ in any set generic extension
${\rV[G]}$ of the universe $\rV$.
\ele

Note that then by Shoenfield the equality 
$C^{\rL[G]}=[T,d]$ also holds in any generic
extension ${\rL[G]}$ of $\rL$ such that $\xi<\omi^{\rL[G]}$.

\bpf[Lemma]
Let a map $f:\om\na\om^\rL_{\ro+1}$ be collapse generic 
over $\rV$.
Let $C^{\rV[f]}\in\rV[f]$ be the extension of $C$ to 
$\rV[f]$, as above. 
Then $C^{\rV[f]}$ is \ess\  $\DD\xi$ in $\rV[f]$, and hence
by Proposition~\ref{sepp} there is a code 
$\ang{T,d}\in\bc_\ro$ 
such that $C^{\rV[f]}=[T,d]$ in $\rV[f]$.
To prove, that this code witnesses that 
$C$ absolutely belongs to\/ $[\bc_\ro]$, 
consider any generic extension $\rV[G]$.
It can be assumed that $G$ is generic even over $\rV[f]$. 

Let $C^{\rV[G]}\yi C^{\rV[f,G]}$ be the extensions of $C$ 
(a set in $\rV$) 
to resp.\ $\rV[G]\yi\rV[f,G]$ (see Remark~\ref{xt}). 
The code $\ang{T,d}$ is countable in 
$\rV[f]$ and in $\rV[f,G]$ by the choice of $f$. 
Therefore the equality $C^{\rV[f]}=[T,d]$ can be expressed by a 
Shoenfield-absolute formula. 
We conclude that $C^{\rV[f,G]}=[T,d]$ holds in $\rV[f,G]$, 
hence $C^{\rV[G]}=[T,d]$ is true in $\rV[G]$ as well as 
easily $C^{\rV[G]}= C^{\rV[f,G]}\cap\rV[G]$ and 
$[T,d]^{\rV[G]}= [T,d]^{\rV[f,G]}\cap\rV[G]$.
\epF{Lemma}

It follows that there is a code $\ang{T_0,d_0}\in\bc_\rod$ 
such that $[T_0,d_0]=C^{\rV[G]}$ in any extension $\rV[G]$ 
of $\rV$, and hence we obtain by Shoenfield:

\bcor
\lam{resG} 
In any set-generic extension\/ $\rV[G]$ of\/ $\rV$, 
$[T_0,d_0]$ is a\/ \ds\dd{\rE^{\rV[G]}}class   
containing only those\/ \dd{\rE^{\rV[G]}}classes   
presented in\/ $[T_0,d_0]\cap \rV$.\qed 
\ecor

\vyk{
\bpf 
Apply the Shoenfield absoluteness theorem. 
\epf
}
%and does not depend on the parameter $\xi$! 

\vyk{

\punk{The case of $\is12$ sets: resolution of \ds classes}
\las{res}

Here our goal will be to resolve 
\ds classes, as in the ``or'' option of 
%Definition~\ref{ef-2} and implicitly in 
\ref{tk1} of Theorem~\ref{tk}, into 
countable unions of single ``\dd\rL definable'' equivalence 
classes. 

We are going to prove the next theorem in this section. 

}

\vyk{
\bcor
\lam{resTC}
Under the requirements of Theorem~\ref{tk}, if $\ro<\omil$ 
and the relations $\rF_j$ all belong to $\fp0{1+\ro}$ then
\ref{tk1} of Theorem~\ref{tk} can be 
reformulated as follows$:$
\ben
\Renu
\itla{tk1r}
$A$ is\/ \lebo{\ans{\rF_1,\dots,\rF_n}}, 
%below\/ $\rho+2\;;$  
in the sense that for each $x\in A$:
\bit
\item[\rm --] 
either there is\/ $y\in\bn\cap\rL$ such that\/ $x\doo y$, 

\item[\rm --] 
or there is an index\/ $j=1,\dots,n$ and 
a\/ \dd{\rF_j}class\/ $X$ 
which contains\/ $x$ and belongs to\/ $[\bc_\ro]$.\qed
\eit
\een
\ecor 
}

We continue with a few definitions.

If $\ang{T,d}\yd \ang{T',d'}\in\bc$ 
then let $\ang{T,d}\lec\ang{T',d'}$ mean that 
$[T,d]\sq[T',d']$ in any set generic 
extension $\rL[G]$ of $\rL$.
Then, using appropriate collapse extensions, we 
conclude by Shoenfield, that 
$[T,d]\sq[T',d']$ also holds in any set generic 
extension $\rV[G]$ of the ground universe $\rV$, including 
$\rV$ itself.

Say that a code $\ang{T,d}\in \bc$ is \rit{\ene} 
if $[T,d]\ne\pu$ in at least one set-generic extension 
of $\rL$. 
By Shoenfield, this is equivalent to $[T,d]\ne\pu$ 
in some/any extension $\rL[G]$ with 
$\sup T<\omi^{\rL[G]}$.

%\bdf
%\lam{forP}
%If $\la\in \Ord$ then 
Let $\dpl\in\rL$ be the forcing 
notion which consists of all ``essentially non-empty'' 
codes $\ang{T,d}\in\bc$ such that 
$\ang{T,d}\lec\ang{T_0,d_0}$ and $\sup T\le\omr2$. 
We order $\dpl$ by $\lec$, and $\ang{T,d}\lec\ang{T',d'}$ 
is understood as $\ang{T,d}$ being a stronger 
forcing condition.
%\edf

In particular 
%if $\la\ge\om^{\rL}_\rod$ then 
condition $\ang{T_0,d_0}$ itself, as in 
Corollary~\ref{resG}, belongs 
to $\dpl$.

\ble
\lam{resL}
%If\/ $\la\ge\om^{\rL}_\rod$ then\/ 
$\dpl$ forces a real over\/ $\rL$, so that 
if a set\/ $G\sq\dpl$ is generic over\/ $\rL$ 
then the intersection\/ $\bigcap_{\ang{T,d}\in G}[T,d]$ 
contains a single real in\/ $\rL[G]$. 
\ele
\bpf
If $u\in\nse$ is a string of length $n=\lh u$ then let 
$T^u=\ans\La$ and let 
$d^u$ consist of all pairs $\ang{\La,v}$ such that 
$v\in\nse$, $v\ne u$, $\lh v=n$. 
Then $\ang{T^u,d^u}\in\dpl$ and 
$[T^u,d^u]=\ibn u=\ens{a\in\bn}{u\su a}$.  
By the genericity, for any $n$ there is a inuque 
$u=u[n]\in\nse$ such that $\lh {u[n]}=n$ and 
$\ang{T^{u[n]},d^{u[n]}}\in G$, and in addition 
$u[n]\su u[m]$ whenever $n<m$. 
It follows that there is a real $x_G=\bigcup_nu[n]\in\rL[G]$ 
such that $x_G\res n=u[n]$, and hence 
$x_G\in [T^{u[n]},d^{u[n]}]$, $\kaz n$.
We claim that 
\rit{if\/ $\ang{T,d}\in\dpl$ then\/ 
$\ang{T,d}\in G$ iff\/ $x_G\in [T,d]$ in\/ $\rL[G]$}; 
this obviously proves the lemma.

We prove the claim by induction on the rank $\ra T$.

Suppose that $\ra T=0$, so that $T=\ans\La$,  
$d\sq\ans\La\ti\nse$, and 
$[T,d]=\bn\bez\bigcup_{v\in U}\ibn v$, where 
$U=\ens{v\in\nse}{\ang{\La,v}\in d}$.
We assert that 
\ben
\aenu
\itla{resL1}
any $\ang{T',d'}\in\dpl$ is compatible, in $\dpl$, 
either with 
$\ang{T,d}$ or with one of the codes $\ang{T^v,d^v}$, 
where $v\in U$ --- therefore either $\ang{T,d}$ or 
one of the codes $\ang{T^v,d^v}$, $v\in U$, belongs to $G$.
\een
Indeed we have $[T,d]=\bn\bez\bigcup_{v\in U}[T^v,d^v]$ in 
any universe. 

With \ref{resL1} in hands, if $v\in U$ and 
$\ang{T^v,d^v}\in G$ then on the one hand $\ang{T,d}\nin G$ 
by \ref{resL1}, and on the other hand, 
obviously $v=u[n]$, where $n=\lh v$, so that 
$x_G\in [T^v,d^v]$ and $x_G\nin [T,d]$. 
Conversely, if there is no $v\in U$ with $\ang{T^v,d^v}\in G$ 
then on the one hand $\ang{T,d}\in G$ by \ref{resL1}, 
and on the other hand, $x_G\nin \bigcup_{v\in U}[T^v,d^v]$, 
so that $x_G\in [T,d]$. 

To carry out the step, suppose that $\ra T>0$. 
Let $\Xi=\ens{\xi}{\ang\xi\in T}$ 
(where $\ang\xi$ is a one-term string). 
If $\xi\in\Xi$ then let 
$$
T^\xi=\ens{s\in\olom}{\xi\we s\in T}
\qand
d^\xi=\ens{\ang{s,v}}{\ang{\xi\we s,v}\in d}\,.
$$
Thus each $\ang{T^\xi,d^\xi}$ is a code in $\dpl$, 
$\ra{T^\xi}<\ra T$, and 
$[T,d]=\bn\bez\bigcup_{\xi\in\Xi}[T^\xi,d^\xi]$ in any 
universe containing $\ang{T,d}$. 
Similarly to  \ref{resL1} above, we have
\ben
\aenu
\atc
\itla{resL2}
any $\ang{T',d'}\in\dpl$ is compatible, in $\dpl$, 
either with 
$\ang{T,d}$ or with one of the codes $\ang{T^\xi,d^\xi}$, 
where $\xi\in \Xi$ --- therefore either $\ang{T,d}$ or 
one of the codes $\ang{T^\xi,d^\xi}$, $\xi\in \Xi$, 
belongs to $G$.
\een
Now, if $\xi\in\Xi$ and 
$\ang{T^\xi,d^\xi}\in G$ then 
on the one hand $\ang{T,d}\nin G$ 
by \ref{resL2}, and on the other hand, 
$x_G\in [T^\xi,d^\xi]$ by the inductive hypothesis, 
and hence $x_G\nin [T,d]$. 
Conversely, if there is no $\xi\in\Xi$ with 
$\ang{T^\xi,d^\xi}\in G$ 
then on the one hand $\ang{T,d}\in G$ by \ref{resL2}, 
and on the other hand, 
$x_G\nin \bigcup_{\xi\in\Xi}[T^\xi,d^\xi]$, 
by the inductive hypothesis, so that $x_G\in [T,d]$.
\epf

Reals of the form $x_G=$ 
the only element of $\bigcap_{\ang{T,d}\in G}[T,d]$ 
in $\rL[G]$, 
where $G\sq\dpl$ is \dd{\dpl}generic, \eg, over $\rV$, 
will be called 
\rit{\dd{\dpl}generic over\/ $\rV$}, too.
Let $\nx$ be a canonical \dd{\dpl}name for $x_G$. 
Let $\lx$, $\px$ be canonical 
\dd{(\dpl\ti\dpl)}names for 
the left and the right copies of $x_G$. 

Let $\uE$ be a canonical \dd{\dpl}name 
for the extension $\EVG$ or $\ELG$ of $\rE$ to any class like 
$\rL[G]$ or $\rV[G]$, $G$ being generic.

\bdf
\lam{sta}
A code $\ang{T,d}\in\dpl$ is \rit{stable} if 
condition 
$\qar{\ang{T,d}}{\ang{T,d}}$ \ 
\dd{(\dpl\ti\dpl)}forces, 
over $\rL$, that $\lx \uE \px$.
\edf

\ble
\lam{stL1}
If\/ 
%$\la\ge\om^{\rL}_\rod$ and\/ 
$\ang{T,d}\in\dpl$ is stable then there is an element\/ 
$y\in C=[T_0,d_0]\cap\rV$ 
such that\/ $\ang{T,d}$ \dd{\dpl}forces, 
over\/ $\rV$, that\/ $\nx \uE y$.
\ele

\vyk{
We identify $\rF$ (a Borel relation in $\rV$) 
with its extension $\rF^{\rV[G]}$ to any bigger universe 
$\rV[G]$ and with subsequent restriction $\rF^{\rL[G]}$
to $\rL[G]\sq\rV[G]$.
}

\bpf
By Corollary~\ref{resG} the contrary assumption leads to a pair 
of conditions 
$\ang{T',d'}\lec\ang{T,d}$ and $\ang{T'',d''}\lec\ang{T,d}$
in ${\dpl}$ and elements $y',y''\in [T_0,d_0]\cap\rV$ such that 
\begin{center}
$\ang{T',d'}$ \dd{\dpl}forces $\nx \uE y'$, \ and \
$\ang{T'',d''}$ \dd{\dpl}forces $\nx \uE y''$ \ --- \
over $\rV$,  
\end{center}
and $y'\nE y''$.
To get a contradiction consider a set $G'\ti G''$, 
\dd{({\dpl}\ti{\dpl})}generic over $\rV$, and containing 
condition $\qar{\ang{T',d'}}{\ang{T'',d''}}$. 
Then, on the one hand, the generic reals $x_{G'}$ and 
$x_{G''}$ satisfy $x_{G'}\EVGp y'$ and $x_{G''}\EVGpp y''$, 
but on the other hand, $x_{G'}\EVGq x_{G''}$ holds 
by stability.  
Therefore $y'\rE y''$, which contradicts to the choice 
of $y'\yi y''$.
\epf

\vyk{
\bcor
\lam{stC1}
If\/ $\la\ge\om^{\rL}_\ro$ and 
conditions\/ $\ang{T,d}$ and\/ $\ang{T',d'}$ in\/ ${\dpl}$ 
are stable then condition $\qar{\ang{T,d}}{\ang{T',d'}}$ 
either\/ \dd{({\dpl}\ti{\dpl})}forces, over $\rL$, that\/ 
$\lx \uE \px$, 
or\/ \dd{({\dpl}\ti{\dpl})}forces, over $\rL$, that\/ 
$\lx \nuE \px$.\qed
\ecor
}

\ble
\lam{stL2}
%If\/ $\la\ge\om^{\rL}_\rod$ then t
The set of all stable 
conditions\/ $\ang{T,d}\in{\dpl}$ is 
dense in\/ ${\dpl}$.
\ele
\bpf
%Let $\mu=\om^\rL_{\ro+3}<\nu=\om^\rL_{\ro+4}$ 
%(the least \dd\rL cardinals above $\om^\rL_{\ro+2}$). 
By definition $\car{\dpl}=\omr3$ and $\car{\pws{\dpl}}=\omr4$
in $\rL$.
Consider an extension $\rV[g]$ by a collapse-generic 
map $g:\om\na\omr4$. 
Then, in $\rV[g]$, there is an enumeration $\sis{D_n}{n<\om}$ 
of all dense sets $D\sq\dpl\ti\dpl\yt D\in\rL$.

Now suppose towards the contrary that 
$\ang{T^\ast,d^\ast}\in{\dpl}$ 
and there is no stable $\ang{T,d}\lec\ang{T^\ast,d^\ast}$ 
in ${\dpl}$.
Then for any condition $\ang{T,d}\lec\ang{T^\ast,d^\ast}$ 
there are stronger conditions 
$\ang{T',d'}\lec\ang{T,d}$ and 
$\ang{T'',d''}\lec\ang{T,d}$
such that $\qar{\ang{T',d'}}{\ang{T'',d''}}$
\dd{({\dpl}\ti{\dpl})}forces $\neg\;\lx \uE \px$ 
over $\rL$.
This allows to define, in $\rV[g]$, a family 
$\sis{\ang{T(u),d(u)}}{u\in\bse}$ of conditions in $\dpl$  
satisfying $\ang{T(\La),d(\La)}=\ang{T^\ast,d^\ast}$,
and in addition
\ben
\renu
\item\msur
$\ang{T(u\we i),d(u\we i)}\lec\ang{T(u),d(u)}$ for each 
$i=0,1$ and $u\in\nse$, 

\item 
if $u\ne v$ are of length $n+1$ 
then $\qar{\ang{T(u),d(u)}}{\ang{T(v),d(v)}}\in D_n$, 

\item 
if $u\in\bse$ then the condition 
$\qar{\ang{T(u\we 0),d(u\we 0)}}{T(u\we 1),d(u\we 1)}$ \ 
\dd{({\dpl}\ti{\dpl})}forces 
$\neg\;\lx \uE \px$ over $\rL$.
\een
Then, in $\rV[g]$, if $a\in\dn$ then the intersection 
$\bigcap_n[T(a\res n),d(a\res n)]$ contains a single point 
$x_a\in [T^\ast,d^\ast]$ by Lemma~\ref{resL},  
and if $a\ne b$ then $\neg\;(x_a\EVg x_b)$. 
But by construction 
$[T^\ast,d^\ast]\sq [T_0,d_0]$ in $\rV[g]$, 
so that $[T_0,d_0]$ contains uncountably many \dd\EVg classes 
in $\rV[g]$ --- a contradiction to 
Corollary~\ref{resG}.
\epf

Let $H$ be the set of all codes $\ang{T,d}\in\bc_\ro$ such 
that the \dd{\omr4}collapse forcing notion
$\coll(\omr4)=(\omr4){}\lom$  
forces, over $\rL$, that 
\bce
$[T,d]\sq [T_0,d_0]$ and $[T,d]$ is an 
\dd\uE equivalence class.
\ece 
%where $\bg$ is a canonical name for the 
%\dd{\coll(\omr4)}generic map $g:\om\na\omr4$.  

\ble
\lam{stL4}
If\/ $\ang{T,d}\in H$ then it is true in 
%any generic 
%extension\/ $\rL[G]$ or\/ $\rV[G]$, in particular, 
in the ground set universe\/ $\rV$ 
%itself, 
that\/ $[T,d]\sq [T_0,d_0]$ and\/ 
$[T,d]$ is a\/ 
%\dd\ELG class, resp., a\/ 
\de class.
\ele
\bpf
By definition this is true for \dd{\coll(\omr4)}generic
extensions of $\rL$ --- hence by Shoenfield also for 
all generic extensions $\rV[G]$ in which $\omr4$ is countable, 
and then, by quite obvious downward absoluteness, 
for $\rV$.
\epf

\ble
\lam{stL3}             
$H\ne\pu$.
\ele
\bpf
By Lemma~\ref{stL3} there is a stable condition 
$\ang{T',d'}\in{\dpl}$.
%Let $\mu<\nu$ be defined as in the proof of Lemma~\ref{stL2}.
%
Using an \dd{\omr4}enumeration of all dense sets $D\sq\dpl$ 
in $\rL$, we easily get a code 
$\ang{T^\ast,d^\ast}\in\bc$ 
such that $\sup T^\ast\le\omr4$ and 
$$
[T^\ast,d^\ast]=\ens{x\in[T',d']}
{x\text{ is \dd{\dpl}generic over }\rL}
$$
in any class $\rV[G]$. 
Lemma~\ref{stL1} implies that 
all elements $x\in[T^\ast,d^\ast]$ in $\rV[G]$ 
are \dd\EVG equivalent to each other and to some 
$y^\ast\in [T_0,d_0]\cap\rV$. 

Let $g:\om\na\omr4$ be a collapse-generic map.
\rit{We argue in\/ $\rV[g]$}. 
By a simple cardinality argument, $[T^\ast,d^\ast]\ne\pu$
in\/ $\rV[g]$, and $[T^\ast,d^\ast]$ consists of pairwise 
\dd\EVg equivalent elements by the above.
This allows us to define  
$$
Z =\ens{z}{\sus x\in [T^\ast,d^\ast]\:(x\EVg z)}=
\ens{z}{\kaz x\in [T^\ast,d^\ast]\:(x\EVg z)}\,.
$$
in $\rV[g]$, so that it is true in $\rV[g]$ that $Z$ is 
an entire \dd\EVg equivalence class, which includes 
$[T^\ast,d^\ast]$, hence, has a non-empty intersection 
with $[T',d']\sq[T_0,d_0]$, therefore $Z\sq[T_0,d_0]$ 
as $[T_0,d_0]$ is an \ds\dd{\EVg}class in 
$\rV[g]$ by \ref{resF}.

It follows that $Z$ is $\fp0{1+\ro}$ in $\rV[g]$. 
Moreover, by the choice of $g$ it is true in $\rV[g]$ that 
$\ang{T^\ast,d^\ast}\in\rL\cap\hc$, and hence 
$\ang{T^\ast,d^\ast}$ is $\idh1(\et)$ in $\rV[g]$ for an 
ordinal $\et<\omi^{\rV[g]}$.
(Indeed let $\et$ be the first ordinal such that 
$\ang{T^\ast,d^\ast}$ is the \dd\et th set in the G\"odel 
construction of $\rL$.) 
Then $Z$ is $\idh1(\et)$ in $\rV[g]$. 
Therefore by Proposition~\ref{sepp} that there is a code 
$\ang{T,d}\in\bc_\ro$ such that $Z=[T,d]$ in $\rV[g]$. 
Let us demonstrate that $\ang{T,d}\in H$.

Consider a collapse-generic map $g':\om\na\omr4$; we can assume 
that $g'$ is \dd{\coll(\omr4)}generic even over $\rV[g]$.
We have to prove that 
\ben
\fenu
\atc\atc
\itla{g'1}
in $\rL[g']$: $[T,d]\sq [T_0,d_0]$ and $[T,d]$ is an 
\dd\ELgp equivalence class.
\een
Recall that by construction 
$Z=[T,d]\sq [T_0,d_0]$ and $[T,d]$ is an 
\dd\EVg class in $\rV[g]$.
But the Borel codes involved are countable in both classes 
$\rV[g]$ and $\rL[g']$. 
This implies \ref{g'1} by Shoenfield. 
\epf

Now we have gathered everything necessary to end the proof 
of the theorem in a few lines. 
It suffices to prove that 
$C=[T_0,d_0]\sq\bigcup_{\ang{T,d}\in H}[T,d]$ in $\rV$. 
Suppose tovards the contrary that this is not the case.

The set $H\sq\bc_\ro$ belongs to $\rL$ and 
$\car H\le\omr1$ in $\rL$, of course. 
As $\ang{T_0,d_0}\in\bc_\rod$, we can easily define a code 
$\ang{T_1,d_1}\in\bc_\rod$ such that   
$[T_1,d_1]=[T_0,d_0]\bez\bigcup_{\ang{T,d}\in H}[T,d]$ 
in any universe, 
and hence $[T_1,d_1]\ne\pu$ in $\rV$ by the contrary 
assumption, and still 
$[T_1,d_1]$ is a \ds\de class in $\rV$ since so is 
$C=[T_0,d_0]$ while each $[T,d]\yt \ang{T,d}\in H$, 
is a \de class by Lemma~\ref{stL4}.

In other words, the code $\ang{T_1,d_1}\lec\ang{T_0,d_0}$ 
has the same properties 
(see Corollary~\ref{resG}) 
as $\ang{T_0,d_0}$ does.
In fact by exactly the same arguments as above, but with 
the forcing notion 
$\dpd=\ens{\ang{T,d}\in\dpl}{\ang{T,d}\lec\ang{T_1,d_1}}$,
the set $H_1$ of all codes $\ang{T,d}\in H$ such 
that the $\coll(\omr4)$ forces \lap{$[T,d]\sq [T_0,d_0]$} 
over $\rL$, is non-empty. 
(Compare Lemma~\ref{stL3}.) 
Let $\ang{T,d}\in H_1$.
Then $\pu\ne [T,d]\sq [T_1,d_1]$ in $\rV$ 
(similarly to Lemma~\ref{stL4}), 
which contradicts to the definition of $\ang{T_1,d_1}$.\vom

\epF{Theorem~\ref{resT}}

\qeDD{Theorem~\ref{tk}, full version}

\vyk{
the following sentences 
\ben
\aenu
%\atc
\itla{stL31}\msur 
$
[T^\ast,d^\ast]=\ens{x\in[T,d]}
{x\text{ is \dd{\dpl}generic over }\rL}
$, and

\itla{stL32}
all elements $x\in[T^\ast,d^\ast]$ are \dd\EVg equivalent, 
resp., \dd\ELg equivalent  
to each other and to some $y^\ast\in C$ (so $y^\ast\in\rV$) 
by Lemma~\ref{stL1},
\een
%\ref{stL31} and \ref{stL32} 
hold in any extension $\rL[G]$ or $\rV[G]$, and hence 
\ben
\aenu 
\atc\atc
\itla{stL33}\msur 
$[T^\ast,d^\ast]\ne\pu$ in any extension $\rL[G]$ 
of $\rL$ with $\omr4<\omi^{\rL[G]}$.
\een
Let $g:\om\na\omr4$ be a collapse-generic map.
\rit{Arguing in\/ $\rV[g]$}, we define  
$$
Z =\ens{z}{\sus x\in [T^\ast,d^\ast]\:(x\EVg z)}=
\ens{z}{\kaz x\in [T^\ast,d^\ast]\:(x\EVg z)}\,.
$$
Then $Z$ is a \dd\EVg class in $\rV[g]$, and $y^\ast\in Z$, 
hence, $Z\sq [T_0,d_0]$ in $\rV[g]$.
It follows that $Z$ is $\fp0{1+\ro}$ in $\rV[g]$. 
Moreover, by the choice of $g$ it is true in $\rV[g]$ that 
$\ang{T^\ast,d^\ast}\in\rL\cap\hc$, and hence 
$\ang{T^\ast,d^\ast}$ is $\idh1(\et)$ in $\rV[g]$ for an 
ordinal $\et<\omi^{\rV[g]}$.
(Indeed let $\et$ be the first ordinal such that 
$\ang{T^\ast,d^\ast}$ is the \dd\et th set in the G\"odel 
construction of $\rL$.) 
Then $Z$ is $\idh1(\et)$ in $\rV[g]$. 
It follows by Proposition~\ref{sepp} that there is a code 
$\ang{T,d}\in\bc_\ro$ such that $Z=[T,d]$ in $\rV[g]$. 
Let us demonstrate that $\ang{T,d}\in H$.
}

\vyk{
\punk{A theorem on $\is12$ sets}
\las{ts12}

\bte
\lam{mq}
Suppose that\/ $A\sq\bn$ is a\/ $\is12$ set.

\ete
\bpf
\vyk{
If $\om\le\la\in\Ord$ then there exists a bijection 
$\pib\la:\la\ti\la\na\la$. 
Moreover we can choose $\pib\la$ for all 
$\om\le\la\le\omi$ so that\vom
 
1) the sequence $\sis{\pib\la}{\om\le\la\le\omi}$ is 
$\id\HC1$, and \vom 
 
2) for $\la=\om$, if $m=2^i(2j+1)-1\in\dN$ then 
$\pib\om(i,j)=m$.\vom

\noi
Now f
}%
Fix an ordinal $\la$, $\om\le\la<\omi$. 
\vyk{
We write $\pai\xi\eta$ instead of $\pib\la(\xi,\et)$, 
and if $\za=\pib\la(\xi,\et)<\la$ then we put 
$\lev\za=\xi$ and $\pra\za=\et$.
Further, if $f\in\lao$ and $\xi
}%
Consider a language containing \rit{variables} 
of type $\la$ (with domain $\la$), 
usually denoted by $\xi,\eta$
(and also $k,m,n$ if explicitly bounded by $\om$, see 
Remark~\ref{bom}), 
and of type $\lao$ (with domain $\lao$), 
denoted by $f,g,h$. 
\rit{Terms} of type $\la$ are allowed, so that 
\ben
\tenu{(t\arabic{enumi})}
\itla{t1} 
variables of type $\la$ are terms, 

%\itla{t1+} 
%if $s,t$ are term then so are $\pai st$, $\lev t$, $\pra t$; 

\itla{t2} 
if $t$ is a term then $\an t$ is a term interpreted 
so that if $t<\dN$ then $\an t=t$, and otherwise 
$\an t=0$; 

\itla{t3} 
if $t$ is a term and $f$ is a 
variable of type $\lao$ then $f[t]$ is a term, 
interpreted as $f(\an t)$;

\itla{t4} 
if $t_1,\dots,t_m$ are terms and $r:\dN^m\to\dN$ a 
recursive function then $r[t_1,\dots,t_m]$ is a term, 
interpreted as $r(\an{t_1},\dots,\an{t_n})$;

\itla{t5} 
in particular, if $s,t$ are terms then 
$t+s$, $t\cdot s$, $t^s$ are terms, 
interpreted as, resp., 
${\an t}+{\an s}$, ${\an t}\cdot{\an s}$, ${\an t}^{\an s}$ 
(with $0^0:=0$),
while $\lh \ns t$ and $\bin ts$ are terms, interpreted 
as, resp., 
$\lh\ns{\an t}$ and $\bin{\an t}{\an s}$, 
where: 

$\sis{\ns k}{k<\dN}$ is a fixed recursive enumeration of 
all strings in $\nse$, 

$\lh u$ is the length of a string $u\in\nse$,

and $\bin mi=\ns m(i)$ provided $i<\lh s$, but $\bin mi=0$ 
otherwise.
\een 

The relational symbols are: 
\ben
\tenu{(r\arabic{enumi})}
\itla{r1} 
the equality of terms $t_1=t_2$; 

\itla{r2} 
the ordering of terms $t_1<t_2$ 
(interpreted as the natural order on $\la$); 

\itla{r3} 
the unary predicate $\bom(t)$ interpreted as $t\in\dN$; 

\itla{r4} 
all recursive relations $R[t_1,\dots,t_m]$ of any arity $m$, 
interpreted as $R(\an{t_1},\dots,\an{t_n})$.
\een

If $m=2^i(2j+1)-1<\om$ then let $\lev m=i$, 
$\pra m=j$, and $m=\pai ij$.
If $f\in\lao$ and $i<\om$ then let $\bin fi\in\lao$ 
be defined so that $\bin fi(j)=f(\pai ij)$ for all $j<\om$.

\bre
\lam{term2}
We allow to use terms of type $\lao$ of the form $\bin ft$, 
where $t$ is a term of type $\la$, as above --- as shorthands. 
That is, any occurrence of, say, $\bin ft(s)$ is just 
a shorthand  for $f(\pai{\an t}{\an s})$.
\ere

An formula of this language is 
\bde
\item[\rit{elementary,}] 
if it is  
made of relational symbols (properly filled in with terms) 
by means of propositional connectives and 
\rit{bounded quantifiers}, \ie, those of the form 
$\sus \xi<\an t$ and $\kaz \xi<\an t$;\snos
{The domain of a bounded quantifier is 
a finite set, an initial segment of $\dN$.
}%

\item[\rit{arithmetic,}] 
if in addition quantifiers $\sus\xi<\om$ and 
$\kaz\xi<\om$ are allowed;

\item[\rit{analytic,}] 
if in addition to the above quantifiers over $\lao$ 
are allowed.
\ede

\bdf
\lam{rla}
By definition 
{\ubf we do not allow}
 quantifiers over the full domain 
$\la=\ens{\xi}{\xi<\la}$.
\edf

\bre
\lam{bom}
We'll use $k,m,n$ for variables of type $\la$ 
explicitly bounded by $\om$, that is, $\sus k\:\Phi(k)$ 
is a shorthand for $\sus\xi<\om\:\Phi(\xi)$.
\ere

Define classes $\is1n$ and $\ip1n$ of analytic formulas 
usual way. 
For instance a $\is13$ formula is any formula of the form 
$$
\sus f\:\kaz g\:\sus h\:\kaz\xi<\om\:\Phi(f,g,h,\xi),
$$
where $\Phi$ is an elementary formula. 
Accordingly, classes $\isa1n$ or $\ipa1n$ contain all subsets 
of spaces of the form $\la^m\ti(\lao){}^k$, definable by a 
$\is1n$ formula, resp., $\ip1n$ formula, with arbitrary 
parameters $\xi<\la$ 
(not only natural numbers!) 
allowed.

The restriction of Definition~\ref{rla}  
allows us to prove all 
ordinary quantifier transformation rules. 
For instance, 
$$
\kaz\xi<\om\:\sus f\:\Phi(\xi,f)
\leqv
\sus f\:\kaz\xi<\om\:\Phi(\xi,\bin f\xi).
$$ 
It follows that classes $\isa1n$ or $\ipa1n$ satisfy the 
same basic closure properties as classes $\is1n$ or $\ip1n$
in the usual sense do.
There are some differences, of course. 
For instance, there is no universal $\isa1n$ 
sets $U\sq\om\ti \lao$.  

\epf
}

\punk{OD sets in Solovay's model: generalizations}
\las{sm1}

It is a rather common practice that results like Theorems 
\ref{mt'} and \ref{mt} generalize this or another way 
in the Solovay model.\snos 
{By \rit{the Solovay model} we'll always mean a model of 
$\ZFC$, 
a generic extension of $\rL$ introduced in \cite{sol}, 
in which all projective sets 
are Lebesgue measurable, rather than the other 
model of \cite{sol}, 
in which only $\ZF+\DC$ holds but all sets 
of reals are measurable.} 
We'll prove below the following such generalizations 
of our main results.
Recall that $\od$ means ordinal-definable, and also denotes 
the class of all ordinal-definable sets.

\bte
\lam{ms'}
The following is true in the Solovay model. 
If\/ $A\sq\bn$ is an\/ $\od$ set then one and only one of 
the next two claims\/ \ref{ms'1}, \ref{ms'2} holds$:$
\ben
\tenu{{\rm(\Roman{enumi})}}
\tenu{{\rm(\Roman{enumi})}}
\itla{ms'1}\msur
$A$ is\/ \ddd{\od}effectively \bou, so that 
there exists an\/ $\od$ 
sequence\/ $\sis{T_\xi}{\xi<\omi^\rL}$ of compact trees\/ 
$T_\xi\sq\nse$ such that\/ $A\sq\bigcup_\xi[T_\xi]\,;$ 
%--- therefore\/ 
%$A$ is covered by the union\/ $U$ of all sets of the form\/ 
%$[T]$, where\/ $T\sq\nse$ is a compact\/ $\od$ tree$;$ 

%\itla{ms'1}\msur
%$A$ is covered by the union\/ $U$ of all sets of the form\/ 
%$[T]$, where\/ $T\sq\nse$ is a compact\/ $\od$ tree$;$ 
%%{\rm(then $U$ is \ddd\fsg compact)} 

\itla{ms'2}
there is a superperfect set\/ $Y\sq A$.  
\een
\ete

Here conditions \ref{ms'1} and \ref{ms'2} are incompatible. 
Indeed, the union $U$ in \ref{ms'1} is countable, hence the 
set $U$ is \ddd\fsg compact. 
Therefore if $Y$ is a set as is \ref{mt'2} then $Y$ cannot 
be covered by $U$. 

Note that condition \ref{ms'1} cannot be strengthened 
to the form that 
\rit{there is an\/ $\od$ 
sequence\/ $\sis{T_n}{n\in\dN}$ of compact trees\/ 
$T_n\sq\nse$ such that\/ $A\sq\bigcup_n[T_n]$}. 
For a counterexample take $A=\bn\cap\rL$ 
(all constructible reals). 
This is a countable set in the Solovay model, 
hence \ref{ms'1} of Theorem~\ref{ms'} holds and \ref{ms'2} fails, 
but the existence of an $\od$ 
(hence, constructible) 
sequence of trees as indicated is clearly impossible.   

\bte
\lam{ms}
The following is true in the Solovay model. 
If\/ $A\sq\bn$ is an\/ $\od$ set then one and only one of 
the next two claims\/ \ref{ms1}, \ref{ms2} holds$:$
\ben
\tenu{{\rm(\Roman{enumi})}}
\itla{ms1}\msur
$A$ is\/ \ddd{\od}effectively \sik, so that 
there exists an\/ $\od$ sequence\/ 
$\sis{T_\xi}{\xi<\omi^\rL}$ of compact trees\/ 
$T_\xi\sq\nse$ such that\/ $A=\bigcup_\xi[T_\xi]\,;$ 
%--- therefore\/
%$A$ is equal to the union\/ $U$ of all sets of the form\/ $[T]$, 
%where\/ $T\sq\nse$ is a compact\/ $\od$ tree and\/ $[T]\sq A\:;$

%\itla{ms1}\msur
%$A$ is equal to the union\/ $U$ of all sets of the form\/ $[T]$, 
%where\/ $T\sq\nse$ is a compact\/ $\od$ tree and\/ $[T]\sq A\,;$
%{\rm(and then $A$ is \ddd\fsg compact)}$;$ 
%--- and moreover there is an\/ $\od$ 
%sequence\/ $\sis{T_n}{n\in\dN}$ of compact trees\/ 
%$T_n\sq\nse$ such that\/ $A=\bigcup_n[T_n]\,;$  

\itla{ms2}
there is a set\/ $Y\sq A$ homeomorphic to\/ $\bn$ and relatively 
closed in\/ $A$.
\een
%In addition, in case\/ \ref{mt1}, if\/ $A\ne\pu$ then\/ $A$ 
%contains a\/ \ddd{\id11}element.
\ete

The proof of both theorems follows in the next section.

\vyk{
is similar to the proof of 
Theores~\ref{ms'} and \ref{ms}, with the 
major difference being that the Gandy -- Harrington topology 
is replaced by the forcing by $\od$ (ordinal-definable) sets 
in the Solovay model.
}

Let's start with some definitions and a couple of special 
results related to the Solovay model. 
If $\Om$ is an ordinal then by \osm\ we denote 
the following sentence:
\lap{$\Om=\omi$, 
$\Om$ is strongly inaccessible in $\rL$, the 
constructible universe,   and
the whole universe $\rV$ is 
a generic extension of $\rL$ via a known collapse forcing 
$\text{\rm Coll}(\om,{}<\Om)$, as in \cite{sol}}. 
Thus \osm\ says that the universe is a Solovay-type 
extension of $\rL$.

\ble
[assuming \osm]
\lam{oO}
If\/ $X$ is a countable\/ $\od$ set then there exist an 
ordinal\/ $\la<\Om$ and an\/ $\od$ $1-1$ map\/ $f:\la\onto X$.
\ele
\bpf
Let $F:\Ord\onto\od$ be a canonical $\od$ map. 
Recall that under \osm\ the universe is a homogeneous generic 
extension of $\rL$. 
Therefore the relations $F(\xi)\in X$ and $F(\xi)=F(\eta)$ 
(with arguments $\xi,\eta$)
are $\od$.
%The rest of the proof is obvious.
\epf

\bdf
[assuming \osm]
\lam{dod}
Let $\OD$ be the collection of all {\ubf non-empty} 
$\od$ sets $Y\sq \bn$.
We consider $\OD$ as a forcing notion 
(smaller sets are stronges conditions).
A set $G\sq \OD$ is \rit{\ddd{\OD}generic over $\od$} 
if it non-emptily intersects every $\od$ dense set $D\sq\OD$.
\edf

\vyk{
A set $W\sq \OD$ is:
\bde
\item[\psur\quad\rit{dense}\rm,] 
iff 
for every $Y\in\OD$ there exists $Z\in W$, $Z\sq Y$;

\item[\quad\rit{\ddd{\OD}generic}\rm,] iff 
1) 
if $X,Y\in W$ then there is a set $Z\in W$, $Z\sq X\cap Y$, 
and 
2) 
if $D\sq\OD$ is $\od$ and dense then $W\cap D\ne\pu$.
\qed 
\ede
}   

\bpro
[see, \eg, \cite{ksol}] 
\lam{*gen}
Assuming \osm, 
if a set\/ $G\sq\OD$ is\/ \ddd\OD generic then the 
intersection\/ $\bigcap G=\ans{a_G}$ consists of a 
single real.\qed
\epro

%Unlike the set of all $\is11$ sets, 
As the set $\OD$ is definitely uncountable, 
the existence of \ddd\OD generic sets does not 
immediately follow from \osm\ by a cardinality argument. 
Yet fortunately $\OD$ is \rit{locally countable}, in a sense. 

\bdf
[assuming \osm]
\lam{odik}
A set $X\in\od$ is \rit{\odk} if the $\od$ power set 
$\pwod X = \pws X\cap \od$ is at most countable.\snos
{Then the set $\pwod X$ is 
%countable and \od, but 
not necessarily \ddd\od countable. 
Take for instance $X=\dN$.}

Let $\odi$ be the set of all \odk\ sets $X\in\OD$.
\edf

For instance, assuming \osm, the set 
$X=\bn\cap\od=\bn\cap\rL$ 
of all $\od$ reals belongs to $\odi$. 
Indeed  
$\pwod X = \pws X\cap \od= \pws X\cap \rL$, and hence 
$\pwod X$ admits an $\od$  
bijection onto the ordinal $\om_2^{\rL}<\omi$.

The set $\koh$ of all reals $x\in\bn$ 
Cohen generic over $\rL$ belongs to $\odi$ as well. 
Indeed if $Y\sq\koh$ is \od\ and $x\in Y$ then ``$x\in Y$'' 
is Cohen-forced 
over $\rL$. 
It follows that there is a set $S\sq\nse\yt S\in\rL$, such 
that $Y=X\cap\bigcup_{t\in S}\ibn t$.
But the collection of all such sets $S$ belongs to $\rL$ 
and has cardinality 
$\omi$ in $\rL$, hence, is countable under \osm.

%Note that if $X\in\od$ then the set $\pwod X$ belongs to 
%$\od$ either. 
%Therefore if 

\bpro
\lam{*den}
Assuming \osm, $\odi$ is dense in\/ $\OD$, that is, 
if\/ $X\in\OD$ then there is a set\/ $Y\in\odi$ 
such that\/ $Y\sq X$.\qed
\epro
\bpf
[sketch, see details in \cite{ksol}] 
For any ordinal $\la$, let $\koh_\la$ be the set of all 
elements $f\in\la^\om$, \ddd{(\la^{<\om})}generic over $\rL$. 
Suppose that $X\in\OD$. 
Then by definition $X\ne\pu$, hence, there is 
a real $x\in X$. 
Then it follows from \osm\ that there exist: 
an ordinal $\la<\omi=\Om$, 
an element $f\in\koh_\la$, 
and an $\od$ map $H:\la^\om\to\bn$, 
such that $x=H(f)$. 
The set $P=\ens{f'\in\koh_\la}{H(f')\in X}$ is then \od\ 
and non-empty (contains $f$), and hence so is its image 
$Y=\ens{H(f')}{f'\in P}\sq X$ (contains $x$). 

It remains to prove that $Y\in\odi$. 
As $H$ is an \od\ map, it is sufficient to show that 
$\koh_\la$ is $\odi$. 
But this is true by the same reasons as for the set $\koh$ 
(see just before Proposition~\ref{*den}). 
\epf

\bre
\lam{fg}
One may want to know whether Theorem~\ref{fm} also 
admits a 
version similar to Theorem~\ref{ms'} --- that is, for 
a finite sequence of $\od$ \er s $\rF_j$ and an $\od$ 
set $A$ in the Solovay model.

But here we have a grave obstacle just   
from the beginning. 
Indeed, coming back to the derivation of \ref{pssC} from 
Theorem~\ref{pss}, we'll have to prove that, 
in the Solovay model, 
any $\od$ set $E\sq\bn\ti\bn$ with \bou\ sections splits 
into a countable union of $\od$ sets with bounded 
sections. 
But this claim fails even for sets with  
\rit{countable} sections: consider \eg\ the $\is11$ set 
$E=\ens{\ang{x,y}}{y\in\rL[x]}$. 

Whether a more modest version holds in the Solovay model, 
with still $\id11$ relations $\rF_j$ and an $\od$ set $A$, 
remains to be seen.
\ere

\punk{OD sets in Solovay's model: proofs}
\las{sm2}

Here we prove Theorems \ref{ms'} and \ref{ms}.
The proofs strongly resemble those in Section~\ref{d2} 
and Section~\ref{gahaSK}, hence we skip some details. 
There are two notable differences. 
First, the \gh\ type of arguments is replaced by the $\od$ 
forcing, and second, various niceties related to classes 
$\is11$ and $\id11$ become obsolete as $\od$ is a more robust 
definability class.

\bpf
[Theorem~\ref{ms'}]
\rit{We argue in the Solovay model}, that is, we assume 
\osm. 
Consider an arbitrary $\od$ set $A\sq\bn$. 
Let $U$ be the union of all sets of the form\/ 
$[T]$, where\/ $T\sq\nse$ is a compact\/ $\od$ tree. 
Clearly the set $U$ and the difference $A'=A\dif U$ are $\od$.

\ble
\lam{skm*}
Under the conditions of Theorem~\ref{ms'}, 
if\/ $Y\sq A'$ is a non-empty\/ $\od$ set then 
its topological closure\/ $\clo Y$ in\/ $\bn$ 
is not compact. 
\ele
\bpf
If $\clo Y$ is compact  
then $T=\der Y$ is a compact $\od$ tree, hence  
$Y\sq\clo Y=[T]\sq U$, a contradiction 
to the assumption $Y\sq A'$.
\epF{Lemma}

\rit{Case 1}: $A'=\pu$, that is, $A\sq U$. 
%Thus \ref{mt'1a} of Theorem~\ref{ms'} holds. 
To check \ref{ms'1} of Theorem~\ref{ms'}, 
note that under \osm\ 
$\od$ reals are the same as constructible reals, and 
hence there is an $\od$ enumeration of all $\od$ trees 
by ordinals $\xi<\omi^\rL$.\vom 

\rit{Case 2}: the set
%\ref{mt'1} fails, \ie, 
$A'=A\dif U$ is non-empty. 
By Proposition~\ref{*den}, there is a set $A''\sq A'$,  
$A''\in\odi$.
Then the power set 
$P=\pwod{A''}= \pws{A''}\cap \od$ is at most countable. 
By Lemma~\ref{oO}, there exist an ordinal $\la<\Om$ 
and an $\od$ map $f:\la\onto P$.
But the power set $\pwod\la$ is obviously countable, 
therefore so is $\pwod P$. 
Fix an arbitrary enumeration $\sis{\cD^\od_n}{n\in\om}$ of 
all $\od$ sets $\cD\sq P=\pwod{A''}$, 
dense in $\odi$ below $A''$.
We assert that then there is a system of non-empty 
$\od$ sets $Y_s\sq A'$ satisfying conditions 
\ref{gan1}, \ref{gan2}, \ref{gan3}, \ref{han5} 
in Section~\ref{d2}, along with the following condition 
instead of \ref{gan4}:
\ben
\itsep
\tenu{$(\arabic{enumi}^{\text{\sc od}})$}
\atc\atc\atc 
\itla{san4} 
if $s\in\nse$ then $Y_s\in\cD^\od_{\lh s}$.
\een
If such a construction is accomplished then 
%\ref{san4} implies 
$\bigcap_mY_{a\res m}=\ans{f(a)}$ 
for each $a\in\bn$ by Proposition~\ref{*gen}, and 
$f:\bn\na Y=\ens{f(a)}{a\in\bn}$ is   
a homeomorphism. 
Moreover the set $Y$ is closed in $\bn$ by exactly the 
same reasons as in Section~\ref{d2},  
and hence we have \ref{mt'2} of Theorem~\ref{ms'}. 

The construction of sets $Y_s$ goes on exactly as in 
Section~\ref{d2}, with the only difference that   
$\is11$ and Lemma~\ref{tkm*} 
are replaced by $\od$ and Lemma~\ref{skm*}.\vom

\epF{Theorem~\ref{ms'}}

\bpf
[Theorem~\ref{ms}]
Assuming \osm, consider any $\od$ set $A\sq\bn$. 
Let $U$ be the union of all sets $[T]$, 
where\/ $T\sq\nse$ is a compact\/ $\od$ tree and $[T]\sq A$. 
The set $U$ and the difference $A'=A\dif U$ are $\od$.

By Theorem~\ref{ms'}, we can \noo\ assume 
that $A$ is \bou, and hence 
if $F\sq A$ is a closed set then $F$ is \ddd\fsg compact.

\ble
\lam{skm-l}
If\/ $F\sq A'$ is a non-empty\/ $\od$ set then\/ 
$\clo F\not\sq A$. 
% is not\/ \ddd\fsg compact.
\ele

Recall that $\clo F$ is the closure of a set $F\sq\bn$.

\bpf
Suppose towards the contrary that 
$\pu\ne F\sq A'$ is an $\od$ set but $\clo F\sq A$. 
By the \noo\ assumption above,  
$\clo F=\bigcup_nF_n$ is \ddd\fsg compact, 
where all $F_n$ are compact. 
There is a Baire interval $\ibn s$ such that the set 
$X=\ibn s\cap \clo F$ is non-empty and $X\sq F_n$ for some $n$. 
Thus $X\sq A$ is a non-empty compact $\od$ set, hence 
by definition $X\sq U$ and $A'\cap X=\pu$. 
In other words, $\ibn s\cap \clo F\cap A'=\pu$.
It follows that $\ibn s\cap F=\pu$ (because $F\sq A'$), 
which contradicts to $X=\ibn s\cap \clo F\ne \pu$. 
\epF{Lemma}

We come back to the proof of Theorem~\ref{ms}.\vom

\rit{Case 1}: $A'=\pu$, that is, $A=U$. 
This implies \ref{ms1} of the theorem.\vom
      
\rit{Case 2}: 
$A'\ne\pu$. 
As in the proof of Theorem~\ref{ms'}, 
choose a set $A''\sq A'$, $A''\in\odi$, and 
fix an arbitrary enumeration $\sis{\cD^\od_n}{n\in\om}$ of 
all $\od$ sets $\cD\sq P=\pwod{A''}$, 
dense in $\odi$ below $A''$.
To get a set $Y\sq A''$, \rit{relatively} 
closed in $A$ and homeomorphic to $\bn$, 
%as in \ref{mt2} of the theorem, 
we make use of  
a system of non-empty $\od$ sets $Y_s\sq A''$ satisfying 
conditions \ref{gan1}, \ref{gan2}, \ref{gan3}  
in Section~\ref{d2}, 
\ref{san4} as in the proof of Theorem~\ref{ms'}, 
and \ref{gan5} in Section~\ref{gahaSK}.

If such a system of sets is defined , then the associated 
map $f:\bn\to A''$ is $1-1$ and is a homeomorphism 
from $\bn$ onto its full image 
$Y=\ran f =\ens{f(a)}{a\in\bn}\sq A''$. 
In addition, the set $Y$ is relatively closed in $A$ 
by the same arguments 
(based on condition \ref{gan5}) 
as in Section~\ref{gahaSK},  
and hence we have \ref{ms2} of Theorem~\ref{ms}.  
The construction of sets $Y_s$ also goes on as in 
Section~\ref{gahaSK}, but we have to apply 
Lemma~\ref{skm-l} instead of Lemma~\ref{tkm-l}.\vom

\epF{Theorem~\ref{ms}}

\vyk{

\punk{Infinite-case generalization}
\las{!!i}

What about countably many \er s?
The following is a conjecture.

\bte
\lam{ifm}
Suppose that\/  
$\rF_0,\rF_1,\rF_2,\dots$ are\/ Borel \er s on\/ $\bn$.
Then one and only one of 
the following two claims holds$:$
\ben
\tenu{{\rm(\Roman{enumi})}}
\itla{ifm1}
the domain\/ $\bn$ is\/ \sm{\ans{\rF_0,\rF_1,\rF_2,\dots}}\/  
{\rm(so $\bn$ is a countable 
union of\/ \ddd{\rF_0}classes, \ddd{\rF_1}classes, 
\ddd{\rF_2}classes, \dots , plus a \bou\/ set)}$;$

\itla{ifm2}
there is a \sps\/ $P\sq \bn,$ a\/
\pis{\rF_i} 
%pairwise\/ \den i 
for each\/ $i<\om$.  
\een
\ete

\bpf Not yet known

\epf

\vyk{

\bpf[tentative]
Assume that the whole sequence $\sis{\rF_n}{n<\om}$ 
is $\id11$.\vom

\rit{Case 1}: 
the set 
$\bU=\ens{x\in\bn}{\fn x 0\,
\text{ is \lar{\ans{\rF_1,\rF_2,\dots}}}}$ 
has at most countably many \de classes. 
Then $\bU$ and 
$$
\bC=\bn\bez\bU=
\ens{x\in\bn}{\fn x 0\,
\text{ is \sm{\ans{\rF_1,\rF_2,\dots}}}}
$$ 
are Borel sets. 
If the set $\bC$ is \sm{\ans{\rF_0,\rF_1,\dots}} then 
we immediately have \ref{fm1}, so suppose that 
$\bC$ is \lar{\ans{\rF_0,\rF_1,\dots}}. 
Our goal will be to find a \spp\ \pis{\rF_0} 
$X\sq\bC$ which is still \lar{\ans{\rF_1,\rF_2,\dots}}.

In our assumptions, the set 
$\bC$ is Borel, and as usual we suppose, 
 for the sake of brevity, that $\bC$ is $\id11$.
%
%If $x\in\bC$ then the \ddd{\rF_0}class $\fn x0\sq\bC$ is 
%\sm{\ans{\rF_1,\rF_2,\dots}}. 
Suppose that $x\in\bC$. 
Then the \ddf0class $\fn x 0$ is an 
\sm{\ans{\rF_1,\rF_2,\dots}} $\id11(x)$ set. 
Therefore by Lemma~\ref{ico} 
$\fn x 0$ is covered by the union of all $\id11(x)$ sets 
each of which is either a compact set or a \ddf nclass 
for some $n\ge1$. 

Now let $R$ be the set of all triples $\ang{x,y,m}$ such that 
$x\yi y\in\bC$, $x\rF_0 y$, $m$ is a code of a $\id11(x)$ 
set $B_m$, $y\in B_m$, 
and $B_m$ is either compact or an \ddf nclass 
for some $n\ge1$. 
Then $R$ is $\ip11$, 
and the projection
$$
\dom R= \ens{\ang{x,y}}{\sus n\:(\ang{x,y,m}\in R)}
$$
coincides with the set 
$\bcd=
%\cfo\bC{\rF_0}=
\ens{\ang{x,y}}{x,y\in \bC\land x\rF_0 y}$ 
by the remark just above. 
By $\ip11$ uniformization, there is a $\id11(x)$ map 
$F:\bcd\to\om$ such that 
$\ang{x,y,F(x,y)}\in R$ for all 
$\ang{x,y}\in \bcd$.
Then by definition if $\ang{x,y}\in \bcd$ then 
$m=F(x,y)$ is a code of a $\id11(x)$ set 
$\ind xy=B_m$, $y\in \ind xy$, 
and $\ind xy$ is either compact or an \ddf nclass for some $n\ge1$.

We then   
\bit 
\item
let $\cor 1x$ be the set of all\/ $\id11(x)$ 
\pfs 1s $X\sq \fn x0$ 
(which we take instead of \ddf 1classes) and put 
$D_1(x)=\bigcup_{Y\in\cor1x}Y$, 

\item
let $\cor 2x$ be the set of all\/ $\id11(x)$ 
\pfs 2s $X\sq \fn x0\bez D_1(x)$ and put 
$D_2(x)=\bigcup_{Y\in\cor2x}Y$, 

\item
let $\cor 3x$ be the set of all\/ $\id11(x)$ 
\pfs 2s $X\sq \fn x0\bez {(D_1(x)\cup D_2(x))}$ and put 
$D_3(x)=\bigcup_{Y\in\cor3x}Y$, 

\item
and so on. 
\eit
Accordingly to an obvious generalization 
of Lemma~\ref{2sil}, $\fn x0$ is covered by the    
union of all\/ $\id11(x)$ \pfs ns $X\sq \fn x0$ 
for various $n\ge1$, 
and the union $K_x$ of all\/ $\id11(x)$ compact sets 
--- thus $\fn x0\sq K_x\cup \bigcup_{n\ge1}D_n(x)$. 

If $y\in \fn x0$ then $D_y$ and $K_y$ may be different 
from resp.\ $D_x$ and $K_x$, yet the sets $D_x$ with $x$ 
in one and the same \ddf0class have a sufficient 
invariant core!
Therefore, accordingly to an obvious generalization 
of Lemma~\ref{2sil}, $\fn x0$ is covered by the    
union $D_x$ of all\/ $\id11(x)$ \pfs ns $X\sq \fn x0$ 
(which we take instead of \ddf nclasses) 
for various $n\ge1$, 
and the union $K_x$ of all\/ $\id11(x)$ compact sets. 
If $y\in \fn x0$ then $D_y$ and $K_y$ may be different 
from resp.\ $D_x$ and $K_x$, yet the sets $D_x$ with $x$ 
in one and the same \ddf0class have a sufficient 
invariant core!

\bcl
\lam{!!1c1}
If\/ $x\rF_0 y$, $n\ge1$, and\/ $Y\sq \fn x0$ 
is a\/ non-\bou\/ $\id11(x)$ \pfs n\/ 
{\rm(so that\/ $Y\sq D_x$)} 
then\/ $Y$ is\/ $\id11(y)$ as well.
\ecl
\bpf

\epF{Claim}

Let $\cor nx$ be the collection of 
all \pfs ns $Y\sq \fn x0$ which belong to $\id11(y)$ for 
\rit{every} $y\in \fn x0$.
Then $\cor nx$ is \ddf0invariant, we have 
$\bigcup\cor nx\sq D_x$, 

\epf
}
}

\vyk{

\punk{*}
\las{*}

\ble
[assuming \osm]
\lam{icoS}
Suppose that\/  $\sis{\rF_n}{n<\om}$ 
is an\/ $\od$ sequence of\/ $\id11$ \er s on\/ $\bn$, 
and an\/ $\od$ set\/ $X\sq\bn$ is\/ \sm{\sis{\rF_n}{n<\om}}.
Then\/ $X$ is\/ 
\ddd{\od}effectively\/ \sm{\sis{\rF_n}{n<\om}}, in 
the sense that it is covered by the union of all\/ $\od$ 
compact sets and the union of all\/ $\od$ \ec es of the 
relations\/ $\rF_j$.
\ele

\bpf
The set $C=\ct\cap\od$ of all $\od$ compact trees
is $\od$, and hence so is $K=\bigcup_{T\in C}\bod{T}$. 
If $n<\om$ then let $U_n$ be the union of all $\od$ 
\ddf nclasses. 
The set $U={\bigcup_nU_n}$ is $\od$ either.\vom

\rit{Case 1}: 
$X\sq K\cup U$.
Then\/ $X$ is\/ 
\ddd{\od}effectively\/ \sm{\sis{\rF_n}{n<\om}}.\vom

\rit{Case 2}:  
$A=X\bez{(K\cup U)}\ne \pu$. 
Then $A$ is a non-empty $\od$ set. 
We are going to derive a contradiction.
By definition, we have 
$X\sq \bigcup_kC_k\cup \bigcup_n\bigcup_k E_{nk}$, 
where each $C_k$ is compact and each $E_{nk}$ is an 
\ddf nclass.
Let $M$ be a countable elementary substructure of a 
sufficiently large structure, containing, in particular, 
the whole sequence of covering sets $C_k$ and $E_{nk}$. 
Below ``generic'' will mean \gh\ generic over $M$.

As $A\ne\pu$ is $\is11$, there is a perfect set $P\sq A$ 
of points both generic and pairwise generic. 
It is known that then $P$ is a \pis{\rF_n} for every $n$, 
hence, definitely a set not covered by a countable union 
of \ddf nclasses for all $n<\om$. 
Thus to get a contradiction it suffices to prove that 
$P\cap C_k=\pu$ for all $k$. 
In other words, we have to prove that if $k<\om$ and 
$x\in A$ is any generic real then $x\nin C_k$. 

Suppose towards the contrary that a non-empty $\is11$ 
condition $Y\sq A$ forces that $\ja\in C_k$, where 
$\ja$ is a canonical name for the \gh\ generic real. 
We claim that $Y$ is not \bou. 
Indeed otherwise we have $Y\sq \bigcup_n \bod{T_n}$ by 
Theorem~\ref{mt'}, where all trees $T_n\sq\nse$ are 
$\id11$ and compact, which contradicts the fact that $A$ 
%(and hence $X$ as well) 
does not intersect any compact $\id11$ set. 

Therefore $Y\not\sq C_k$. 
Then there is a point $x\in Y$ and a number $m$ such 
that the set $I=\ens{y\in\bn}{y\res m=x\res m}$ 
does not intersect $C_k$. 
But then the $\is11$ condition $Y'=Y\cap I$ 
forces that $\ja\nin C_k$, a contradiction.
\epf

}

\end{document}